 \documentclass[numbers,webpdf,imaiai]{ima-authoring-template}%
\pdfoutput=1
\usepackage{subfigure}
\usepackage{multirow, booktabs}
\graphicspath{{Fig/}}
\usepackage[figuresright]{rotating}
\usepackage{mathrsfs}

\theoremstyle{thmstyletwo}%
\newtheorem{theorem}{Theorem}
\newtheorem{proposition}[theorem]{Proposition}%
\newtheorem{remark}{Remark}%
\newtheorem{lemma}{Lemma}%
\newtheorem{definition}{Definition}

\numberwithin{equation}{section}

\usepackage{enumitem}

\begin{document}

 \DOI{}
\copyrightyear{2023}
\firstpage{1}


\title[ADC-siDCA for SPLRMR]{An efficient asymptotic DC method for sparse and low-rank matrix recovery}

\author{Mingcai Ding
\address{\orgdiv{College of sciences}, \orgname{Shihezi University}, \orgaddress{\street{Xiangyang Street}, \postcode{832003}, \state{Xinjiang}, \country{People’s Republic of China}}}}
\author{Xiaoliang Song*
\address{\orgdiv{School of Mathematical Sciences}, \orgname{Dalian University of Technology}, \orgaddress{\street{Lingshui Street}, \postcode{116081}, \state{Liaoning}, \country{People’s Republic of China}}}}
\author{Bo Yu
\address{\orgdiv{ National Center for Applied Mathematics in Chongqing}, \orgname{Chongqing Normal University}, \orgaddress{\street{}, \postcode{401331}, \state{Chongqing}, \country{People’s Republic of China}}}
\address{\orgdiv{School of Mathematical Sciences}, \orgname{Dalian University of Technology}, \orgaddress{\street{Lingshui Street}, \postcode{116081}, \state{Liaoning}, \country{People’s Republic of China}}}
}


\corresp[*]{Corresponding author: \href{email:songxiaoliang@dlut.edu.cn}{songxiaoliang@dlut.edu.cn}}

\received{14}{11}{2023}


\abstract{The optimization problem of sparse and low-rank matrix recovery is considered, which involves a least squares problem with a rank constraint and a cardinality constraint. To overcome the challenges posed by these constraints, an asymptotic difference-of-convex (ADC) method that employs a Moreau smoothing approach and an exact penalty approach is proposed to transform this problem into a DC programming format gradually.  To solve the gained DC programming, by making full use of its DC structure, an efficient inexact DC algorithm with sieving strategy (siDCA) is introduced. The subproblem of siDCA is solved by an efficient dual-based semismooth Newton method. The convergence of the solution sequence generated by siDCA is proved. To illustrate the effectiveness of ADC-siDCA, matrix recovery experiments on nonnegative and positive semidefinite matrices. The numerical results are compared with those obtained using a successive DC approximation minimization method and a penalty proximal alternating linearized minimization approach. The outcome of the comparison indicates that ADC-siDCA surpasses the other two methods in terms of efficiency and recovery error. Additionally, numerical experiments on sparse phase retrieval demonstrate that ADC-siDCA is a valuable tool for recovering sparse and low-rank Hermitian matrices.}
\keywords{Matrix recovery, Sparse phase retrieval,
Rank constraint, Cardinality constraint, Inexact difference-of-convex algorithm, Semismooth Newton method}


\maketitle

\section{Introduction}\label{intro}
In various signal processing and machine learning applications, matrices often possess multiple underlying structures \cite{Ref_oymak2015simultaneously}. The jointly sparse and low-rank structured model has been effectively implemented in various domains such as sub-wavelength imaging \cite{Ref_shechtman2011sparsity}, hyperspectral image unmixing \cite{Ref_giampouras2016simultaneously,Ref_tsinos2017distributed}, feature coding \cite{Ref_Zhang2016Joint}, covariance matrix estimation \cite{Ref_richard2012estimation,Ref_chen2014estimation}, and sparse phase retrieval \cite{Ref_ohlsson2012compressive,Ref_ohlsson2012cprl,Ref_li2013sparse}. The sparse and low-rank matrix estimation problem  \cite{Ref_shechtman2011sparsity, Ref_giampouras2016simultaneously,Ref_tsinos2017distributed, Ref_richard2012estimation,Ref_chen2014estimation, Ref_ohlsson2012compressive,Ref_li2013sparse,Ref_chen2018simultaneously,Ref_liu2019successive} has garnered significant attention in recent years. This paper focuses on the sparse and low-rank matrix recovery problem that can be defined as follows:
\begin{equation}\label{eq1}
\begin{aligned}
  \min_{\mathbf{U}\in\mathbb{U}_+}&\; \ell(\mathbf{U})=\frac{1}{2}\|\mathcal{A}(\mathbf{U})-\boldsymbol{b}\|^{2}\\
  \text{s.t.} &\;\operatorname{rank}(\mathbf{U})\leq r,\\
  &\;\|\mathbf{U}\|_{0}\leq s.
\end{aligned}
\end{equation}
where $\mathcal{A}: \mathbb{U}\rightarrow \mathbb{R}^N$ is a linear operator, which can be explicitly expressed as
\[\mathcal{A}(\mathbf{U}) = [\langle \mathbf{A}_1, \mathbf{U}\rangle,\langle \mathbf{A}_2, \mathbf{U}\rangle, \cdots, \langle \mathbf{A}_N, \mathbf{U}\rangle],\,\mathbf{A}_i\in\mathbb{U}, i = 1,\cdots,N,\]
where $\mathbf{A}_i$ denotes the measurement matrix. The rank constraint $\operatorname{rank}(\mathbf{U})\leq r$ stipulates that the rank of matrix $\mathbf{U}$ must not exceed a given positive integer $r$. The cardinality constraint $\|\mathbf{U}\|_{0}\leq s$ requires that the number of nonzero entries of $\mathbf{U}$ does not exceed a given positive integer $s$. The symbol $\mathbb{U}_+$ represents the vector space $\mathbb{R}_+^{m\times n}$, which contains all nonnegative $m\times n$ matrices. Alternatively, $\mathbb{U}_+$ can also denote the vector spaces $\mathbb{S}_+^n$ and $\mathbb{H}_+^n$ of all $n\times n$ positive semidefinite and positive semidefinite Hermitian matrices, respectively. Let $\mathbb{C}_r=\{\mathbf{U}\in\mathbb{U}: \operatorname{rank}(\mathbf{U})\leq r\}$ and $\mathbb{C}_s=\{\mathbf{U}\in\mathbb{U}: \|\mathbf{U}\|_0\leq s\}$ denote the sets of solutions satisfying rank and cardinality constraints, respectively.\par

Solving problems with either cardinality or rank constraints is usually difficult since these problems are known to be NP-hard. To address this challenge, most existing methods for sparse and low-rank matrix estimation problems use a convex relaxation technique, leading to a tractable convex optimization problem. Several previous studies, such as \cite{Ref_richard2012estimation,Ref_chen2014estimation,Ref_ohlsson2012compressive,Ref_li2013sparse}, have employed nuclear norm regularization and $l_1$-norm regularization to replace the rank and cardinality constraints, respectively. As a result, the following convex problem was obtained:
\begin{equation}\label{eq4}
  \min_{\mathbf{U}\in\mathbb{U}_+}\; f(\mathbf{U})+\alpha\|\mathbf{U}\|_{*}+\beta\|\mathbf{U}\|_{1},
\end{equation}
where $\|\cdot\|_{*}$ denotes the nuclear norm of matrix, $\|\cdot\|_1$ represents the $l_1$-norm of matrix (the sum of the absolute values of the matrix elements), $\alpha>0$ and $\beta>0$ are the penalty parameters. Since the proximal map of $\|\mathbf{U}\|_{*}$ or $\|\mathbf{U}\|_1$ has a closed-form solution, problem \eqref{eq4} can be efficiently solved using some first-order methods. However, it is important to note that although this convex relaxation model promotes sparse and low-rank solutions, there is a significant difference between it and the original problem with rank and cardinality constraints in a general setting. This is because the nuclear norm and the $l_1$-norm are convex, while the rank function and cardinality function are nonconvex. Additionally, the convex relaxation problem is unable to generate solutions with the desired structure. For instance, the solutions may not satisfy the rank or cardinality constraints imposed on the problem. \par

In addition to the convex relaxation approach, a popular class of nonconvex methods that rely on some tractable proximal projection can be used to handle rank and cardinality constraints. One such method is the penalty decomposition (PD) method \cite{Ref_lu2012Sparse,Ref_lu2014Penalty}, which can be used to directly solve problems with rank or cardinality constraints. Additionally, by utilizing the advantages of PD method and  a proximal alternating linearized minimization (PALM) method \cite{Ref_bolte2014proximal}, Teng et al. \cite{Ref_teng2017penalty} proposed a penalty PALM (PPALM) method to solve a class of cardinality constrained problems. Notably, the PPALM method can also be used to directly solve \eqref{eq1}. However, to ensure its convergence, the proximal parameter must be greater than the gradient Lipschitz constant of the quadratic penalty subproblem, which can affect its convergence rate. \par

 Here, we introduce another novel nonconvex method to handle  rank  and cardinality constraints. 
In our previous work \cite{Ref_ding2022inexact}, the rank constraint was  transformed into  an equivalent difference-of-convex (DC) equality constraint, i.e.,
\[\operatorname{rank}(\mathbf{U})\leq r\iff \|\mathbf{U}\|_{*}-\|\mathbf{U}\|_{*(r)}=0,\]
where $\|\mathbf{U}\|_{*(r)}= \sum_{i=1}^r \sigma_i(\mathbf{U})$ denotes Ky Fan $r$-norm of matrix $\mathbf{U}$ and $\sigma_i(\mathbf{U})$ is the $i$-th largest singular value of $\mathbf{U}$. This method was also employed in some other literatures \cite{Ref_gao2010majorized,Ref_bi2016error,Ref_liu2020exact} for rank constrained optimization problems.  Likewise, the cardinality constraint can be  transformed into an equivalent DC equality constraint, i.e.,
\[\|\mathbf{U}\|_{0}\leq s\iff \|\mathbf{U}\|_{1}-\|\mathbf{U}\|_{(s)}=0,\]
where
 $\|\mathbf{U}\|_{(s)}$ represents the Ky Fan $s$-norm of matrix $\mathbf{U}$, can be denoted as the sum of the $k$ entries with the greatest absolute value of the matrix. However, due to the nonconvex and nonsmooth nature of these DC equality constraints, solving the resulting problem is inherently challenging. One intuitive approach is to penalize these two DC equality constraints in the objective function and then transform the original problem into a standard DC programming problem. However, since both constraints are penalized into the objective function simultaneously, it is exceedingly difficult or even impossible to obtain their exact penalty results. As such, it is theoretically unsuitable to utilize the exact penalty method to handle the rank and cardinality constraints simultaneously.\par

 In a recent work, Liu et al. \cite{Ref_liu2019successive} introduced a successive DC approximation method (SDCAM) to tackle nonconvex nonsmooth problems, which includes a problem of sparse and low-rank matrix estimation:
\begin{equation}\label{eq_5}
\begin{aligned}
\min_{\mathbf{U}\in\mathbb{R}^{m\times n}}&\;\frac{1}{2}\|\mathbf{U}-\mathbf{V}\|_F^2+\delta_{\mathbb{C}_s}(\mathbf{U})+\delta_{\mathbb{C}_r}(\mathbf{U}).
\end{aligned}
\end{equation}
 In each step of SDCAM, either $\delta_{\mathbb{C}_s}(\mathbf{U})$ or $\delta_{\mathbb{C}_r}(\mathbf{U})$ was replaced by its Moreau envelope. For a proper function $g:\mathbb{U}\rightarrow \mathbb{R}\cup\{\infty\}$ ($\operatorname{dom}g \neq\emptyset$), its Moreau envelope for any given $\mu>0$ can be defined as 
 \[M_{\mu, g}(\mathbf{U})=\inf_{\mathbf{V}\in \operatorname{dom} g} \left\{\frac{1}{2\mu}\|\mathbf{V}-\mathbf{U}\|_F^2+g(\mathbf{V})\right\}.\]
 By setting \[\mathbb{C}_{r,\tau}=\left\{\mathbf{U}\in\mathbb{C}_{r}: \sigma_{\max}(\mathbf{U})\leq \tau\right\},\quad \mathbb{C}_{s,\tau}=\left\{\mathbf{U}\in\mathbb{C}_{s}:  \max\{\vert\mathbf{U}_{i,j}\vert\}\leq \tau\right\},\]
  where $\tau \gg \max\{\max_{i,j}\vert\overline{\mathbf{U}}^*_{i,j}\vert, \sigma_{\max}(\overline{\mathbf{U}}^*)\}$ and $\overline{\mathbf{U}}^*$ is a global optimal solution of 
 \eqref{eq_5}.  Then \eqref{eq_5} was reformulated as the following forms:
 \begin{equation}\label{eq_3}
\begin{aligned}
\min_{\mathbf{U}\in\mathbb{R}^{m\times n}}\;\frac{1}{2}\|\mathbf{U}-\mathbf{V}\|_F^2+\delta_{\mathbb{C}_{s,\tau}}(\mathbf{U})+M_{\mu, \delta_{\mathbb{C}_r}}(\mathbf{U}),
\end{aligned}
\end{equation}
 \begin{equation}\label{eq_4}
\begin{aligned}
\min_{\mathbf{U}\in\mathbb{R}^{m\times n}}\;\frac{1}{2}\|\mathbf{U}-\mathbf{V}\|_F^2+\delta_{\mathbb{C}_{r,\tau}}(\mathbf{U})+M_{\mu, \delta_{\mathbb{C}_s}}(\mathbf{U}),
\end{aligned}
\end{equation}
where $M_{\mu, \delta_{\mathbb{C}_r}}(\mathbf{U})$ and $M_{\mu, \delta_{\mathbb{C}_s}}(\mathbf{U})$ are the Moreau envelope of $\delta_{\mathbb{C}_r}(\mathbf{U})$ and $\delta_{\mathbb{C}_s}(\mathbf{U})$, respectively.  Notice that $M_{\mu, g}(\mathbf{U})$
 can be rewritten into the following DC form:
\begin{equation}\label{eq12}
M_{\mu, g}(\mathbf{U})=\frac{1}{2 \mu}\|\mathbf{U}\|_F^{2}-\underbrace{\sup _{\mathbf{V} \in \operatorname{dom} g}\left\{\frac{1}{\mu}\langle \mathbf{U}, \mathbf{V}\rangle-\frac{1}{2 \mu}\|\mathbf{V}\|_F^{2}-g(\mathbf{V})\right\}}_{D_{\mu, g}(\mathbf{U})}.
\end{equation} 
This implies that Moreau smoothing method can be used to achieve the DC approximation of $\delta_{\mathbb{C}_{r,\tau}}(\mathbf{U})$ or $\delta_{\mathbb{C}_{s,\tau}}(\mathbf{U})$ for a given smoothing parameter $\mu>0$.
 Nevertheless, simply using this iterative Moreau smoothing approach to  $\delta_{\mathbb{C}_s}(\mathbf{U})$ or $\delta_{\mathbb{C}_r}(\mathbf{U})$ would result in  neither \eqref{eq_3} nor \eqref{eq_4} being a DC programming since nonconvex functions $\delta_{\mathbb{C}{r,\tau}}(\mathbf{U})$ lack an explicit DC decomposition. 
 In \cite{Ref_liu2019successive}, a nonmonotone proximal gradient method with majorization was proposed to solve the nonconvex subproblems of SDCAM. As a result, SDCAM did not fully utilize the DC structure of problem \eqref{eq_5}, and this might affect its numerical efficiency. Moreover, in theory, it is not feasible to adopt the Moreau smoothing method to approximate $\delta_{\mathbb{C}_s}(\mathbf{U})$ and $\delta_{\mathbb{C}_r}(\mathbf{U})$ simultaneously.\par

To address the challenges presented by the rank and cardinality constraints, we utilize the Moreau smoothing method for the indicator function of one of the constraints and the exact penalty method for the other constraint.
Let \[g(\mathbf{U})= \delta_{\mathbb{C}_{s,\tau}}(\mathbf{U}),\; P_1(\mathbf{U})=\|\mathbf{U}\|_{*},\;P_2(\mathbf{U})=\|\mathbf{U}\|_{*(r)},\] or 
\[g(\mathbf{U}) =\delta_{\mathbb{C}_{r,\tau}}(\mathbf{U}),\;P_1(\mathbf{U}) =\|\mathbf{U}\|_{1},\;P_1(\mathbf{U}) =\|\mathbf{U}\|_{(s)},\]
where $\tau \gg \max\{\max_{i,j}\vert\mathbf{U}^*_{i,j}\vert, \sigma_{\max}(\mathbf{U}^*)\}$ and $\mathbf{U}^*$ is a global optimal solution of 
 \eqref{eq1}.
Then \eqref{eq1}  can be reformulated into the following problem (with the same optimal solution):
  \begin{equation}\label{eq5}
    \begin{aligned}
     \min_{\mathbf{U}\in\mathbb{U}_+}\;& \ell(\mathbf{U})+g(\mathbf{U})\\
      \text{s.t.}\;& P_1(\mathbf{U})-P_2(\mathbf{U})=0.
    \end{aligned}
  \end{equation}
 Firstly,  we replace $g(\mathbf{U})$  with its Moreau envelope $M_{\mu, g}(\mathbf{U})$ with a given smoothing constant $\mu>0$. As a result, the following problem is obtained:
  \begin{equation}\label{eq6}
    \begin{aligned}
   \min_{\mathbf{U}\in\mathbb{U}_+}\;& J_{\mu}(\mathbf{U}) := \ell(\mathbf{U})+M_{\mu, g}(\mathbf{U})\\
      \text{s.t.}\;& P_1(\mathbf{U})-P_2(\mathbf{U})=0.
    \end{aligned}
  \end{equation}
  Notice that the difficulty caused by the DC constraint $P_1(\mathbf{U})-P_2(\mathbf{U})=0$ still exists. To tackle this issue, we employ a penalty approach to penalize $P_1(\mathbf{U})-P_2(\mathbf{U})$ into the objective function of \eqref{eq6}
 with a given penalty parameter $c>0$,
 \begin{equation}\label{eq7}
\min_{\mathbf{U}\in\mathbb{U}_+}\; J_{\mu,c}(\mathbf{U}) :=\ell(\mathbf{U})+M_{\mu, g}(\mathbf{U})+c P_1(\mathbf{U})-c P_2(\mathbf{U}).
 \end{equation}
 It is worth noting that the penalty method is exact, which can be proven. As shown in \eqref{eq12}, $M_{\mu, g}(\mathbf{U})$ can be expressed as a DC form. By doing so, \eqref{eq7} can be transformed into a standard DC programming, which can be solved under the framework of DC algorithm (DCA). Furthermore, we devise an asymptotic strategy to adjust the smoothing parameter $\mu>0$ and penalty parameter $c>0$, and propose an asymptotic DC method to effectively solve \eqref{eq1}. \par
 In the case of large-scale DC programming problems like $\min_{\mathbf{U}} f_1(\mathbf{U})-f_2(\mathbf{U})$, it may be impractical to solve the convex subproblems exactly, or it may not be required. It is important to note that solving an optimization problem with high precision using an iterative algorithm is computationally expensive. Therefore, it may not be efficient to solve all the subproblems of DCA to a high accuracy. Additionally, solving subproblems to high precision at the initial iteration of DCA may not be necessary. To mitigate the computational cost, one popular approach is to design a termination criterion for solving subproblems, ensuring that the entire algorithm framework is both convergent and numerically efficient. The
subproblem of DCA in step $k$ can be expressed as \begin{equation}\nonumber
\min_{\mathbf{U}}\, G(\mathbf{U}):= f_1(\mathbf{U})-\langle \mathbf{W}^k,\mathbf{U} \rangle,
\end{equation}
where $\mathbf{W}^k\in\partial f_2(\mathbf{U}^k)$, $f_1(\mathbf{U})$ is assumed to be strongly convex with parameter $\alpha>0$.
Let $\mathbf{U}^{k+1}$ be an approximate solutions of $\min_{\mathbf{U}}G(\mathbf{U})$, then there exists a $\Delta^{k+1}$ such that $\Delta^{k+1}\in \partial G(\mathbf{U}^{k+1})$.
To ensure the convergence of the general inexact DCA, $\mathbf{U}^{k+1}$ and $\Delta^{k+1}$ were required to satisfy the termination criteria similar to $\|\Delta^{k+1}\|_F\leq \eta_k\|\mathbf{U}^{k+1}-\mathbf{U}^k\|_F$ with a given constant $\frac{\alpha}{2}\geq\eta_k\geq 0$, see \cite{Ref_souza2016global}. However, because the vector $\Delta^{k+1}$ is implicitly related to $\mathbf{U}^{k+1}$, the algorithm for solving the subproblem may fail to terminate even if the subproblem is solved to high precision. As a result, implementing the termination criterion of $\|\Delta^{k+1}\|_F\leq \eta_k\|\mathbf{U}^{k+1}-\mathbf{U}^k\|_F$ can be challenging in practical applications, especially for general optimization algorithms. \par 
Building upon our previous work \cite{Ref_ding2022inexact}, we propose an efficient inexact DC algorithm with a sieving strategy (siDCA) to solve \eqref{eq7}. In siDCA, we first apply a simple and readily applicable termination criterion to the subproblem, after which we employ a "sieve" post-processing technique to refine the obtained solution. We prove that the sequence generated by siDCA converges globally to a stationary point of \eqref{eq7}. To ensure the efficiency of the overall inexact DCA, it is crucial to solve the subproblems of siDCA efficiently. We exploit the semismooth properties in the dual subproblems and the superlinear convergence of the semismooth Newton (SSN) method to adopt a dual-based SSN method for solving the subproblem in siDCA. \par
 
\subsection{Notations and preliminaries}
Below are some common notations to be used in this paper. Let $\boldsymbol{e}_i$ be the $i$-th standard unit vector. Given
an index set $\mathbb{L} \subset \{1, \cdots, q\}$, $\vert\mathbb{L}\vert$ denotes the size of $\mathbb{L}$. We denote the vector and matrix of all ones by $\boldsymbol{1}$ and $\mathbf{E}$, respectively. We denote the identity matrix by $\mathbf{I}$. We use $\mathbb{S}^q$ to denote the linear subspace of all $q\times q$ real symmetric matrices. We use $\mathbb{H}^q$ to denote the linear subspace of all $q\times q$ Hermitian matrices. According to $\mathbb{U}_+$, $\mathbb{U}$ can denote $\mathbb{R}^{m\times n}$, $\mathbb{S}^q$ and $\mathbb{H}^q$. Let $\mathbf{U}^{\top}$ be the conjugate transpose of $\mathbf{U}$. Notice that when $\mathbf{U}\in\mathbb{S}^q$, the conjugate transpose of $\mathbf{U}$ equals to its non-conjugate transpose.  Let $\|\cdot\|_F$ denote the Frobenius norm of matrices. Let $\|\mathbf{U}\|_*$ denote the nuclear norm of matrix $\mathbf{U}$. $\|\cdot\|$ is used to represent the $l_2$ norm of vectors and matrices, respectively.  $\langle\mathbf{U},\mathbf{A}\rangle = \sum_{i=1}^m\sum_{j=1}^n\mathbf{U}_{i,j}\mathbf{A}_{i,j}$ is used to denote the inner product between matrix $\mathbf{U}$ and matrix $\mathbf{A}$. We use `$\otimes$' to denote the matrix Kronecker product. Let $\lambda_1\geq\lambda_2\geq\cdots\geq\lambda_q$ be the eigenvalues of $\mathbf{U}\in \mathbb{S}^q$ ($\mathbf{U}\in \mathbb{H}^q$) being arranged in nonincreasing order. Let $\boldsymbol{\alpha}:=\{i\vert\lambda_i>0,i = 1,\cdots,q\}$ and $\boldsymbol{\beta}:=\{i\vert\lambda_i<0,i = 1,\cdots,q\}$ be denoted as the index set of positive eigenvalues and the index set of negative eigenvalues, respectively. As a result, the spectral decomposition of $\mathbf{U}$ is given as $\mathbf{U} = \mathbf{Q}\Lambda\mathbf{Q}^{\top}$ with $\Lambda = \operatorname{diag}(\Lambda_{\boldsymbol{\alpha}},\mathbf{0},\Lambda_{\boldsymbol{\beta}})$, where $\Lambda_{\boldsymbol{\alpha}}$ ($\Lambda_{\boldsymbol{\beta}}$) is a diagonal matrix with positive (negative) eigenvalues as diagonal elements.
Then the real-valued (complex-valued) positive semidefinite matrix cone projection of $\mathbf{U}$ is represented as \[\Pi_{\mathbb{S}_+^q}(\mathbf{U})=\mathbf{Q}_{\boldsymbol{\alpha}}\Lambda_{\boldsymbol{\alpha}}\mathbf{Q}_{\boldsymbol{\alpha}}^{\top} \quad (\Pi_{\mathbb{H}_+^q}(\mathbf{U})=\mathbf{Q}_{\boldsymbol{\alpha}}\Lambda_{\boldsymbol{\alpha}}\mathbf{Q}_{\boldsymbol{\alpha}}^{\top}),\]
where $\mathbf{Q}_{\boldsymbol{\alpha}}\in\mathbb{R}^{q\times\vert\boldsymbol{\alpha}\vert }$  is the sub-matrix of $\mathbf{Q}$ indexed by index set $\boldsymbol{\alpha}$.\par

\section{An asymptotic DC method for sparse and low-rank matrix recovery}\label{sec:2}
From \eqref{eq12} ,  it holds that problem \eqref{eq7} can be reformulated as follows:
\begin{equation}\label{eq.13}
  \begin{aligned}
 \min_{\mathbf{U}\in \mathbb{U}_+}\; J_{\mu,c}(\mathbf{U}) = \ell(\mathbf{U})+\frac{1}{2 \mu}\|\mathbf{U}\|_F^{2}-D_{\mu, g}(\mathbf{U})+cP_1(\mathbf{U})-cP_2(\mathbf{U}).
  \end{aligned}
\end{equation}
\begin{remark}\label{reamrk_1}
For \eqref{eq.13}, we choose $g(\mathbf{U})$ and $ P(\mathbf{U})$ according to $\mathbb{U}_+$ as follows: 
\begin{itemize}
  \item When $\mathbb{U}_+ = \mathbb{S}_+^n$ $(\mathbb{U}_+ = \mathbb{H}_+^n)$, it holds that $\|\mathbf{U}\|_* = \langle \mathbf{U}, \mathbf{I}\rangle$. Then we set $g(\mathbf{U}) = \delta_{\mathbb{C}_{s,\tau}}(\mathbf{U}),$
  \[P_1(\mathbf{U})=\|\mathbf{U}\|_{*},\quad P_2(\mathbf{U}) = \|\mathbf{U}\|_{*(r)},\quad P(\mathbf{U}) = P_1(\mathbf{U})-P_2(\mathbf{U}).\]
  \item When $\mathbb{U}_+ = \mathbb{R}_+^{m\times n}$, we have $\|\mathbf{U}\|_1 = \langle \mathbf{U}, \mathbf{E}\rangle$.
  Then we set $g(\mathbf{U}) = \delta_{\mathbb{C}_{r,\tau}}(\mathbf{U})$, \[P_1(\mathbf{U}) = \|\mathbf{U}\|_{1},\quad P_2(\mathbf{U}) = \|\mathbf{U}\|_{(s)},\quad P(\mathbf{U}) = P_1(\mathbf{U})- P_2(\mathbf{U}).\]
\end{itemize} 
\end{remark}

Since $\frac{1}{2 }\|\mathbf{V}\|_F^{2}+\mu g(\mathbf{V})$ may be nonconvex,  $D_{\mu, g}(\mathbf{U})$ may be a nonsmooth complex function and its subdifferential is difficult to compute. Fortunately, a subset of the subdifferential of $D_{\mu, g}(\mathbf{U})$ can be obtained by the following conclusion \cite{Ref_liu2019successive}.\par

\begin{proposition}\label{proP_1}
Let $D_{\mu, g}(\mathbf{U})$ be defined in \eqref{eq12}, then for any $\mu>0$,
 \begin{equation}\nonumber
   \frac{1}{\mu} \operatorname{Prox}_{\mu g}(\mathbf{U}) \subseteq \partial D_{\mu, g}(\mathbf{U}).
 \end{equation}
\end{proposition}
\subsection{Exact penalty property}
Based on Lipschitz continuity of the objective function, Bi et al. \cite{Ref_bi2016error} and Liu et al. \cite{Ref_liu2020exact} have studied the exact penalty properties of some rank constrained optimization problems. However,  since $\mathbb{U}_+$ is unbounded, then $J_{\mu}(\mathbf{U})$ is not global Lipschitz continuous on $\mathbb{U}_+$. Then exact penalty results in \cite{Ref_bi2016error} and \cite{Ref_liu2020exact} can not be applied to \eqref{eq.13}.\par
Let $\widehat{\mathbf{U}}_{\mu}^*$ be a global optimal solution of the following problem:
 \begin{equation}\label{eq.14}
   \begin{aligned}
    \min_{\mathbf{U}\in\mathbb{U}_+}\;J_{\mu}(\mathbf{U})= \ell(\mathbf{U})+M_{\mu,g}(\mathbf{U}).
   \end{aligned}
 \end{equation}
Let $\mathbf{U}_{\mu}^*$ be a global optimal solution of \eqref{eq6} and $\mathbf{U}_{\mu,c}^*$ be a global optimal solution of \eqref{eq.13}. We 
 first introduce  an exact penalty result in Proposition \ref{proP_2}, which was also used in \cite{Ref_bai2015constrained} and \cite{Ref_gao2010majorized}.\par

\begin{proposition}\label{proP_2}
Let $\mathbf{U}_f$ be a feasible solution of \eqref{eq6}. Assume that for a given $\epsilon>0$,  $c>0$ is a chosen constant such that $(J_{\mu}(\mathbf{U}_f)-J_{\mu}(\widehat{\mathbf{U}}_{\mu}^*))/c\leq \epsilon$. Then 
 \begin{equation}\label{eq_10}
P(\mathbf{U}_{\mu,c}^*)\leq\epsilon \quad \text{ and }\quad J_{\mu}(\mathbf{U}_{\mu,c}^*)\leq J_{\mu,c}(\mathbf{U}_{\mu,c}^*)\leq J_{\mu}(\mathbf{U}_{\mu}^*)= \bar{J}.
\end{equation}
\end{proposition}
 From Proposition~{\upshape\ref{proP_2}}, it is easy to see that an $\epsilon$-optimal solution to \eqref{eq6} in the sense of \eqref{eq_10} is guaranteed by solving \eqref{eq.13} with a chosen penalty parameter $c$. This result provides the rationale to replace the DC constraint in \eqref{eq.13} by the penalty term $c P(\mathbf{U})$.  Let $\mathbb{R}=\left\{\mathbf{U}\in\mathbb{U}: P(\mathbf{U})=0\right\}$, $\mathbb{F}=\{\mathbf{U}\in\mathbb{U}_+:P(\mathbf{U})=0\}.$
  For any $\mathbf{U}\in \mathbb{U}_+$, we have the following results about the orthogonal projection $\Pi_{\mathbb{R}}(\mathbf{U})$. \par

  \begin{lemma}\label{lemma_1}
  For any $\mathbf{U}\in \mathbb{U}_+$, let $\Pi_{\mathbb{R}}(\mathbf{U})$ be an orthogonal projection of $\mathbf{U}$ on $\mathbb{R}$, then $\Pi_{\mathbb{R}}(\mathbf{U})\in \mathbb{F}$ and
  \begin{equation}\label{eq_11}
 \begin{aligned}
  \|\mathbf{U}-\Pi_{\mathbb{R}}(\mathbf{U})\|_F\leq P(\mathbf{U}).
  \end{aligned}
  \end{equation}
  \end{lemma}
For a given $\epsilon >0$, let $\overline{c}$ be the chosen constant such that \[(J_{\mu}(\Pi_{\mathbb{R}}(\widehat{\mathbf{U}}_{\mu}^*))-J_{\mu}(\widehat{\mathbf{U}}_{\mu}^*))/\overline{c}\leq \epsilon.\]
Based on the $\epsilon$ and $\overline{c}$, we give an exact penalty parameter as follows:
\begin{equation}\label{eq17}
 \hat{c} = \max\left\{\bar{c},\frac{1}{2}\|\mathbf{A}\|_F^2\epsilon+ \sqrt{2\bar{J}}\|\mathbf{A}\|_F+\frac{1}{2\mu}\epsilon+\phi\tau\right\},
\end{equation}
where $\mathbf{A}$ is the matrix representation of the linear operator $\mathcal{A}$ and $\overline{J} = J_{\mu}(\mathbf{U}_{\mu}^*)$, $\tau>0$ is a given constant such that $\tau\gg \max\{\max_{i,j}\vert\mathbf{U}_{i,j}^*\vert,\sigma_{\max}(\mathbf{U}^*)\}$ with $\mathbf{U}^*$ being a global optimal solution of \eqref{eq1}, $\phi$ is the number of all entries of the matrix $\mathbf{U}\in\mathbb{U}$, $\phi=n^2$ when $\mathbb{U} = \mathbb{S}^n$ ($\mathbb{U} = \mathbb{H}^n$) and $\phi=mn$ when $\mathbb{U} = \mathbb{R}^{m\times n}$.
Based the exact penalty parameter $\hat{c}$ defined in \eqref{eq17}, we have the following conclusions.

 \begin{proposition}\label{prop_3}
 Let  $\hat{c}$ be the constant defined in \eqref{eq17}. Then for any $c>\hat{c}$, 
\begin{equation}\label{eq_13}
  J_{\mu}(\Pi_{\mathbb{R}}(\mathbf{U}_{\mu,c}^*))\leq J_{\mu}(\mathbf{U}_{\mu,c}^*)
+\hat{c}P(\mathbf{U}_{\mu,c}^*).
\end{equation}
\begin{proof}
From Proposition \ref{proP_1}, it holds that $P(\mathbf{U}_{\mu,c}^*)\leq\epsilon$ and \[ J_{\mu}(\mathbf{U}_{\mu,c}^*)=\frac{1}{2}\|\mathcal{A}(\mathbf{U}_{\mu,c}^*)- \boldsymbol{b}\|^{2}+\frac{1}{2 \mu}\|\mathbf{U}_{\mu,c}^*\|_F^{2}-D_{\mu, g}(\mathbf{U}_{\mu,c}^*)\leq \bar{J},\,\forall\, c>\hat{c}\geq\overline{c}.\]
It follows from  $M_{\mu,g}(\mathbf{U})= \frac{1}{2 \mu}\|\mathbf{U}\|_F^{2}-D_{\mu, g}(\mathbf{U})\geq 0$ that
\begin{equation}\label{eq.19}
  \|\mathcal{A}(\mathbf{U}_{\mu,c}^*)- \boldsymbol{b}\|\leq\sqrt{2\bar{J}}.
\end{equation}
Since $\Pi_{\mathbb{R}}(\mathbf{U}_{\mu,c}^*)\in \mathbb{F}$ is an orthogonal projection of $\mathbf{U}_{\mu,c}^*$ on $\mathbb{R}$, then 
\begin{equation}\label{eq_16}
  \langle\Pi_{\mathbb{R}}(\mathbf{U}_{\mu,c}^*)-\mathbf{U}_{\mu,c}^*,\mathbf{U}_{\mu,c}^*\rangle = 0.
\end{equation}
From the defination of $g(\mathbf{U})$, it holds that $\forall\,\overline{\mathbf{W}}\in \operatorname{Prox}_{\mu g}(\mathbf{U}_{\mu,c}^*)$, $\|\overline{\mathbf{W}}\|_F< N\tau$. This, together with the convexity of $D_{\mu, g}(\mathbf{U})$, yields that
\begin{equation}\label{eq.21}
\begin{aligned}
     D_{\mu, g}(\Pi_{\mathbb{R}}(\mathbf{U}_{\mu,c}^*))-D_{\mu, g}(\mathbf{U}_{\mu,c}^*)
     \geq -\phi\tau\|\Pi_{\mathbb{R}}(\mathbf{U}_{\mu,c}^*)-\mathbf{U}_{\mu,c}^*\|_F.
\end{aligned}
\end{equation}
Evidently,
\begin{equation}\label{eq.22}
  \begin{aligned}
   \|\mathcal{A}(\Pi_{\mathbb{R}}(\mathbf{U}_{\mu,c}^*)-\mathbf{U}_{\mu,c}^*)\|
   \leq \|\mathbf{A}\|_F\|\Pi_{\mathbb{R}}(\mathbf{U}_{\mu,c}^*)-\mathbf{U}_{\mu,c}^*\|_F.
  \end{aligned}
\end{equation}
It holds from Lemma \ref{lemma_1} that $\Pi_{\mathbb{R}}(\mathbf{U}_{\mu,c}^*) \in\mathbb{F}$ and
\begin{equation}\label{eq.23}
\begin{aligned}
  \|\Pi_{\mathbb{R}}(\mathbf{U}_{\mu,c}^*)-\mathbf{U}_{\mu,c}^*\|_F\leq P(\mathbf{U}_{\mu,c}^*).
\end{aligned}
\end{equation}
Consequently,
 \begin{equation}\nonumber
\begin{aligned}
&J_{\mu}(\Pi_{\mathbb{R}}(\mathbf{U}_{\mu,c}^*))\\ =& \frac{1}{2}\|\mathcal{A}(\Pi_{\mathbb{R}}(\mathbf{U}_{\mu,c}^*))- \boldsymbol{b}\|^{2}+\frac{1}{2 \mu}\|\Pi_{\mathbb{R}}(\mathbf{U}_{\mu,c}^*)\|_F^{2}-D_{\mu, g}(\Pi_{\mathbb{R}}(\mathbf{U}_{\mu,c}^*))\\
 =& J_{\mu}(\mathbf{U}_{\mu,c}^*)+\frac{1}{2 \mu}\|\Pi_{\mathbb{R}}(\mathbf{U}_{\mu,c}^*)-\mathbf{U}_{\mu,c}^*\|_F^{2}+D_{\mu, g}(\mathbf{U}_{\mu,c}^*)-D_{\mu, g}(\Pi_{\mathbb{R}}(\mathbf{U}_{\mu,c}^*))\\
&+\frac{1}{2}\|\mathcal{A}(\Pi_{\mathbb{R}}(\mathbf{U}_{\mu,c}^*)-\mathbf{U}_{\mu,c}^*)\|^{2} +
\langle\mathcal{A}(\mathbf{U}_{\mu,c}^*)- \boldsymbol{b}, \mathcal{A}(\Pi_{\mathbb{R}}(\mathbf{U}_{\mu,c}^*)-\mathbf{U}_{\mu,c}^*)\rangle\\
 \leq& J_{\mu}(\mathbf{U}_{\mu,c}^*) +\frac{1}{2 \mu}\|\Pi_{\mathbb{R}}(\mathbf{U}_{\mu,c}^*)-\mathbf{U}_{\mu,c}^*\|_F^{2}+N\tau\|\mathbf{U}_{\mu,c}^*-\Pi_{\mathbb{R}}(\mathbf{U}_{\mu,c}^*)\|_F\\
&+\frac{1}{2}\|\mathcal{A}(\Pi_{\mathbb{R}}(\mathbf{U}_{\mu,c}^*)-\mathbf{U}_{\mu,c}^*)\|^{2}+
\|\mathcal{A}(\mathbf{U}_{\mu,c}^*)- \boldsymbol{b}\|\|\mathcal{A}(\Pi_{\mathbb{R}}(\mathbf{U}_{\mu,c}^*)-\mathbf{U}_{\mu,c}^*)\| \\
\leq& J_{\mu}(\mathbf{U}_{\mu,c}^*)
+(\tfrac{1}{2}\|\mathbf{A}\|_F^2\epsilon+ \sqrt{2\bar{J}}\|\mathbf{A}\|_F+\tfrac{1}{2\mu}\epsilon+\phi\tau)P(\mathbf{U}_{\mu,c}^*) \\
\leq& J_{\mu}(\mathbf{U}_{\mu,c}^*)
+\hat{c}P(\mathbf{U}_{\mu,c}^*),
\end{aligned}
\end{equation}
where the second equality follows from \eqref{eq_16}, the first inequality holds from \eqref{eq.21}, the second inequality holds from \eqref{eq.19}, \eqref{eq.22}, \eqref{eq.23} and $P(\mathbf{U}_{\mu,c}^*)\leq \epsilon$, the last inequality is due to the definition of $\hat{c}$.
This completes the proof.
 \end{proof}
 \end{proposition}

 Based on Proposition \ref{prop_3}, we have the following conclusion about the exact penalty property of \eqref{eq.13}.
\begin{theorem}\label{thm_1}
Let $\hat{c}>0$ be the constant defined in \eqref{eq17}.
Then $\forall\;c>\hat{c}$, $\overline{\mathbf{U}}_{\mu,c}$ is a global minimizer of \eqref{eq6} if and only if $\overline{\mathbf{U}}_{\mu,c}$ is a global minimizer of \eqref{eq.13}.
\end{theorem}

Notice that for any $\mathbf{U},\widetilde{\mathbf{U}}\in\mathbb{B}(\overline{\mathbf{U}},\varepsilon)$ with $\varepsilon>0$, by setting $\widetilde{\mathbf{W}}\in\operatorname{Prox}_{\mu g}(\widetilde{\mathbf{U}})$, then it holds that
\begin{equation}\nonumber
  \begin{aligned}
   &J_{\mu}(\mathbf{U})-J_{\mu}(\widetilde{\mathbf{U}})\\
   =& \tfrac{1}{2}\|\mathcal{A}(\mathbf{U})-\boldsymbol{b}\|^2+\tfrac{1}{2 \mu}\|\mathbf{U}\|_F^{2}-D_{\mu, g}(\mathbf{U})-(\tfrac{1}{2}\|\mathcal{A}(\widetilde{\mathbf{U}})-\boldsymbol{b}\|^2+\tfrac{1}{2 \mu}\|\widetilde{\mathbf{U}}\|_F^{2}-D_{\mu, g}(\widetilde{\mathbf{U}}))\\
   \leq& \max_{ \mathbf{V}\in\mathbb{B}(\overline{\mathbf{U}},\varepsilon)}\|\mathcal{A}^{*}(\mathcal{A}(\mathbf{V})-\boldsymbol{b})+\tfrac{1}{\mu}\mathbf{V}\|_F \|\mathbf{U}-\widetilde{\mathbf{U}}\|_F+\|\widetilde{\mathbf{W}}\|_F\|\mathbf{U}-\widetilde{\mathbf{U}}\|_F\\
    \leq& \max_{ \mathbf{V}\in\mathbb{B}(\overline{\mathbf{U}},\varepsilon)}\|\mathcal{A}^{*}(\mathcal{A}(\mathbf{V})-\boldsymbol{b})+\tfrac{1}{\mu}\mathbf{V}\|_F \|\mathbf{U}-\widetilde{\mathbf{U}}\|_F+\phi\tau\|\mathbf{U}-\widetilde{\mathbf{U}}\|_F\\
   =& L_{J_{\mu}}^{\varepsilon}(\overline{\mathbf{U}})\|\mathbf{U}-\widetilde{\mathbf{U}}\|_F,
  \end{aligned}
\end{equation}
where $L_{J_{\mu}}^{\varepsilon}(\overline{\mathbf{U}}) = \max_{ \mathbf{V}\in\mathbb{B}(\overline{\mathbf{U}},\varepsilon)}\|\mathcal{A}^{*}(\mathcal{A}(\mathbf{V})-\boldsymbol{b})+\tfrac{1}{\mu}\mathbf{V}\|_F+\phi\tau$,
 the first inequality is due to  the convexity of $D_{\mu, g}(\mathbf{U})$, the second inequality is due to $\widetilde{\mathbf{W}}\in\operatorname{Prox}_{\mu g}(\widetilde{\mathbf{U}})$ and $\|\widetilde{\mathbf{W}}\|_F\leq \phi\tau$. In addition, we have
 $J_{\mu}(\widetilde{\mathbf{U}})-J_{\mu}(\mathbf{U})\leq L_{J_{\mu}}^{\varepsilon}(\overline{\mathbf{U}})\|\mathbf{U}-\widetilde{\mathbf{U}}\|_F.$
 As a result, it holds that $\vert J_{\mu}(\mathbf{U})- J_{\mu}(\widetilde{\mathbf{U}})\vert\leq L_{J_{\mu}}^{\varepsilon}(\overline{\mathbf{U}})\|\mathbf{U}-\widetilde{\mathbf{U}}\|_F$, which means $J_{\mu}(\mathbf{U})$ is local Lipschitz continuous with constant  $L_{J_{\mu}}^{\varepsilon}(\overline{\mathbf{U}})$. Then we give the following local exact penalty property of \eqref{eq.13}. \par
\begin{theorem}\label{thm_2} 
For any local minimizer $\overline{\mathbf{U}}_{\mu}\in \mathbb{F}$ of \eqref{eq6},  $\forall\,c>L_{J_{\mu}}^{\varepsilon}(\overline{\mathbf{U}}_{\mu})$, $\overline{\mathbf{U}}_{\mu}$ is also a local minimizer of \eqref{eq.13}. Conversely, if  $\overline{\mathbf{U}}_{\mu}\in \mathbb{U}_+$ is a local minimizer of \eqref{eq.13} with $c>0$ and $P(\overline{\mathbf{U}}_{\mu})=0$, then $\overline{\mathbf{U}}_{\mu}$ is a local minimizer of \eqref{eq6}.
\end{theorem}
Following the proof of \cite[Theorem 1 and Theorem 2]{Ref_ding2022inexact}, one can easily prove  Theorem \ref{thm_1} and Theorem \ref{thm_2}.
\subsection{An asymptotic DC method}
 Given a fixed $\mu>0$, obtaining the optimal solutions $\mathbf{U}_{\mu}^*$ and $\widehat{\mathbf{U}}_{\mu}^*$ for \eqref{eq.13} is a challenging task, which makes it impractical to directly compute $\hat{c}$ for \eqref{eq.13}. However, by utilizing Proposition \ref{proP_2} and Theorem \ref{thm_1}, we can adopt a strategy of gradually increasing the penalty parameter $c$. It is worth noting that when $c$ is sufficiently large ($c>\hat{c}$), any global optimal solution of \eqref{eq.13} is also a global optimal solution of \eqref{eq6}. Nevertheless, obtaining a global or local optimal solution for general DC problems is quite difficult. Then, for given $\mu>0$ and $c>0$, we approximate the solution of \eqref{eq.13}, and then gradually decrease $\mu$ while increasing $c$. This way, we propose an asymptotic DC method (ADC) to solve \eqref{eq5}.\par

\begin{algorithm}[ht]
\caption{An asymptotic DC method for (\ref{eq5})}
\label{alg_ ADC}
\begin{algorithmic}
\State\begin{enumerate}[leftmargin = 2.8em]
\item[\textbf{Step 0}] Pick four sequences of positive numbers with $\varepsilon_t\downarrow 0$, $\mu_t\downarrow 0$, $\epsilon_t\downarrow 0$ and $c_t$. Choose a feasible solution $\mathbf{U}^{0} \in\mathbb{F}$. Give $\rho>1$. Let $t = 0$.
\item[\textbf{Step 1}]  If $J_{\mu_t}(\Pi_{\mathbb{R}}(\mathbf{U}^t))\leq J_{\mu_t}(\mathbf{U}^{0})$, set $\mathbf{U}^{t+1,0}= \Pi_{\mathbb{R}}(\mathbf{U}^t)$. Else, set $\mathbf{U}^{t+1,0}=\mathbf{U}^{0}$. 
\item[\textbf{Step 2}] Set $c_{t,0}= c_t$ and $i=0$. Compute $\mathbf{U}^{t+1}$ as follows:
\begin{enumerate}[leftmargin=1.8em]
   \item[\textbf{2.1}] If $J_{\mu_t,c_{t,i}}(\mathbf{U}^{t+1,i})\leq J_{\mu_t,c_{t,i}}(\mathbf{U}^{t+1,0})$, set $\mathbf{U}^{t+1,i,0}:= \mathbf{U}^{t+1,i}$. Else, set $\mathbf{U}^{t+1,i,0}:=\mathbf{U}^{t+1,0}$.
    \item[\textbf{2.2}] Compute $\mathbf{U}^{t+1,i,k_i+1} \approx \mathop{\arg}\mathop{\min}\limits_{\mathbf{U}\in\mathbb{U}_+} J_{\mu_t,c_{t,i}}(\mathbf{U})$ by using Algorithm \ref{alg_siDCA} so that the following conditions hold: 
    \begin{align}
    \|\mathbf{U}^{t+1,i,k_i+1}- \mathbf{U}^{t+1,i, k_i}\|_F\leq \varepsilon_t, \label{eq.24}\\
    J_{\mu_t,c_{t,i}}(\mathbf{U}^{t+1,i,k_i})\leq J_{\mu_t,c_{t,i}}(\mathbf{U}^{t+1,i,0}).\label{eq.25}
    \end{align}
    \item[\textbf{2.3}]  Set $\mathbf{U}^{t+1,i+1}:=\mathbf{U}^{t+1,i,k_i}$. If $P(\mathbf{U}^{t+1,i,k_i})\leq \epsilon_t$, go to \textbf{Step 3}. Else, set $c_{t,i+1}:=\rho c_{t,i}$, $i:= i+1$ and go to \textbf{2.1}.
\end{enumerate}
\item[\textbf{Step 3}] Set $\mathbf{U}^{t+1}:=\mathbf{U}^{t+1,i+1}$, $t:= t+1$ and go to \textbf{Step 1}.
\end{enumerate}
\end{algorithmic}
\end{algorithm}





\begin{remark}\label{remark_1}
It should be pointed that the sequence of objective function values generated by Algorithm \ref{alg_siDCA} is monotonically decreasing, then the condition in \eqref{eq.25} is easy to satisfy. Let $\underline{\ell}:=\min_{\mathbf{U}\in\mathbb{U}_+} \ell(\mathbf{U})$, without loss of generality, we assume that $\ell(\mathbf{U}^{0})>\underline{\ell}$. If this is not the case, i.e., $\ell(\mathbf{U}^{0})\leq\underline{\ell}$, it holds from $\mathbf{U}^{0}\in\mathbb{F}$ that $\mathbf{U}^{0}$ is a global optimal solution of \eqref{eq5}.
\end{remark}

\begin{proposition}
For any $t>0$, let $\{\mathbf{U}^{t,i}\}_{i=1}^{\infty}$ be the sequence generated by inner iteration of Algorithm \ref{alg_ ADC}, then there exists \[i_0\leq\max\left\{\left\lfloor\frac{\log(\ell(\mathbf{U}^{0})-\underline{\ell})-\log(\epsilon_t c_t)}{\log(\rho)}+1\right\rfloor,1\right\}\] so that $P(\mathbf{U}^{t+1,i_0+1})\leq\epsilon_t$.
\begin{proof}
According to the definition of $\mathbf{U}^{t+1,i,0}$, we have that for a given $t\geq 0$, $\forall\, i\geq 0$, 
\[ J_{\mu_t,c_{t,i}}(\mathbf{U}^{t+1,i,0})\leq J_{\mu_t,c_{t,i}}(\mathbf{U}^{t+1,0}).\] 
From the definition of $\mathbf{U}^{t+1,0}$ and the feasibility of $\mathbf{U}^{0}$ to \eqref{eq5}, it holds that $P(\mathbf{U}^{t+1,0})=0$ and $M_{\mu_t, g}(\mathbf{U}^0)\leq g(\mathbf{U}^0)=0$. It follows from \eqref{eq.25} that for a given $t\geq 0$, $\forall\, i\geq 0$,
\begin{equation}\label{eq_4_27}
\begin{aligned}
J_{\mu_t,c_{t,i}}(\mathbf{U}^{t+1,i,k_i})\leq J_{\mu_t,c_{t,i}}(\mathbf{U}^{t+1,0})=J_{\mu_t}(\mathbf{U}^{t+1,0})\leq J_{\mu_t}(\mathbf{U}^{0})=\ell(\mathbf{U}^0).
\end{aligned}
\end{equation}
This, together with the nonnegativity of $M_{\mu_t,g}(\mathbf{U})$ and $c_{t,i}= c_t\rho^{i}$, yields that
\[P(\mathbf{U}^{t+1,i,k_i})\leq \frac{\ell(\mathbf{U}^{0})-J_{\mu_t}(\mathbf{U}^{t+1,i,k_i})}{c_{t,i}} \leq \frac{\ell(\mathbf{U}^{0})-\ell(\mathbf{U}^{t+1,i,k_i})}{c_{t,i}} \leq\frac{\ell(\mathbf{U}^{0})-\underline{\ell}}{c_t\rho^{i}},\]
where $\underline{\ell}= \min_{\mathbf{U}\in\mathbb{U}_+} \ell(\mathbf{U})$.
Then it follows from $\mathbf{U}^{t+1,i+1}= \mathbf{U}^{t+1,i,k_i}$ that when 
\[i\geq \frac{\log(\ell(\mathbf{U}^{0})-\underline{\ell})-\log(\epsilon_t c_t)}{\log(\rho)},\]
$P(\mathbf{U}^{t+1,i+1})\leq \epsilon_t$. This implies that there exists 
\[i_0\leq\max\left\{\left\lfloor\frac{\log(\ell(\mathbf{U}^{0})-\underline{\ell})-\log(\epsilon_t c_t)}{\log(\rho)}+1\right\rfloor,1\right\}\] 
so that $P(\mathbf{U}^{t+1,i_0+1})\leq\epsilon_t$. This completes the proof.
\end{proof}
\end{proposition}
\par

To ensure the objective function $J_{\mu,c}$ in \eqref{eq.13} is level-bounded, we assume $\mathcal{A}$ satisfies the Restricted Isometry Property (RIP) condition \cite{Ref_candes2008restricted}, which is also used as one of the most standard assumptions in the low-rank matrix recovery literatures \cite{Ref_cai2014sparse,Ref_luo2020recursive,Ref_bhojanapalli2016global}. Let $d := n$ when $\mathbb{U}=\mathbb{S}^n (\mathbb{H}^n)$ and $d:= \min(m,n)$ when $\mathbb{U}=\mathbb{R}^{m\times n}$.
\begin{definition}[RIP]\label{def_1}
 Let $\mathcal{A}:\mathbb{U}\rightarrow \mathbb{R}^N$ be a linear map.
For each $1\leq s \leq d$, define $s$-restricted isometry constant to be the
smallest number $R_s$ such that
\[(1-R_{s})\|\mathbf{U}\|_F^{2} \leq \|\mathcal{A}(\mathbf{U})\|^{2} \leq(1+R_{s})\|\mathbf{U}\|_F^{2}\]
holds for all $\mathbf{U}$ of rank at most $s$. And $\mathcal{A}$ is said to satisfy the {\rm $s$-RIP }condition if $0 \leq R_{s} < 1$. 
\end{definition}

\begin{theorem}\label{theorem_7} Let $\{\mathbf{U}^t\}$ be the sequence generated by Algorithm \ref{alg_ ADC}. By assuming that $\mathcal{A}$ satisfies {\rm $d$-RIP}, then the following statements hold.
\begin{itemize}
\item[$(1)$] $J_{\mu,c}(\mathbf{U})$ is lower bounded and level-bounded.
\item[$(2)$] The sequence $\{\mathbf{U}^t\}$ is bounded.
\item[$(3)$] Let $\mathbf{U}^*$ be an accumulation point of $\{\mathbf{U}^t\}$, then $\mathbf{U}^*\in \mathbb{F}$.
\end{itemize}
\begin{proof}
For statement (1), notice that $J_{\mu,c}(\mathbf{U})$ is the sum of nonnegative function $\frac{1}{2}\|\mathcal{A}(\mathbf{U})- \boldsymbol{b}\|^2$, $M_{\mu, g}(\mathbf{U})$ and $cP(\mathbf{U})$, then $J_{\mu,c}(\mathbf{U})$ is lower bounded. This, together with the {\rm $d$-RIP} of $\mathcal{A}$, yields that $\forall\,\mathbf{U}\in\mathbb{U}$,
\begin{equation}\nonumber
\begin{aligned}
J_{\mu,c}(\mathbf{U})& = \frac{1}{2}\|\mathcal{A}(\mathbf{U})- \boldsymbol{b}\|^2+M_{\mu, g}(\mathbf{U})+cP(\mathbf{U})\\
&\geq \frac{1}{2}\|\mathcal{A}(\mathbf{U})\|^2-\| \mathcal{A}(\mathbf{U})\|\|\boldsymbol{b}\|+\frac{1}{2}\|\boldsymbol{b}\|^2\\
&\geq \frac{1}{2}(1-R_d)^2\|\mathbf{U}\|_F^2-(1+R_d)\|\mathbf{U}\|_F\|\boldsymbol{b}\|+\frac{1}{2}\|\boldsymbol{b}\|^2.
\end{aligned}
\end{equation}
This implies that for any $\rho\in\mathbb{R}$, $\{\mathbf{U}\in\mathbb{U}\vert J_{\mu,c}(\mathbf{U})\leq \rho\}$ is bounded, i.e., $J_{\mu,c}(\mathbf{U})$ is level-bounded. \par
For statement (2), according to \eqref{eq_4_27} and $\mathbf{U}^{t+1}= \mathbf{U}^{t+1,i,k_i}$, it holds that
\begin{equation}\label{eq_4_31}
J_{\mu_t,c_{t,i}}(\mathbf{U}^{t+1})\leq \ell (\mathbf{U}^{0}).
\end{equation}
It follows from this and the level-boundedness of $J_{\mu_t,c_{t,i}}(\mathbf{U})$ that $\{\mathbf{U}^t\}$ is bounded.\par
For statement (3), since $\mathbf{U}^*$ is an accumulation point of $\{\mathbf{U}^t\}$, then there exists $\{\mathbf{U}^t\}_{t\in\mathbb{I}}$ so that \[\lim_{t\in\mathbb{I}}\mathbf{U}^t=\mathbf{U}^*.\] From the compactness of $\mathbb{U}_+$ and $\forall\,t\geq 0, \mathbf{U}^t\in\mathbb{U}_+$, we have $\mathbf{U}^*\in\mathbb{U}_+$. In addition, due to $\lim_{t\in\mathbb{I}}\epsilon_t = 0$, it holds that $P(\mathbf{U}^*)=0$. Since $g(\mathbf{U})$ is nonnegative, we get $\forall\,\mathbf{U}\in\mathbb{U}_+$,
\begin{equation}\nonumber
\begin{aligned}
&\frac{1}{2}\operatorname{dist}^2(\mathbf{U}, \operatorname{dom} g)\\
\leq& \inf_{\mathbf{V}\in\operatorname{dom} g} \frac{1}{2}\|\mathbf{V}-\mathbf{U}\|_F^2\leq  \inf_{\mathbf{V}\in\operatorname{dom} g} \frac{1}{2}\|\mathbf{V}-\mathbf{U}\|_F^2+\mu_t g(\mathbf{V})=\mu_t M_{\mu_t,g}(\mathbf{U}).
\end{aligned}
\end{equation}
This, together with \eqref{eq_4_31} and $c_{t,i}P(\mathbf{U})$, implies that
\begin{equation}\nonumber
\begin{aligned}
\frac{1}{2\mu_t}\operatorname{dist}^2(\mathbf{U}^{t+1}, \operatorname{dom} g)\leq M_{\mu_t,g}(\mathbf{U}^{t+1})\leq J_{\mu_t,c_{t,i}}(\mathbf{U}^{t+1})\leq \ell(\mathbf{U}^{0}).
\end{aligned}
\end{equation}
Due to $\lim_{t\in\mathbb{I}}\mu_t= 0$, by passing limitation to the above inequality along $t\in\mathbb{I}$, we get $\operatorname{dist}(\mathbf{U}^*,\operatorname{dom} g)=0$, which means that $\mathbf{U}^*\in\operatorname{dom} g$. Therefore, we can obtain $\mathbf{U}^*\in \mathbb{F}$.
This completes the proof.
\end{proof}
\end{theorem}

\section{Inexact DCA with sieving strategy}\label{sec:3}
Problem \eqref{eq.13} can be reformulated into a standard DC programming:
 \begin{equation}\label{eq.43}
\min_{\mathbf{U}\in\mathbb{U}_+} J_{\mu,c}(\mathbf{U}) = f_{\mu,c}(\mathbf{U})- h_{\mu,c}(\mathbf{U}),
\end{equation}
where \[f_{\mu,c}(\mathbf{U}) = \ell(\mathbf{U})+\frac{1}{2 \mu}\|\mathbf{U}\|_F^{2}+cP_1(\mathbf{U}),\quad h_{\mu,c}(\mathbf{U}) = D_{\mu_t, g}(\mathbf{U}) +cP_2(\mathbf{U})).\] Evidently, both $f_{\mu,c}(\mathbf{U})$ and $h_{\mu,c}(\mathbf{U})$ are convex functions.
Notice from Remark \ref{remark_1} that $f_{\mu,c}(\mathbf{U})$ is smooth. Moreover, $f_{\mu,c}(\mathbf{U})$ is a strongly convex function with parameter $\alpha = \frac{1}{\mu}$. 
Then \eqref{eq.43} can be formulated into the following form:
 \begin{equation}\label{eq.44}
 \min_{\mathbf{U}} f_{\mu,c}(\mathbf{U})+\delta_{\mathbb{U}_+}(\mathbf{U})-h_{\mu,c}(\mathbf{U}).
\end{equation}
The classical DCA for solving DC problem \eqref{eq.44} are presented in Algorithm \ref{alg_DCA}.

\begin{algorithm}[ht]
\caption{ DCA for \eqref{eq.44}}
\label{alg_DCA}
\begin{algorithmic}
\State
\begin{enumerate}[leftmargin=2.8 em]
\item[\textbf{Step 0}] Give penalty parameter $c>0$, smooth parameter $\mu>0$ and tolerance error $\varepsilon\geq 0$. Initialize $\mathbf{U}^{0} \in \mathbb{U}_+$. Choose $\mathbf{W}^{0}\in \partial h_{\mu,c}(\mathbf{U}^0)$. Set $k = 0$.
\item[\textbf{Step 1}] Compute $\mathbf{U}^{k+1}$ by solving
\begin{equation}\label{eq.45}
\min_{\mathbf{U}} G^k_{\mu,c}(\mathbf{U}):= f_{\mu,c}(\mathbf{U})+\delta_{\mathbb{U}_+}(\mathbf{U})-\langle \mathbf{U}, \mathbf{W}^{k}\rangle.
\end{equation}
\item[\textbf{Step 2}] If $\|\mathbf{U}^{k+1}-\mathbf{U}^{k}\|_F \leq \varepsilon$, stop and return $\mathbf{U}^{k+1}$.
\item[\textbf{Step 3}] Choose $\mathbf{W}^{k+1} \in\partial h_{\mu,c}(\mathbf{U}^{k+1})$, set $k:=k+1$ and go to $\textbf{Step 1}$.
\end{enumerate}
\end{algorithmic}
\end{algorithm}


\begin{remark}\label{remark_3}
We distinguish the following two cases to choose $\mathbf{W}^{k+1} \in\partial h_{\mu,c}(\mathbf{U}^{k+1})$:\\
\textbf{Case 1:} When $\mathbb{U}_+ = \mathbb{S}_+^n (\mathbb{H}_+^n)$,  $P_2(\mathbf{U}) = \|\mathbf{U}\|_{*(r)}$. By assuming that $\mathbf{U}^{k+1}$ has spectral decomposition $\mathbf{U}^{k+1} = \sum_{i= 1}^n \lambda_i\boldsymbol{q}_i\boldsymbol{q}_i^{\top}$ and $\mathbf{Y}^{k+1}\in\operatorname{Prox}_{\mu g}(\mathbf{U}^{k+1})$, then $\mathbf{W}^{k+1}$ can be chosen as $\mathbf{W}^{k+1}=c\sum_{i = 1}^r \boldsymbol{q}_i\boldsymbol{q}_i^{\top}+\frac{1}{\mu}\mathbf{Y}^{k+1}$.\\
\textbf{Case 2:} When $\mathbb{U}_+ = \mathbb{R}_+^{m\times n}$,  $P_2(\mathbf{U}) = \|\mathbf{U}\|_{(s)}$. Let $\mathbb{I}^{k+1}$ be an index set of corresponding to  $s$ largest elements of $\mathbf{U}$ in absolute value. Then $\mathbf{W}^{k+1}=\mathbf{X}^{k+1}+\frac{1}{\mu}\mathbf{Y}^{k+1}$, where $\mathbf{Y}^{k+1}\in\operatorname{Prox}_{\mu g}(\mathbf{U}^{k+1})$ and  $\mathbf{X}^{k+1}\in \partial cP_2(\mathbf{U}^{k+1})$ is chosen as
\begin{equation}\nonumber
  \mathbf{X}_{i,j}^{k+1} = \left\{\begin{array}{ll}
    c\operatorname{sign}(\mathbf{U}^{k+1}_{i,j}), & \text{if}\; (i,j)\in\mathbb{I}^{k+1} \\
    0, & otherwise
  \end{array}\right.  \quad 1\leq i\leq m, 1\leq j\leq n. 
\end{equation}
\end{remark}

It is important to note that problem \eqref{eq.45} cannot  obtain a closed-form solution and requires an iterative method for its solution. However, achieving a high level of accuracy in the solution through an iterative algorithm can be computationally expensive. Therefore, solving \eqref{eq.45} to a high degree of accuracy at every iteration of DCA can result in high computational costs and time consumption. Hence, it is necessary to establish a terminal condition for the subproblems to ensure that the whole algorithm framework is numerically efficient and convergent. Let $\mathbf{U}^{k+1}$ be an approximate solution of \eqref{eq.45}, i.e.,
\begin{equation}\nonumber
\mathbf{U}^{k+1}\approx\mathop{\arg}\mathop{\min}\limits_{\mathbf{U}} G^k_{\mu,c}(\mathbf{U}),
\end{equation}
then there exists $\Delta^{k+1}$ such that $\Delta^{k+1}\in \partial G_{\mu,c}^k(\mathbf{U}^{k+1})$.
To design an effective inexact DCA (iDCA), we need to establish a termination condition for its subproblems, to ensure the theoretical convergence of iDCA and efficient resolution of its subproblems.\par
\subsection{An efficient iDCA}\label{sec:3_1}
For some convex problems, it is only required that $\|\Delta^{k+1}\|_F\leq \epsilon_{k+1}$ and $\sum_{k=0}^{\infty}\epsilon_{k}<\infty$ in some traditional inexact algorithms \cite{Ref_jiang2012An,Ref_sun2015An}. Although this inexact strategy is simple and numerically implementable, it cannot guarantee the convergence of the corresponding iDCA. 
 An inexact condition similar to that of \cite{Ref_souza2016global,Ref_wang2019task},
\begin{equation}\label{eq.47}
 \|\Delta^{k+1}\|_F \leq \eta_k\|\mathbf{U}^{k+1}-\mathbf{U}^{k}\|_F,
\end{equation}
can ensure that the corresbounding iDCA is convergent, where $0\leq\eta_k \leq \frac{1}{2\mu}$. However, since $\Delta^{k+1}$ is related to $\mathbf{U}^{k+1}$ implicitly, then condition \eqref{eq.47} may not hold even if subproblem \eqref{eq.45} is solved to the higher accuracy. This means that termination condition \eqref{eq.47} is difficult to numerically implement of appropriate termination criterion for subproblem \eqref{eq.45}.\par
Based on our previous work \cite{Ref_ding2022inexact},  we first use $\|\Delta^{k+1}\|_F\leq\epsilon_{k+1}$ as the termination criterion for solving \eqref{eq.45}. Then we use  a sieving condition similar to \eqref{eq.47}
as a `sieve' to perform a post-processing on the approximate solution, where $\kappa\in (0,1)$ is a sieving parameter. Therefore, an inexact DCA with sieving strategy (siDCA) is given, see Algorithm \ref{alg_siDCA} for more details. \par

\begin{algorithm}[htb]
\caption{Inexact DCA with sieving strategy for (\ref{eq.44})}
\label{alg_siDCA}
\begin{algorithmic}
\State
\begin{enumerate}[leftmargin=2.8em] 
\item[\textbf{Step 0}] Give $c>0$, tolerance error $\varepsilon\geq 0$, nonnegative sequence $\epsilon_{k}\downarrow 0$, smooth parameter $\mu>0$ and sieving parameter $\kappa\in (0,1)$. Initialize $\mathbf{U}^{0} = \mathbf{V}^{0}\in \mathbb{U}_+$. Choose $\mathbf{W}^{0}\in \partial h_{\mu,c}(\mathbf{U}^0)$. Set $k = 0$.
\item[\textbf{Step 1}] Compute $\mathbf{V}^{k+1}$ by approximately solving
\begin{equation}\nonumber
\min_{\mathbf{U}} G^k_{\mu,c}(\mathbf{U}),
\end{equation}
so that there exists $\Delta^{k+1}\in \partial G^k_{\mu,c}(\mathbf{V}^{k+1})$, $\|\Delta^{k+1}\|_F\leq \epsilon_{k+1}$.
\item[\textbf{Step 2}] If $\|\mathbf{V}^{k+1}-\mathbf{U}^{k}\|_F \leq \varepsilon$ and $\|\Delta^{k+1}\|_F\leq \varepsilon$, stop and return $\mathbf{V}^{k+1}$. 
\item[\textbf{Step 3}] If sieving condition
\begin{equation}\label{eq.48}
\|\Delta^{k+1}\|_F < \frac{1}{2\mu}(1-\kappa)\|\mathbf{V}^{k+1}-\mathbf{U}^{k}\|_F
\end{equation}
holds, set $\mathbf{U}^{k+1} := \mathbf{V}^{k+1}$, choose $\mathbf{W}^{k+1} \in\partial h_{\mu,c}(\mathbf{U}^{k+1})$. Else, set $\mathbf{U}^{k+1} := \mathbf{U}^{k}$, $\mathbf{W}^{k+1} :=\mathbf{W}^{k}$. Set $k:=k+1$ and go to $\textbf{Step 1}$.
\end{enumerate}
\end{algorithmic}
\end{algorithm}
 


\subsection{A dual-based semismooth Newton method}
By noting that the main computation of siDCA is in solving subproblem \eqref{eq.45}, then we need to introduce an efficient method to solve \eqref{eq.45}. 
Let $\Phi_c^{k}=\mathbf{W}^k-c\mathbf{I}$ when $\mathbb{U}_+=\mathbb{S}_+^n (\mathbb{H}_+^n)$ and  $\Phi_c^{k}=\mathbf{W}^k-c\mathbf{E}$ when $\mathbb{U}_+= \mathbb{R}_+^{m\times n}$.
Evidently, \eqref{eq.45} can be reformulated as
\begin{equation}\label{eq48}
   \min_{\mathbf{U}}\frac{1}{2}\|\mathcal{A}(\mathbf{U})-\boldsymbol{b}\|^2+\frac{1}{2\mu}\|\mathbf{U}\|_F^2-\langle \Phi_c^{k},\mathbf{U}\rangle+\delta_{\mathbb{U}_+}(\mathbf{U}).
\end{equation}
Then the dual problem of \eqref{eq48} can be equivalently formulated as the following minimization problem:
\begin{equation} \label{eq.50}
\min_{\boldsymbol{z}\in\mathbb{R}^N,\mathbf{Y}\in\operatorname{dom}\delta_{\mathbb{U}_+}^*} H_{\mu,c}^k(\boldsymbol{z},\mathbf{Y})= \frac{1}{2} \|\boldsymbol{z}\|^{2}+\boldsymbol{z}^{\top}\boldsymbol{b} +\delta_{\mathbb{U}_+}^*(\mathbf{Y})+\frac{\mu}{2}\|\Phi_c^{k}-\mathcal{A}^*(\boldsymbol{z})-\mathbf{Y}\|_F^2.
\end{equation}
The KKT conditions of \eqref{eq.50} can be given as follows: 
  \begin{align}
\boldsymbol{z}+\boldsymbol{b} -\mu\mathcal{A}(\Phi_c^{k}-\mathcal{A}^*(\boldsymbol{z})-\mathbf{Y})=\boldsymbol{0},&\label{eq.51}\\
 \mathbf{0}\in \partial \delta_{\mathbb{U}_+}^*(\mathbf{Y})-\mu(\Phi_c^{k}-\mathcal{A}^*(\boldsymbol{z})-\mathbf{Y}). &\label{eq.52}
  \end{align}
  Let \[\Theta_{\mu,c}^k(\boldsymbol{z}):= \inf_{\mathbf{Y}\in\operatorname{dom} \delta_{\mathbb{U}_+}^*} H_{\mu,c}^k(\boldsymbol{z},\mathbf{Y}).\]
  Then the optimal solution $(\boldsymbol{z}^{\star},\mathbf{Y}^{\star})$ of \eqref{eq.50} can be obtained by
\begin{equation}\label{eq.53}
  \boldsymbol{z}^{\star} = \mathop{\arg}\mathop{\min}\limits_{\boldsymbol{z}\in\mathbb{R}^N} \Theta_{\mu,c}^k(\boldsymbol{z}), \;\mathbf{Y}^{\star} = \operatorname{Prox}_{\mu^{-1}\delta_{\mathbb{U}_+}^*}(\boldsymbol{t}^k(\boldsymbol{z}^{\star})), 
\end{equation}
  where $\boldsymbol{t}^k(\boldsymbol{z}) = \Phi_c^{k}-\mathcal{A}^*(\boldsymbol{z})$.
Then it holds that
 \begin{equation}\label{eq.55}
  \Theta_{\mu,c}^k(\boldsymbol{z})  = \frac{1}{2}\|\boldsymbol{z}\|^2 +\boldsymbol{z}^{\top}\boldsymbol{b}+ M_{\delta_{\mathbb{U}_+}^*}^{\mu^{-1}}(\boldsymbol{t}^k( \boldsymbol{z} )),
 \end{equation} 
 where $M_{\delta_{\mathbb{U}_+}^*}^{\mu^{-1}}(\boldsymbol{t}^k(\boldsymbol{z}))$ is Moreau envelope of $\delta_{\mathbb{U}_+}^*$ with parameter $\mu^{-1}$.
 From the extended Moreau decomposition $\boldsymbol{x} = \lambda\operatorname{Prox}_{\lambda^{-1}f^*}(\frac{\boldsymbol{x}}{\lambda})+\operatorname{Prox}_{\lambda f}(\boldsymbol{x})$, we have \[\operatorname{Prox}_{\mu^{-1} \delta_{\mathbb{U}_+}^*}(\boldsymbol{t}^k(\boldsymbol{z} )) = \boldsymbol{t}^k(\boldsymbol{z}) - \frac{1}{\mu}\operatorname{Prox}_{\mu \delta_{\mathbb{U}_+}}(\mu\boldsymbol{t}^k(\boldsymbol{z}))=\boldsymbol{t}^k(\boldsymbol{z}) - \frac{1}{\mu}\Pi_{\mathbb{U}_+}(\mu\boldsymbol{t}^k(\boldsymbol{z})).\]
 From \cite[Theorem 6.60 ]{Ref_beck2017first}, it holds that the Moreau envelope $M_{\delta_{\mathbb{U}_+}^*}^{\mu^{-1}}(\boldsymbol{t}^k( \boldsymbol{z} ))$ is a continuously differentiable function with gradient
 \[\nabla M_{\delta_{\mathbb{U}_+}^*}^{\mu^{-1}}(\boldsymbol{t}^k(\boldsymbol{z} )) = \mu(\boldsymbol{t}^k( \boldsymbol{z} )- \operatorname{Prox}_{\mu^{-1}\delta_{\mathbb{U}_+}^*}(\mu\boldsymbol{t}^k( \boldsymbol{z} )))=\Pi_{\mathbb{U}_+}(\mu\boldsymbol{t}^k(\boldsymbol{z})).\]
 Then the gradient of $\Theta_{\mu,c}^k(\boldsymbol{z})$ can be computed as
 \begin{equation}\nonumber
\nabla\Theta_{\mu,c}^k(\boldsymbol{z}) =\boldsymbol{z} +\boldsymbol{b}-\mathcal{A}(\Pi_{\mathbb{U}_+}(\mu\boldsymbol{t}^k(\boldsymbol{z}))).
\end{equation}
 By noting that $\Pi_{\mathbb{U}_+}(\mu\boldsymbol{t}^k(\boldsymbol{z}))$ is strongly semismooth, then it holds that $ \nabla\Theta_{\mu,c}^k(\boldsymbol{z})$
  is strongly semismooth. Hence, we employ an efficient semismooth Newton (SSN) method to solve \eqref{eq.50}. To this end, we define $\widehat{\partial}^2\Theta_{\mu,c}^k$ as follows: 
 \begin{equation}\nonumber
\widehat{\partial}^2\Theta_{\mu,c}^k(\boldsymbol{z}) :=\mathbf{I}+\mu\mathcal{A}\partial\Pi_{\mathbb{U}_+}(\mu\boldsymbol{t}^k(\boldsymbol{z}))\mathcal{A}^{*}, \forall\, \boldsymbol{z}\in \mathbb{R}^N,
\end{equation}
where $\partial\Pi_{\mathbb{U}_+}(\mu\boldsymbol{t}^k(\boldsymbol{z}))$ is the Clarke subdifferential of $\Pi_{\mathbb{U}_+}$ at point $\mu\boldsymbol{t}^k( \boldsymbol{z})$. From Hiriart-Urruty et al. \cite{Ref_hiriart1984Generalized}, it states that
\[\widehat{\partial}^2\Theta_{\mu,c}^k(\boldsymbol{z})\boldsymbol{d}_z= \partial^2\Theta_{\mu,c}^k(\boldsymbol{z})\boldsymbol{d}_z,\forall\, \boldsymbol{d}_z\in\
\mathbb{R}^N,\]
where $\partial^2\Theta_{\mu,c}^k(\boldsymbol{z})$ is the generalized Hessian of $\Theta_{\mu,c}^k(\boldsymbol{z})$ at $ \boldsymbol{z}$.
For any $\boldsymbol{z}^j\in\mathbb{R}^N$, we define a matrix $\mathbf{S}^j$ as follows:
 \begin{equation}\label{eq.56}
\mathbf{S}^j:=\mu\mathcal{A}\mathbf{H}^j\mathcal{A}^{*},
\end{equation}
where $\mathbf{H}^j \in\partial\Pi_{\mathbb{U}_+}(\mu\boldsymbol{t}^k( \boldsymbol{z^j}))$.
Consequently, it follows that $\mathbf{I}+\mathbf{S}^j\in \widehat{\partial}^2\Theta_{\mu,c}^k(\boldsymbol{z}^j)$ and $\mathbf{I}+\mathbf{S}^j\in \partial^2\Theta_{\mu,c}^k(\boldsymbol{z}^j)$. 
Thus, we present a SSN method to solve \eqref{eq.50}, see Algorithm \ref{algo_SSN}.\par

\begin{algorithm}[htb]
\caption{Semismooth Newton method for (\ref{eq.50})}
\label{algo_SSN}
\begin{algorithmic}
\State 
\begin{enumerate}[leftmargin=2.8em]
\item[\textbf{Step 0}] Initialize $\boldsymbol{z}^0\in \mathbb{R}^N$.  Give tolerance error $\chi_{k+1}$, line search parameter $\eta\in(0,1)$ and $\rho\in(0,1/2)$. Set $j = 0$.
\item[\textbf{Step 1}] Choose $\mathbf{H}^j \in\partial\Pi_{\mathbb{U}_+}(\mu\boldsymbol{t}^k( \boldsymbol{z^j}))$. Let $\mathbf{S}^j$ be defined in \eqref{eq.56}. Compute $\boldsymbol{d}_z^j$ by solving
\begin{equation}\label{eq.57}
   (\mathbf{I}+\mathbf{S}^j)\boldsymbol{d}_z^j = -\nabla\Theta_{\mu,c}^k(\boldsymbol{z}^j).
\end{equation}
\item[\textbf{Step 2}] (Line search) Find the least
nonnegative integer $u$ satisfying
\[\Theta_{\mu,c}^k(\boldsymbol{z}^j+\eta^u\boldsymbol{d}_z^j) \leq \Theta_{\mu,c}^k(\boldsymbol{z}^j)+\rho\eta^u\langle\nabla\Theta_{\mu,c}^k(\boldsymbol{z}^j),\boldsymbol{d}_z^j\rangle.\]
 Set $l_j = \eta^u$ and $\boldsymbol{z}^{j+1}= \boldsymbol{z}^j+l_j\boldsymbol{d}_z^j$.
\item[\textbf{Step 3}] If $\|\nabla\Theta_{\mu,c}^k(\boldsymbol{z}^{j+1})\|< \chi_{k+1}$, stop and return $\boldsymbol{z}^{k+1}:= \boldsymbol{z}^{j+1}$. Else, set $j := j+1$ and go to \textbf{Step 1}.
\end{enumerate}
\end{algorithmic}
\end{algorithm}



From primal-dual relationship, we can compute a feasible solution of \eqref{eq.45} at $j$-th iteration of Algorithm \ref{algo_SSN} as follows:
\begin{equation}\label{eq.58}
  \mathbf{U}^{(j)} = \mu(\Phi_c^{k}-\mathcal{A}^*(\boldsymbol{z}^{j})-\mathbf{Y}^{j}).
\end{equation}
Notice that $\mathbf{Y}^j=\Pi_{\mathbb{U}_+}(\mu\boldsymbol{t}^k(\boldsymbol{z}^j)),$ then we have 
\begin{equation}\label{eq.59}
  \mathbf{U}^{(j)} = \Pi_{\mathbb{U}_+}(\mu\boldsymbol{t}^k(\boldsymbol{z}^j))= \Pi_{\mathbb{U}_+}(\mathbf{U}^{(j)}+\mu\mathbf{Y}^{(j)}).
\end{equation}
Due to $\mathbf{Y} = \boldsymbol{t}^k(\boldsymbol{z})- \frac{1}{\mu} \Pi_{\mathbb{U}_+}(\mu\boldsymbol{t}^k(\boldsymbol{z}))$, then the KKT condition in \eqref{eq.52} holds for any $j>0$.  Let $\boldsymbol{\gamma}^j$ be defined as
 \begin{equation}\label{eq.60}
\boldsymbol{\gamma}^j := \nabla\Theta_{\mu,c}^k(\boldsymbol{z}^j) =\boldsymbol{z}^j +\boldsymbol{b}-\mathcal{A}(\Pi_{\mathbb{U}_+}(\mu\boldsymbol{t}^k( \boldsymbol{z}^j))).
\end{equation}
Since the residual of KKT conditions of \eqref{eq.50} is used to terminate Algorithm \ref{algo_SSN}, when $\|\boldsymbol{\gamma}^j\|\leq\chi_{k+1}$, we terminate Algorithm \ref{algo_SSN}.

\subsection{Convergence analysis of siDCA}\label{sec:3.4}
 A feasible point $\mathbf{U}\in\mathbb{U}_+$ is said to be a stationary point of DC problem \eqref{eq.44} if
\begin{equation}\nonumber
(\nabla f_{\mu,c}(\mathbf{U})+\partial \delta_{\mathbb{U}_+}(\mathbf{U}))\cap \partial h_{\mu,c}(\mathbf{U})\neq\emptyset.
\end{equation} 
The following results on the convergence of siDCA for solving \eqref{eq.44} follows from the basic convergence theorem of DCA \cite{Ref_tao1997convex}. Firstly, we can obtain the following conclusion about $\{\mathbf{U}^{k}\}$ and $\{J_{\mu,c}(\mathbf{U}^{k})\}$.\par

\begin{proposition}\label{prop_4.1} Let $\{\mathbf{U}^{k}\}$ be the stability center sequence generated by {\rm siDCA} for solving \eqref{eq.44}, then the following statements hold.
\begin{itemize}
\item[$(1)$] The sequence $\{J_{\mu,c}(\mathbf{U}^{k})\}$ is nonincreasing.
\item[$(2)$] The sequence $\{\mathbf{U}^{k}\}$ is bounded.
\end{itemize}
\end{proposition}

If $\mathbf{V}^{k+1}$ satisfies the sieving condition in \eqref{eq.48}, we say that a serious step is performed in siDCA, and say that $\mathbf{V}^{k+1}$ is a stability center. Else, we say that a null step is performed. When $\varepsilon$ is set to $0$, the following two situations would occur: (1) only finite serious steps are performed, then infinite null steps are performed or siDCA terminates after finite steps; (2) infinite serious steps are performed.\par

\subsubsection{Finite serious steps are performed in siDCA}
For the situation that only finite serious steps are performed in siDCA for solving \eqref{eq.44}, we have the following convergence results.
\begin{theorem}
\label{theo_4.1}
Set the tolerance error $\varepsilon = 0$. Suppose that only finite serious steps are performed in {\rm siDCA} for solving \eqref{eq.44}. Then the following statements hold:
\begin{itemize}
\item[$(1)$] If {\rm siDCA} terminates in finite steps, i.e., there exists a $\bar{k}>0$ such that $\mathbf{V}^{\bar{k}+1}= \mathbf{U}^{\bar{k}}$ and $\Delta^{\bar{k}+1} = \mathbf{0}$, then the stability center $\mathbf{U}^{\bar{k}}$ is a stationary point of \eqref{eq.44}.
\item[$(2)$] If after $\hat{k}$-th iteration of {\rm siDCA}, only null step is performed, then the stability center $\mathbf{U}^{\hat{k}+1}$ is a stationary point of \eqref{eq.44}.
\end{itemize}
\end{theorem}
\subsubsection{Infinite serious steps are performed in siDCA}
As siDCA involves infinite number of serious steps, there exist only a finite number of null steps between any two adjacent serious steps, and the stability center in the null step is a repetition of the stability center in the previous serious step. By removing the $\mathbf{U}^k$ generated in null steps from the set $\{\mathbf{U}^k\}$, we obtain a subsequence denoted by $\{\mathbf{U}^{k_l}\}$. We can also derive sequences $\{\mathbf{W}^{k_{l}}\}$ and $\{\Delta^{k_{l}}\}$ that correspond to $\mathbf{U}^{k_{l}}$. With these, we can prove the following global subsequential convergence results.

\begin{theorem}[Global subsequential convergence of siDCA]\label{thm:1}
Set the tolerance error $\varepsilon = 0 $. Let $\{\mathbf{U}^{k_l}\}$ be the stability center sequence generated in serious steps of {\rm siDCA} for solving \eqref{eq.44}. Then the following statements hold:
\begin{itemize}
\item[$(1)$] $\lim_{l\rightarrow \infty}\|\mathbf{U}^{k_l}-\mathbf{U}^{k_{l+1}}\|_F =0$.
\item[$(2)$] Any accumulation point $\overline{\mathbf{U}}\in\{\mathbf{U}^{k_l}\}$ is a stationary point of \eqref{eq.44}.
\end{itemize}
\end{theorem}

In order to display that $\{\mathbf{U}^k\}$ actually converges to a stationary point of \eqref{eq.44} when infinite serious steps are performed in Algorithm \ref{alg_siDCA}, we construct the following auxiliary function:
\begin{equation} \label{eq.76}
\begin{aligned}
&\Psi(\mathbf{U},\mathbf{W},\mathbf{V}, \mathbf{Z})\\
=& f_{\mu,c}(\mathbf{U})+\delta_{\mathbb{U}_+}(\mathbf{U})-\langle \mathbf{U},\mathbf{W}\rangle + h^*_{\mu,c}(\mathbf{W})+\frac{1}{2\mu}\|\mathbf{U} -\mathbf{V}\|^2_F -\langle \mathbf{Z}, \mathbf{U}-\mathbf{V}\rangle.
\end{aligned}
\end{equation}
where $h^*_{\mu,c}$ is the convex conjugate of $h_{\mu,c}$.
 According to \cite{Ref_attouch2009convergence,Ref_bolte2016majorization,Ref_bolte2014proximal}, it follows that semialgebraic functions satisfy Kurdyka-Łojaziewicz (KŁ) property. Evidently, the function
$\Psi$ belongs to the class of semialgebraic functions. 
 Based on $\Psi$, we have the following properties about $\{\Psi(\mathbf{U}^{k_{l+1}},\mathbf{W}^{k_l},\mathbf{U}^{k_l},\Delta^{k_{l+1}})\}$.
\begin{proposition}\label{prop_4.2}
 Let $\Psi$ be the function defined in \eqref{eq.76}. Let $\{\mathbf{U}^{k_l}\}$, $\{\Delta^{k_{l}}\}$ and $\{\mathbf{W}^{k_l}\}$ be the subsequences generated in serious steps of {\rm siDCA} for solving \eqref{eq.44}. Suppose that infinite serious steps are performed in {\rm siDCA} for solving \eqref{eq.44}, then the
following statements hold.\\
$(1)$ For any $l\geq 1$,
 \begin{equation} \label{eq.77}
J_{\mu,c}(\mathbf{U}^{k_{l+1}})
\leq \Psi(\mathbf{U}^{k_{l+1}},\mathbf{W}^{k_l},\mathbf{U}^{k_l},\Delta^{k_{l+1}}).
\end{equation}
$(2)$ For any $l\geq 1$, 
\begin{equation} \label{eq.78}
 \begin{aligned}
 &\Psi(\mathbf{U}^{k_{l+1}},\mathbf{W}^{k_l},\mathbf{U}^{k_l},\Delta^{k_{l+1}})-\Psi(\mathbf{U}^{k_l},\mathbf{W}^{k_{l-1}},\mathbf{U}^{k_{l-1}},\Delta^{k_{l}})\\
 &\leq -\frac{\kappa}{2\mu}\|\mathbf{U}^{k_l} -\mathbf{U}^{k_{l-1}}\|^2_F.
 \end{aligned}
\end{equation}
$(3)$ The set of accumulation points of the sequence $\{(\mathbf{U}^{k_{l+1}}, \mathbf{W}^{k_l}, \mathbf{U}^{k_l},\Delta^{k_{l+1}})\}$, denoted by $\Gamma$, is a nonempty compact set.
\\
$(4)$ The limit $\Upsilon = \lim_{l\rightarrow\infty}\Psi(\mathbf{U}^{k_{l+1}},\mathbf{W}^{k_l},\mathbf{U}^{k_l},\Delta^{k_{l+1}})$ exists and $\Psi \equiv \Upsilon$ on $\Gamma$.
\\
$(5)$ There exists a constant $\rho > 0$ such that for any $l \geq 1$, 
 \begin{equation} \label{eq.79}
\operatorname{dist}(\mathbf{0},\partial \Psi(\mathbf{U}^{k_{l+1}},\mathbf{W}^{k_l},\mathbf{U}^{k_l},\Delta^{k_{l+1}}))\leq \rho\|\mathbf{U}^{k_{l+1}}-\mathbf{U}^{k_l}\|_F.
\end{equation}
\end{proposition}
\begin{theorem}\label{thm:6}
Set the tolerance error $\varepsilon = 0$. Suppose that infinite serious steps are performed in {\rm siDCA} for solving \eqref{eq.44}. Let $\{\mathbf{U}^{k_l}\}$ be the stability center sequence generated in serious steps of {\rm siDCA} for solving \eqref{eq.44}. Then $\{\mathbf{U}^{k_l}\}$ converges to a stationary point of \eqref{eq.44}. Moreover, $\sum_{l=0}^{\infty}\|\mathbf{U}^{k_{l+1}}-\mathbf{U}^{k_{l}}\|_F<\infty$.
\end{theorem}

\begin{theorem}[Global sequential convergence of siDCA]\label{thm:3}
Set the tolerance error $\varepsilon = 0$. Suppose that infinite serious steps are performed in {\rm siDCA} for solving \eqref{eq.44}. Let $\{\mathbf{U}^{k}\}$ be the stability center sequence generated by {\rm siDCA} for solving \eqref{eq.44}. Then $\{\mathbf{U}^{k}\}$ converges to a stationary point of \eqref{eq.44}. Moreover, $\sum_{k=0}^{\infty}\|\mathbf{U}^{k+1}-\mathbf{U}^{k}\|_F<\infty$.
\end{theorem}

\section{Numerical experiments}\label{sec:5}
To illustrate the effectiveness of our ADC-siDCA for solving \eqref{eq1} in $\mathbb{R}^{m\times n}$ and $\mathbb{S}^{n}$, we compare its numerical performance with that of PPALM method and the SDCAM. In addition, we apply ADC-siDCA to solve sparse phase retrieval problem, which is an optimization problem in $\mathbb{H}_+^{n}$.\par
\vspace{0.2cm}
\noindent\textbf{Experimental environment:} All experiments are performed in Matlab 2020a on a 64-bit PC with double Intel(R) Xeon(R) CPU E5-2609 v2 (2.50GHz, 56GB of RAM).\par
\vspace{0.2cm}
\noindent\textbf{Problem scaling:}
To improve the performance of the algorithms, we perform a problem scaling: $\mathbf{A}_i = \mathbf{A}_i/\|\mathbf{A}_i\|_{F}, \boldsymbol{b}_i = \boldsymbol{b}_i/\|\mathbf{A}_i\|_{F}, i = 1,\cdots,p.$\par
\vspace{0.2cm}
To facilitate the use of PPALM and SDCAM, when $\mathbb{U}_+=\mathbb{S}_+^n (\mathbb{H}_+^n)$, we set \[\widehat{\mathbb{C}}_{r,\tau} = \left\{\mathbf{U}\in\mathbb{C}_r\cap\mathbb{U}_+: \sigma_1(\mathbf{U})\leq \tau\right\}, \quad \widetilde{\mathbb{C}}_{r} =\mathbb{C}_r\cap\mathbb{U}_+, \quad \widehat{\mathbb{C}}_{s,\tau} = \mathbb{C}_{s,\tau},\quad \widetilde{\mathbb{C}}_{s} =\mathbb{C}_s,\]  and set  
\[\widehat{\mathbb{C}}_{s,\tau} = \{\mathbf{U}\in\mathbb{C}_s\cap\mathbb{U}_+: \text{max}_{i,j}(\vert\mathbf{U}_{i,j}\vert)\leq \tau\}, \quad \widetilde{\mathbb{C}}_{s} =\mathbb{C}_s\cap\mathbb{U}_+, \;\widehat{\mathbb{C}}_{r,\tau} = \mathbb{C}_{r,\tau}, \quad \widetilde{\mathbb{C}}_{r} = \mathbb{C}_r\] when $\mathbb{U}_+=\mathbb{R}_+^{m\times n}$.\par
\vspace{0.2cm}
\noindent\textbf{SDCAM:}
It is easily to see that \eqref{eq1} can be rewritten into the form of $\min_{\mathbf{U}}\ell(\mathbf{U})+g(\mathbf{U})+P_0(\mathbf{U})$ (with same optimal value) in the following two ways:
\begin{equation}\label{eq.84}
\min_{\mathbf{U}}\quad F(\mathbf{U}) = \ell(\mathbf{U})+\underbrace{\delta_{\widetilde{\mathbb{C}}_{s}}(\mathbf{U})}_{g(\mathbf{U})}+\underbrace{\delta_{\widehat{\mathbb{C}}_{r,\tau}}(\mathbf{U})}_{P_0(\mathbf{U})},
\end{equation}
\begin{equation}\label{eq.85}
\min_{\mathbf{U}}\quad F(\mathbf{U}) = \ell(\mathbf{U})+\underbrace{\delta_{\widetilde{\mathbb{C}}_{r}}(\mathbf{U})}_{g(\mathbf{U})}+\underbrace{\delta_{\widehat{\mathbb{C}}_{s,\tau}}(\mathbf{U})}_{P_0(\mathbf{U})}.
\end{equation}
By replacing  $g(\mathbf{U})$ with its Moreau
envelope, we can get an auxiliary function
 \begin{equation}\nonumber
  F_{\mu}(\mathbf{U}) = \underbrace{\ell(\mathbf{U})+\frac{1}{2 \mu}\|\mathbf{U}\|_F^{2}}_{f_{\mu}(\mathbf{U})} +P_0(\mathbf{U})-D_{\mu, g}(\mathbf{U}),
\end{equation}
 where $f_{\mu}(\mathbf{U})$ is the smooth part of $F_{\mu}(\mathbf{U})$. The algorithm framework of SDCAM for solving \eqref{eq.84} and \eqref{eq.85} is presented in Algorithm \ref{alg_SDCAM}.\par

\begin{algorithm}[ht]
\caption{SDCAM for (\ref{eq.84}) and (\ref{eq.85})}
\label{alg_SDCAM}
\begin{algorithmic}
\State
\begin{enumerate}[leftmargin=2.8em] 
\item[\textbf{Step 0}] Pick two sequences of positive numbers with $\epsilon_t\downarrow 0$ and $\mu_t\downarrow 0$, a feasible point $\mathbf{U}^{0}\in \operatorname{dom} P_0\cap\operatorname{dom}g$. Initialize $\mathbf{U}^{0} \in \operatorname{dom}P_{0}$. Set $t = 0$.
\item[\textbf{Step 1}] If $F_{\mu_t}(\mathbf{U}^{t})\leq F_{\mu_t}(\mathbf{U}^{0})$, set $\mathbf{U}^{t,0}:=\mathbf{U}^{t}$. Else, set $\mathbf{U}^{t,0}:=\mathbf{U}^{0}$.
\item[\textbf{Step 2}] Approximately minimize $F_{\mu_t}(\mathbf{U})$, starting at $\mathbf{U}^{t,0}$, and terminating at $\mathbf{U}^{t,l_t}$ when
\begin{equation}\nonumber
\begin{aligned}
   &\operatorname{dist}(\mathbf{0},\nabla f(\mathbf{U}^{t,l_t})+\partial P_0 (\mathbf{U}^{t,l_{t}+1})-\tfrac{1}{\mu_t}\operatorname{Prox}_{\mu_t, g}(\mathbf{U}^{t,l_t}))\leq \epsilon_{t},\\
   &\|\mathbf{U}^{t,l_t}-\mathbf{U}^{t,l_{t}+1}\|_F\leq\epsilon_t \quad \text{and}\quad F_{\mu_t}(\mathbf{U}^{t})\leq F_{\mu_t}(\mathbf{U}^{t,0}).
\end{aligned}
\end{equation}
\item[\textbf{Step 3}] Update $\mathbf{U}^{t+1,0} := \mathbf{U}^{t,l_{t}}$. Set $t:=t+1$ and go to $\textbf{Step 1}$.
\end{enumerate}
\end{algorithmic}
\end{algorithm}

 
Liu et al. \cite{Ref_liu2019successive} called the method
based on \eqref{eq.84} SDCAM$_r$ and the method based on \eqref{eq.85} SDCAM$_s$. As presented in their numerical experiments,  SDCAM$_r$ was more efficient than SDCAM$_s$. An intuitive explanation is that rank constraint is a more complicated
constraint than the cardinality constraint to approximate via ‘subgradient’. As a result, we choose SDCAM$_r$ as one of comparison algorithms.
To solve the nonconvex nonsmooth subproblem of SDCAM:
\begin{equation}\label{eq.86}
\min_{\mathbf{U}}  F_{\mu_t}(\mathbf{U}) = f_{\mu_t}(\mathbf{U})+P_0(\mathbf{U})-D_{\mu, g}(\mathbf{U}),
\end{equation}
 a nonmonotone proximal gradient method with majorization (NPG$_{\text{\rm{major}}}$) \cite{Ref_liu2019successive} is used.\par
 \vspace{0.2cm}
\noindent\textbf{PPALM:}
 The PPALM method \cite{Ref_teng2017penalty} for \eqref{eq1} is presented in Algorithm \ref{alg_PPALM}.\par

 \begin{algorithm}[ht]
\caption{ PPALM method for (\ref{eq1}) }
\label{alg_PPALM}
\begin{algorithmic}
\State
\begin{enumerate}[leftmargin=2.8em]
\item[\textbf{Step 0}] Pick a positive sequence $\varepsilon_{k}\downarrow 0$. Let $L$ be the gradient Lipschitz constant of $\ell(\mathbf{U})$. Give $\rho_0>0$, $\sigma>1$, $\gamma_1, \gamma_2>1$, Initialize $\mathbf{U}^0\in\widetilde{\mathbb{C}}_{s}$ and $\mathbf{V}^0\in\widetilde{\mathbb{C}}_{r}$. Let $k = 0$.
\item[\textbf{ Step 1}] Set $l = 0$, $L_1^k = L+\rho_k$, $L_2^k = \rho_k$, $t^k_1 = \gamma_1 L_1^k$, $t^k_2=\gamma_2L_2^k$.  Applying PALM method by performing step \textbf{1.1}-\textbf{1.3}
to obtain a approximate critical point $(\mathbf{U}^{k+1}, \mathbf{V}^{k+1})\in\widetilde{\mathbb{C}}_{s}\times\widetilde{\mathbb{C}}_{r}$ of the following penalty subproblem:
\begin{equation}\nonumber
\min_{\mathbf{U},\mathbf{V}}\Phi_{\rho_k}(\mathbf{U},\mathbf{V})= \delta_{\widetilde{\mathbb{C}}_{s}}(\mathbf{U})+\ell(\mathbf{U})+\frac{\rho_k}{2}\|\mathbf{U}-\mathbf{V}\|_F^2+\delta_{\widetilde{\mathbb{C}}_{r}}(\mathbf{V}),
\end{equation}
\begin{enumerate}[leftmargin=1.8em]
\item[\textbf{1.1}] Compute $\mathbf{U}^{k,l+1}= \Pi_{\widetilde{\mathbb{C}}_{r}}(\mathbf{U}^{k,l}-1/ t^{k}_1(\nabla\ell(\mathbf{U}^{k,l})+\rho_k(\mathbf{U}^{k,l}-\mathbf{V}^{k,l})))$.
\item[\textbf{1.2}] Compute $\mathbf{V}^{k,l+1}= \Pi_{\widetilde{\mathbb{C}}_{s}}(\mathbf{V}^{k,l}-1/ t^{k}_2\rho_k(\mathbf{V}^{k,l}-\mathbf{U}^{k,l+1}))$.
\item[\textbf{1.3}] If $(\mathbf{U}^{k,l+1}, \mathbf{V}^{k,l+1})$ satisfies the inner terminal condition, set $(\mathbf{U}^{k+1}, \mathbf{V}^{k+1}) :=(\mathbf{U}^{k,l+1}, \mathbf{V}^{k,l+1})$ and 
go to \textbf{Step 2}. Else, let $l\leftarrow l+1$ and go to step \textbf{1.1}.
\end{enumerate}
\item[\textbf{ Step 2 }] If an outer iteration termination criterion is not met, set $\rho_{k+1} = \sigma\rho_k$ and $k:=k+1$, go to \textbf{Step 1}. Else, stop and return $\mathbf{U}^{k+1}$.
\end{enumerate}
\end{algorithmic}
\end{algorithm}


\par
 
\subsection{Sparse and low-rank nonnegative matrix recovery}\label{sec:5.1}
In this subsection, we perform matrix recovery experiments in $\mathbb{R}_+^{m\times n}$ to illustrate the effectiveness of ADC-siDCA. \par \vspace{0.2cm}
\noindent\textbf{Data generation:} We generate $\boldsymbol{b}$ by $\boldsymbol{b} = \mathcal{A}(\overline{\mathbf{U}})+\eta\boldsymbol{\theta}$, where \[\mathcal{A}(\overline{\mathbf{U}}) = [\langle\mathbf{A}_1,\overline{\mathbf{U}}\rangle, \langle\mathbf{A}_2,\overline{\mathbf{U}}\rangle,\cdots,\langle\mathbf{A}_N,\overline{\mathbf{U}}\rangle], \mathbf{A}_i\in\mathbb{R}^{m\times n}, i = 1,2,\cdots,N\] and $\boldsymbol{\theta}$ is the Gaussian white noise from a multivariate normal distribution $\mathbb{N}_N(\boldsymbol{0},\mathbf{I})$. The entries of $\mathbf{A}_i$ are randomly independently generated from standard Gaussian distribution. We set $(m, n) = (150,120), (200,160), (250,200)$, $N=16\max(m,n)$ and $\sigma = 0.01, 0.10$. We consider three types of sparse and low-rank nonnegative matrices as follows:
\begin{itemize}
  \item Cliques model (Cliq): $\overline{\mathbf{U}} = \operatorname{diag}(\mathbf{V},\mathbf{0}, \mathbf{V}, \mathbf{0}, \mathbf{V}, \mathbf{0}, \mathbf{V})$, where $\mathbf{V}\in\mathbb{R}_+^{m_1\times n_1}$ is a low-rank dense matrix. $\mathbf{V}$ is obtained by $\mathbf{V} = \mathbf{Y}\mathbf{W}^{\top}$, where the entries of $\mathbf{Y}\in\mathbb{R}_+^{m_1\times 3}$ and $\mathbf{W}\in\mathbb{R}_+^{n_1\times 3}$ are randomly generated to have independent identically distributed (i.i.d.) from a uniform distribution with range of $[0,1]$. We set $(m_1, n_1)=(\lfloor\frac{m}{6}\rfloor, \lfloor\frac{n}{6}\rfloor)$. 
  \item Random model \#1 (Rand1): $\overline{\mathbf{U}} = \mathbf{V}\otimes \mathbf{R}$, where $\mathbf{V}\in\mathbb{R}_+^{m_1\times n_1}$ is a low-rank dense matrix, $\mathbf{R}\in\mathbb{R}_+^{m_2\times n_2}$ is a random sparse matrix. $\mathbf{V}$ is obtained by $\mathbf{V} = \mathbf{Y}\mathbf{W}^{\top}$, where $\mathbf{Y}\in\mathbb{R}_+^{m_1\times 2}$ and $\mathbf{W}\in\mathbb{R}_+^{n_1\times 2}$ are randomly generated to have i.i.d. uniform distribution entries with range of $[0,10]$. $\mathbf{R}_{i,j}$ is assigned to be nonzero with probability $0.03$, independently of other elements. The nonzero elements of $\mathbf{R}$ are set to be $1$.  We set $(m_1, n_1) = (5, 4)$ and $(m_2, n_2) = (30, 30)$, $(40, 40)$, $(50, 50)$.
\item Random model \#2 (Rand2): $\overline{\mathbf{U}} = \mathbf{R}\otimes \mathbf{V}$, where $\mathbf{R}\in\mathbb{R}_+^{m_1\times n_1}$ is a random sparse matrix and $\mathbf{V}\in\mathbb{R}_+^{m_2\times n_2}$ is a low-rank dense matrix. 
The $\mathbf{R}_{i,j}$ is set to be nonzero with probability $0.03$, independently
of other elements. The nonzero elements of $\mathbf{R}$ are set to be $1$. $\mathbf{V}$ is obtained by $\mathbf{V} = \mathbf{Y}\mathbf{W}^{\top}$, where $\mathbf{Y}\in\mathbb{R}_+^{m_2\times 1}$ and $\mathbf{W}\in\mathbb{R}_+^{n_2\times 1}$ are randomly generated to have i.i.d. uniform distribution entries with range of $[0,10]$. We set $(m_1, n_1) = (30, 30)$, $(40, 40)$, $(50, 50)$ and $(m_2, n_2) = (5, 4)$.
\end{itemize}\par \vspace{0.2cm}
\noindent\textbf{Parameters setup:} In both SDCAM$_r$ and ADC-siDCA, we set $\mu_t=\mu_0/5^{t}$ with $\mu_0=50$, $\tau = 10^{5}$ and $\varepsilon_t=\varepsilon_0/1.5^t$ with $\varepsilon_0 = 10^{-4}$. In ADC-siDCA, the sieving parameter is set as $\kappa = 0.1$,   penalty parameter sequence is set as $c_t = 4^t c_0$ with $c_0=10^{-2}$. If $\mathbf{U}^k$ satisfies the sieving condition, the inexact error bound is updated by $\epsilon_{k+1}= \epsilon_k\max(\rho_1,\frac{k}{k_0+k})$ with $\rho_1 = 0.9$ and $k_0 = 20$. Else, set $\epsilon_{k+1}= \epsilon_k\max(\rho_2,\frac{k}{k_0+k})$ with $\rho_1 = 0.99$ and $k_0 = 20$. In NPG$_{\rm{major}}$, we set the parameters the same as those in \cite{Ref_liu2019successive}. For PPALM method, we set the parameters as: $\rho_0 = 0.05$, $\sigma =1.5$, $\gamma_1=1.01$, $\gamma_2=1.01$ and $\varepsilon_t=10^{-5}/1.2^t$. It should be noted that the inner tolerance error of PPALM method is set to be different from the other two methods because this parameter settings are mostly benefit to PPALM method.  In addition, to measure the quality of the solution $\widehat{\mathbf{U}}$ obtained by these methods, we define a matrix recovery error as $\text{MRE} = \frac{\|\widehat{\mathbf{U}}-\overline{\mathbf{U}}\|_F}{\max(1,\|\overline{\mathbf{U}}\|_F)}$.\par
 \vspace{0.2cm}
\noindent\textbf{Termination criteria:} We set termination condition for  SDCAM$_r$ and ADC-siDCA as $\max(\textbf{Vio}_r^t,\textbf{Vio}_s^t)\leq 10^{-9}$ or  $\mu_t\leq 10^{-9}$, where \textbf{Vio}$_r^t$ and \textbf{Vio}$_s^t$ are the violation of rank constraint and cardinality constraint, respectively, defined as \[\textbf{Vio}_r^t=\frac{\|\mathbf{U}^{t}-\Pi_{\mathbb{C}_r}(\mathbf{U}^{t})\|_F}{\max(1,\|\mathbf{U}^{t}\|_F)}, \textbf{Vio}_s^t = \frac{\|\mathbf{U}^{t}-\Pi_{\mathbb{C}_s}(\mathbf{U}^{t})\|_F}{\max(1,\|\mathbf{U}^{t}\|_F)}.\] For PPALM method, we terminate it when $\max(\textbf{Vio}_r^t,\textbf{Vio}_s^t)\leq 10^{-9}$ or $\rho_t>10^{9}$ holds. The termination criterion of siDCA is set as $\|\mathbf{V}^{k+1}-\mathbf{U}^{k}\|_F/\max(1,\|\mathbf{U}^{k}\|_F)\leq \varepsilon_{t}$ and $\|\Delta^{k+1}\|_F\leq \varepsilon_{t}$.
The termination criterion of NPG$_{\rm{major}}$ is set as $\|\mathbf{U}^{k+1}-\mathbf{U}^{k}\|_F/{\max(1,\|\mathbf{U}^{k}\|_F)}\leq \varepsilon_{t}$.
The inner termination criterion of PPALM is set as
\[\max \left\{\frac{\|\mathbf{U}^{k,l+1}- \mathbf{U}^{k,l}\|_{F}}{\max (\|\mathbf{U}^{k,l}\|_{F}, 1)}, \frac{\|\mathbf{V}^{k,l+1}- \mathbf{V}^{k,l}\|_{F}}{\max (\|\mathbf{V}^{k,l}\|_{F}, 1)}\right\} \leq \varepsilon_{k}.\]
\par \vspace{0.2cm}
As a comparison, numerical results for total solving time (Time/s), the number of outer iterations (Iter), matrix recovery error (MRE), optimal value of $\ell(\mathbf{U})$ (Obj), rank and sparsity (R,S) of optimal solution, obtained by ADC-siDCA, SDCAM$_r$ and PPALM method are presented in Table \ref{tab:1}. As a result from Table \ref{tab:1}, one can see that our ADC-siDCA outperforms SDCAM$_r$ and PPALM method for solving \eqref{eq1} when $\mathbb{U}_+=\mathbb{R}_+^{m\times n}$. The optimal objective function value and matrix recovery error obtained by ADC-iDCA are smaller than those of other methods in most cases. It should be noted that when $\mathbb{U}_+= \mathbb{R}_+^{m\times n}$, 
rank constraint is approximated by Moreau smoothing method in our ADC-siDCA. Although rank constraint is more difficult to approximate than cardinality constraint, it takes much less solving time and the number of outer iteration for ADC-siDCA than other two methods. Moreover, the solutions generated by ADC-siDCA satisfy rank constraint and cardinality constraint in most situations, while the solutions of SDCAM$_r$ and PPALM method sometimes do not achieve the required sparsity, especially for Rand1 model and Rand2 model. \par

\begin{sidewaystable}\setlength\tabcolsep{2.8pt}\small
\caption{The performance of  ADC-siDCA, SDCAM$_r$ and PPALM for solving (\ref{eq1}) when \texorpdfstring{$\mathbb{U}_+=\mathbb{R}_+^{m\times n}$}{}}\centering
\label{tab:1}    
\begin{center}
\begin{minipage}{\textheight}
\begin{tabular}{cccccccccccccccccccc}
\hline\noalign{\smallskip}
\multirow{2}{*}{m,n}& \multirow{2}{*}{$\eta$}&\multirow{2}{*}{r, s}&\multicolumn{5}{c}{ADC-siDCA}&\multicolumn{5}{c}{SDCAM$_r$}&\multicolumn{5}{c}{PPALM}\\
\cmidrule(lr){4-8}\cmidrule(lr){9-13}\cmidrule(lr){14-18}
\text{Data}&  &&time/s&Iter&MRE&obj&(R, S)&time/s&Iter& MRE & obj&(R, S)& time/s&Iter &MRE & obj&(R, S) \\[4pt]
\hline\noalign{\smallskip}
\multirow{2}{*}{150,120 }&0.01&12,2000&\textbf{52.72}&\textbf{33}&\textbf{6.21e-5}&\textbf{1.47e-1}&12, 2000& 169.00&3677&1.60e-4&4.41e-1&12,2000&137.83&4083&1.77e-4&5.66e-1&12,2000\\[4pt]
\text{Cliq}&0.10&12,2000&\textbf{56.54} &\textbf{42}&\textbf{2.98e-4}&\textbf{1.01e+1}&12,2000&158.45& 3295&2.61e-4&9.75e+0&12,2000&116.48&3814&3.00e-4&1.02e+1&12,2000\\[4pt]
\hline\noalign{\smallskip}
\multirow{2}{*}{200,160}&0.01&12,3432&\textbf{145.43}&\textbf{57}&2.26e-4&2.16e+0&12,3432& 451.53&4401&\textbf{2.76e-5}&\textbf{1.42e-1}&12,3432&283.47&4285&2.07e-4&1.81e+0&12,3432\\[4pt]
\text{Cliq}&0.10&12,3432& \textbf{153.22} &\textbf{58}&2.97e-4&1.50e+1&12,3432&394.99& 3824&\textbf{2.11e-4}&\textbf{1.33e+1}&12,3432&309.88&14138&2.16e-4&1.34e+1&12,3432\\[4pt]
\hline\noalign{\smallskip}
\multirow{2}{*}{250,200}&0.01&12,5412&\textbf{366.91}&\textbf{73}&2.76e-4&6.13+0&12,5412&866.08&4617&4.11e-4&1.17e+1&12,5412&533.82&4161&\textbf{2.66e-4}&\textbf{5.37e+0}&12,5412\\[4pt]
\text{Cliq}&0.10&12,5412&\textbf{283.16} &\textbf{60}&\textbf{2.14e-4}&\textbf{1.81e+1}&12,5412&707.51& 3751&3.14e-4&2.16e+1&12,5411&436.14&3445&3.06e-4&2.19e+1&12,5411\\[4pt]
\hline\noalign{\smallskip}
\multirow{2}{*}{150,120}&0.01&38,780&\textbf{58.41}&\textbf{51}&\textbf{3.30e-5}&\textbf{1.02e-1}&38,780 &222.18&4940&7.60e-4 &7.04e+0&38,1053&347.64&11609&2.57e-4&8.67e-1&40,780\\[4pt]
\text{Rand1}&.010&44,880&\textbf{83.18}&\textbf{51} &\textbf{3.67e-4}&\textbf{9.25e+0} &46,880&258.09&5834& 9.03e-4&1.62e+1&44,1034&529.25&16655&4.56-4&9.86e+0&44,977\\[4pt]
\hline\noalign{\smallskip}
\multirow{2}{*}{200,160}&0.01&48,1180& \textbf{142.56}&\textbf{48}&\textbf{4.16e-5}&\textbf{1.45e-1}&48,1180&2280.96&23914&7.02e-5&2.05e-1&48,1180&1543.48&20624&6.51e-5&1.91e-1&48,1180\\[4pt]
\text{Rand1}&0.10&58,1340& \textbf{241.45} &\textbf{90}&\textbf{2.40e-4}&\textbf{1.17e+1}&68,1340&2726.43& 27254&2.76e-4&1.23e+1&58,1717&2036.18&30195&4.29e-4&1.70e+1&58,1720\\[4pt]
\hline\noalign{\smallskip}
\multirow{2}{*}{250,200}&0.01&64,1540& \textbf{672.71}&\textbf{113}&\textbf{1.20e-5}&\textbf{4.61e-1}&65,1540&5061.06&25649&8.11e-5&1.37e+1&64,1599&3641.44&24648&4.17e-4&3.87e+2&64,1670\\[4pt]
\text{Rand1}&0.10&74,1600& \textbf{1301.05}&\textbf{187}&3.30e-5&1.67e+1&78,1600&8177.57&40475&\textbf{2.86e-5}&\textbf{1.55e+1}&74,1615&4174.59&28798&4.73e-4&7.27e+2&74,2404\\[4pt]
\hline\noalign{\smallskip}
\multirow{2}{*}{150,120}&0.01&17,800 & \textbf{42.97}&\textbf{29}&\textbf{3.48e-5}& 1.05e-1& 17,800&166.53&3763&3.51e-5&\textbf{1.04e-1}&17,800&57.69&1579&6.76e-5&1.29e-1&17,800\\[4pt] 
\text{Rand2}&0.10& 17,640& \textbf{36.20}&\textbf{31}&4.20e-4& \textbf{1.09e+1}&17,640&179.76&3851&4.15e-4&\textbf{1.09e+1}&17,640&49.15&1335&\textbf{4.06e-4}&\textbf{1.09e+1}&17,640\\[4pt] 
\hline\noalign{\smallskip}
\multirow{2}{*}{200,160}&0.01& 29,1300&\textbf{143.63}&\textbf{45}&\textbf{3.56e-5}&\textbf{1.43e-1}&29,1300&615.55&6304&2.78e-4&1.01e+0&29,1300&243.99&2954&1.00e-4&2.41e-1&29,1300\\[4pt]
\text{Rand2}&0.10&24,1180& \textbf{126.65}&\textbf{47}&4.59e-4&\textbf{1.40e+1}&25,1180&520.97&5223&4.95e-4&1.41e+1&24,1180&294.30&3920&\textbf{4.55e-4}&\textbf{1.40e+1}&24,1180\\[4pt]
\hline\noalign{\smallskip}
\multirow{2}{*}{250,200}&0.01&31,1500&  \textbf{355.96}&\textbf{69}&\textbf{6.51e-6}& \textbf{2.45e-1}&33,1500 &2603.93&14975&1.14e-4&2.58e+1&31,1817&1007.71&8357&2.34e-4&1.14e+2&31,1650\\[4pt]
\text{Rand2}
&0.10&33,1600&\textbf{314.65}& \textbf{59}&\textbf{4.67e-5}& \textbf{1.78e+1}&33,1600 &1558.82&8700&1.19e-4&2.68e+1&35,1600&1174.22&9724&2.56e-4&6.88e+1&33,1603\\[4pt]
\hline
\end{tabular}
\end{minipage}
\footnotetext[*]{The marked result is the one with the best numerical performance in the same experimental group (excluding the rank and sparsity of the matrix).}
\end{center}
\end{sidewaystable}

\subsection{Sparse and low-rank positive semidefinite matrix recovery}\label{sec:5.2}
In this subsection, we compare the performance of ADC-siDCA with SDCAM$_r$ and PPALM method for solving \eqref{eq1} when $\mathbb{U}_+ = \mathbb{S}_+^n$.\par
\noindent\textbf{Data generation:} Let $\boldsymbol{b} = \mathcal{A}(\overline{\mathbf{U}})+\eta\boldsymbol{\theta}$ with \[\mathcal{A}(\overline{\mathbf{U}}) = [\langle\mathbf{A}_1,\overline{\mathbf{U}}\rangle, \langle\mathbf{A}_2,\overline{\mathbf{U}}\rangle,\cdots,\langle\mathbf{A}_N,\overline{\mathbf{U}}\rangle], \mathbf{A}_i= \boldsymbol{a}_i\boldsymbol{a}_i^{\top}, \boldsymbol{a}_i\in\mathbb{R}^{n}, i = 1,\cdots,N,\] where $\boldsymbol{\theta}$ is the Gaussian white noise from a multivariate normal distribution $\mathbb{N}_N(\boldsymbol{0},\mathbf{I})$ and $\boldsymbol{a}_i$ is randomly generated from a multivariate normal distribution $\mathbb{N}_n(\boldsymbol{0}, \mathbf{I})$. We set $n= 200, 400, 600$,\;$N=10n$ and $\tau = 0.01, 0.10$. We consider three types of sparse and low-rank positive semidefinite matrices as follows:
\begin{itemize}
  \item Cliques model (Cliq): $\overline{\mathbf{U}} = \operatorname{diag}(\mathbf{V},\mathbf{0}, \mathbf{V},\mathbf{0},\mathbf{V},\mathbf{0},\mathbf{V}, \mathbf{0}, \mathbf{V})$, where $\mathbf{V} \in\mathbb{S}_+^{n_1}$ is a low-rank dense matrix. $\mathbf{V}$ is obtained by $\mathbf{V} = \mathbf{W}\mathbf{W}^{\top}$, where $\mathbf{W}\in\mathbb{R}^{n_1\times 2}$ is randomly generated to have i.i.d. standard Gaussian distribution entries. $n_1$ is set as $n_1 = \lfloor\frac{n}{10}\rfloor$. 
  \item Random model (Rand): $\overline{\mathbf{U}} = \mathbf{R}\otimes\mathbf{V} $, where $\mathbf{R}\in\mathbb{S}_+^{n_2}$ is low-rank sparse, $\mathbf{V}\in\mathbb{S}_+^{n_1}$ is a low-rank dense . The matrix $\mathbf{R}$, with around $0.03n^2$ nonzero entries, is generated by random Jacobi rotations applied to a nonnegative diagonal matrix $\Lambda$. 
  The matrix $\Lambda$ has at most $4$ nonzero diagonal entries, which are generated from a uniform distribution with range $[0,1]$. 
  We acquire $\mathbf{V}$ by $\mathbf{V} = \mathbf{W}\mathbf{W}^{\top}$, where $\mathbf{W}\in\mathbb{R}^{n_1\times 2}$ is randomly generated to have i.i.d. standard Gaussian distribution entries. 
\item Sparse random positive semidefinite model (Spr): $\overline{\mathbf{U}}\in\mathbb{S}_+^n$, with around $0.006n^2$ nonzero entries, is generated by random Jacobi rotations applied to a diagonal matrix with  given nonnegative eigenvalues $
\{\lambda_i\}_{i=1}^{n}$, where $\lambda_i$ is firstly generated from a uniform distribution with range $[0,100]$, independently to other elements. Then $\lambda_i$ is set to be zero with probability 0.98.
\end{itemize}
\par
\noindent\textbf{Parameters setup:} In both SDCAM$_r$ and ADC-siDCA, we set $\mu=\mu_0/5^{t}$ with $\mu_0=100$, $\tau = 10^{5}$ and $\varepsilon_t=\varepsilon_0/1.2^t$ with $\varepsilon_0 = 10^{-4}$.  The other parameters of siDCA, NPG$_{\text{\rm{major}}}$ and PPALM method are set the same as in Subsection \ref{sec:5.1}.\par
\noindent\textbf{Termination criteria:} We set the termination of  SDCAM$_r$ and ADC-siDCA  as $\max(\textbf{Vio}_r^t,\textbf{Vio}_s^t)\leq 10^{-9}$ or $\mu_t\leq 10^{-9}$.
For PPALM method, we terminate it when $\max(\textbf{Vio}_r^t,\textbf{Vio}_s^t)\leq 10^{-9}$ or $\rho_t>10^{9}$. The termination criteria of the inner iteration of ADC-siDCA, SDCAM$_r$ and PPALM method is set the same as Subsection \ref{sec:5.1}.\par

\begin{sidewaystable}\setlength\tabcolsep{2.8pt}\small
\caption{ The performance of  ADC-siDCA, SDCAM$_r$ and PPALM for sloving (\ref{eq1}) when $\mathbb{U}_+ = \mathbb{S}_+^n$}\centering
\label{tab:2}    
\begin{center}
\begin{minipage}{\textheight}
\begin{tabular}{cccccccccccccccccccc}
\hline\noalign{\smallskip}
\multirow{2}{*}{n}& \multirow{2}{*}{$\eta$}&\multirow{2}{*}{r, s}&\multicolumn{5}{c}{ ADC-siDCA}&\multicolumn{5}{c}{SDCAM$_s$}&\multicolumn{5}{c}{PPALM}\\
  \cmidrule(lr){4-8}\cmidrule(lr){9-13}\cmidrule(lr){14-18}
\text{Data}&  &&time/s&Iter&MRE&obj&(R, S)&time/s&Iter& MRE & obj&(R, S)& time/s&Iter &MRE & obj&(R, S) \\[4pt]
\hline\noalign{\smallskip}
\multirow{2}{*}{200}&0.01&10,2000& \textbf{42.37}&\textbf{59}&\textbf{3.86e-5}& \textbf{9.36e-2}&10,2000 &120.74&894&5.94e-5&9.64e-2&10,2000& 83.98&5248&3.74e-4&8.59e-1&10,2000\\[4pt]
\text{Cliq}
&0.10&10,2000& \textbf{42.46}& \textbf{69}&3.69e-4& 8.89e+0&10,2000& 137.01&943&\textbf{3.61e-4}&\textbf{8.85e+0}&10,2000&88.58&5535&5.14e-4&9.46e+0&10,2000\\[4pt]
\hline\noalign{\smallskip}
\multirow{2}{*}{400}&0.01&10,8000& \textbf{221.96}&\textbf{54}&\textbf{1.95e-5}& \textbf{3.10e-1}&10,8000 & 240.57&553&3.95e-4 &6.37e+0&10,8000&760.73&9797&7.81e-4&2.21e+1&10,8000\\[4pt] 
\text{Cliq}&0.10& 10,8000& \textbf{213.15}&\textbf{56}&\textbf{2.33e-4}&\textbf{1.92e+1}&10,8000& 233.55&540& 4.83e-4&2.58e+1&10,8000&624.34&8151&8.17e-4&3.81e+1&10,8000\\[4pt] 
\hline\noalign{\smallskip}
\multirow{2}{*}{600}&0.01& 10,18000& \textbf{699.97}&\textbf{71}&\textbf{7.93e-5}&\textbf{1.13e+0}&10,18000&2154.72&3222&2.33e-3&7.49e+2&10,18000&2255.28&11990&1.14e-3&1.68e+2&10,18000\\[4pt]
\text{Cliq}&0.10&10,18000& \textbf{743.08}&\textbf{72}&\textbf{1.63e-4}&\textbf{2.82e+1}&10,18000&2276.15&3356&2.48e-3&6.60e+2&10,18000&2718.37&13917&1.18e-3&1.61e+2&10,18000\\[4pt]

\hline\noalign{\smallskip}
\multirow{2}{*}{200}&0.01&3,1275& \textbf{17.30}&\textbf{25}&\textbf{2.93e-4}&9.55e-2&3,1275&29.12&261&2.95e-4&\textbf{9.53e-2}&3,1275&31.23&1767&5.73e-4&1.01e-1&3,1275\\[4pt]
\text{Rand}&0.10&3,1400& \textbf{52.95}&\textbf{216}&\textbf{5.82e-3}&9.32e+0&3,1561&111.27&468&6.68e-3&9.28e+0&3,1478&63.81&3661&5.85e-3&\textbf{9.26e+0}&3,1425\\[4pt]
\hline\noalign{\smallskip}
\multirow{2}{*}{400}&0.01&3,5000& \textbf{108.00}&\textbf{48}&\textbf{1.69e-4}&\textbf{1.98e-1}&3,5000 &239.14&654&1.51e-3 &8.34e-1&3,5302&386.66&4451&1.11e-3&5.46e-1&3,5000\\[4pt]
\text{Rand}&0.10&3,5100&\textbf{184.95}&\textbf{153}&2.71e-3&\textbf{1.91e+1} &3,5761&636.63&665& 2.43e-3&1.95e+1&3,6311&655.86&7802&\textbf{2.25e-3}&1.97e+1&3,5315\\[4pt]
\noalign{\smallskip}\hline
\multirow{2}{*}{600}&0.01&3,10750&\textbf{306.43}&\textbf{84}&\textbf{3.73e-4}&\textbf{2.96e-1}&3,10796& 644.26&737&1.30e-3&3.68e-1&3,10838&1317.56&6888&1.20e-3&5.29e-1&3,10766\\[4pt]
\text{Rand}&0.10&3,10525& \textbf{481.70}&\textbf{88}&\textbf{7.69e-5}&\textbf{2.94e+1}&3,10619&2002.02& 761&1.85e-3&4.38e+2&3,20231&1773.86&8542&1.13e-3&1.92e+2&3,10714\\[4pt]
\hline\noalign{\smallskip}
\multirow{2}{*}{200}&0.01&6, 231&\textbf{26.81}&\textbf{34}&\textbf{7.70e-6}&\textbf{1.03e-1}&6, 231& 38.81&555&2.35e-5&1.26e-1&6, 231&72.38&4376&1.28e-3&7.82e+1&6, 237\\[4pt]
\text{Spr}&0.10&4, 234&\textbf{41.45}&\textbf{110}&\textbf{8.45e-5} &\textbf{9.47e+0}&4, 290&204.00& 413&6.56e-5&9.53e+0&4, 754&44.10&2565&3.72e-4&1.40e+1&4, 243\\[4pt]
\hline\noalign{\smallskip}
\multirow{2}{*}{400}&0.01&8, 948& \textbf{235.91}&\textbf{144}&\textbf{3.40e-5}&\textbf{2.50e-1}&8,1058& 954.56&2166&1.06e-3&5.49e+1&8,7578&467.35&5791&1.26e-3&7.96e+1&8,2937\\[4pt]
\text{Spr}&0.10&9, 913&\textbf{246.91} &\textbf{128}&\textbf{1.28e-4}&\textbf{1.93e+1}&9, 940& 673.65&2658 &9.27e-4&7.34e+1&9,6363&452.11&5598&1.38e-3&1.44e+2&9,1966\\[3pt]
\hline\noalign{\smallskip}
\multirow{2}{*}{600}&0.01&9,2074&\textbf{722.98}&\textbf{136}&\textbf{4.81e-5}&\textbf{6.09e-1}&9,2243& 2236.80&3387&1.53e-3&2.89e+1&9,3805&1486.31&7254&1.61e-3&3.27e+2&9,2079\\[4pt]
\text{Spr}&0.10&17,2083&\textbf{901.50} &\textbf{107}&\textbf{1.24e-4}&\textbf{3.21e+1}&13,2212&1667.80& 976&7.17e-4&1.74e+2&17,4128&2425.39&11627&3.49e-3&3.38e+3&17,9335\\[4pt]
\hline
\end{tabular}
\end{minipage}
\footnotetext[*]{The marked result is the one with the best numerical performance in the same experimental group (excluding the rank and sparsity of the matrix).}
\end{center}
\end{sidewaystable}

As a comparison, numerical results obtained by ADC-siDCA, SDCAM$_r$ and PPALM method are presented in Table \ref{tab:2}. As a result from Table \ref{tab:2}, one can see that our ADC-siDCA outperforms SDCAM$_r$ and PPALM method for solving \eqref{eq1} in $\mathbb{S}_+^n$. The optimal objective function value and matrix recovery error obtained by ADC-iDCA are smaller than those of SDCAM$_r$ and PPALM method in most cases. In addition, the total solving time and the number of outer iterations of ADC-siDCA is less than those of other two methods. Moreover, the solutions generated by all these three methods satisfy rank constraint in all situations, while these solutions sometimes do not achieve the required sparsity. In particular, for the complex matrix recovery, such as Rand model and Spr model, the solutions of ADC-siDCA are more sparse than those of SDCAM$_r$ and PPALM method.\par
\subsection{Sparse phase retrieval}\label{sec:5.3}
 Due to the physical constraints, one can only
measure and record intensities of the Fourier coefficients of an optical
object. This gives rise to a problem of recovering a signal $\boldsymbol{x}\in\mathbb{C}^n$ from magnitude measurements, known as phase retrieval \cite{Ref_fienup1982phase}. Specially, we recover an unknown
$\bar{\boldsymbol{x}}\in\mathbb{C}^n$ from a small number of quadratic measurements of the form $\boldsymbol{b}_i = \vert\boldsymbol{a}_i^{\top}\bar{\boldsymbol{x}}\vert^2, i=1, \cdots, N$, i.e.,
\begin{equation}\nonumber
  \min_{\boldsymbol{x}}\sum_{i=1}^N(\vert\boldsymbol{a}_i^{\top}\bar{\boldsymbol{x}}\vert^2- \boldsymbol{b}_i)^2,
\end{equation}
where $\boldsymbol{a}_i\in \mathbb{C}^n$ is the measurement vector. Notice that it is difficult to recover $\bar{\boldsymbol{x}}$ from the above optimization problem, because it is nonconvex.
To address this problem, a lifting approach introduced by Balan et al. \cite{Ref_balan2009painless} is employed.
By letting $\mathbf{U} = \boldsymbol{x}\boldsymbol{x}^{\top}\in\mathbb{H}_+^n$, the above porblem can be reformulated as
 \begin{equation}
   \begin{aligned}
   \min_{\mathbf{U}\in\mathbb{H}_+^n}\,& \frac{1}{2}\|\mathcal{A}(\mathbf{U})-\boldsymbol{b}\|^2\\
   \text{s.t.}\,&\operatorname{rank}(\mathbf{U})\leq r,
   \end{aligned}
 \end{equation}
 where the linear operator $\mathcal{A}:\mathbb{H}_+^n\rightarrow\mathbb{R}^N$ can be explicitly expressed as
 \[\mathcal{A}(\mathbf{U}) = [\langle \mathbf{U},\boldsymbol{a}_1\boldsymbol{a}_1^{\top}\rangle,\langle \mathbf{U},\boldsymbol{a}_2\boldsymbol{a}_2^{\top}\rangle,\cdots,\langle \mathbf{U},\boldsymbol{a}_N\boldsymbol{a}_N^{\top}\rangle].\]
 The gained model is called PhaseLift model. Specifically, suppose that the true signal $\bar{\boldsymbol{x}}$ is sparse with $\bar{s}$ nonzero entries, then $\overline{\mathbf{U}} = \bar{\boldsymbol{x}}\bar{\boldsymbol{x}}^{\top}$ has $s=\bar{s}^2$ nonzero entries. This yields sparse PhaseLift (SPL) model. Evidently, problem \eqref{eq1}
subsumes SPL model as a special case with $s=\bar{s}^2$ and $r=1$. Then SPL model can be solved  by ADC-siDCA, PPALM method and SDCAM$_r$. \par
\noindent\textbf{Convex relaxation model:} To eliminate the difficulty caused by the rank constraint and cardinality constraint, a convex relaxation model of \eqref{eq1} can be employed for sparse phase retrieval, which is formulated as
\begin{equation}\label{eq.90}
  \min_{\mathbf{U}\in\mathbb{H}_+^n}\quad \frac{1}{2}\|\mathcal{A}(\mathbf{U})-\boldsymbol{b}\|^{2}+\alpha_r\|\mathbf{U}\|_{*}+\alpha_s\|\mathbf{U}\|_{1},
\end{equation}
where $\alpha_r$ and $\alpha_s$ are penalty parameters. This phase retrieval model is called convex sparse PhaseLift (CSPL) model. Since $\mathbf{U}\in\mathbb{H}^n_+$, then $\|\mathbf{U}\|_*=\langle \mathbf{U},\mathbf{I}\rangle$.  Then we can  formulate the dual problem of \eqref{eq.90} as the following equivalent minimization problem:
\begin{equation}\label{eq.91}
\begin{aligned}
   \min_{(\boldsymbol{z},\mathbf{X},\mathbf{Y})\in\mathbb{R}^N\times\mathbb{H}^n_+\times\mathbb{H}^n_+}\;&\frac{1}{2}\|\boldsymbol{z}\|^2+\boldsymbol{z}^{\top}\boldsymbol{b}+\delta_{\mathbb{H}_+^n}(\mathbf{X})+p^*(\mathbf{Y})\\
  \text{s.t.}\;& \mathcal{A}^*(\boldsymbol{z})+\alpha_r\mathbf{I}-\mathbf{X}+\mathbf{Y}=\mathbf{0},
\end{aligned}
\end{equation}
where $p(\mathbf{U}) = \alpha_s\|\mathbf{U}\|_1$ and $p^*$ denotes the convex conjugate of function $p$.  The augmented Lagrange function of \eqref{eq.91} can be expressed as
\begin{equation}\nonumber
  \begin{aligned}
   &\Psi_{\rho}(\boldsymbol{z}, \mathbf{X}, \mathbf{Y}; \mathbf{U} )\\=&\frac{1}{2}\|\boldsymbol{z}\|^2+\boldsymbol{z}^{\top}\boldsymbol{b}+\delta_{\mathbb{H}_+^n}(\mathbf{X})+p^*(\mathbf{Y})+\frac{\rho}{2}\left\|\mathcal{A}^*(\boldsymbol{z})+\alpha_r\mathbf{I}-\mathbf{X}+\mathbf{Y}-\frac{\mathbf{U}}{\rho}\right\|_F^2.
  \end{aligned}
\end{equation}
Evidently, problem \eqref{eq.91} belongs to a class of multi-block, nonsmooth, equality constrained convex optimization problems with coupled objective function.  To solve this problem, we presented an alternating direction method of multipliers (ADMM) in Algorithm \ref{algo_ADMM}.\par

\begin{algorithm}[htb]
\caption{Alternating direction method of multipliers method for \eqref{eq.91}}
\label{algo_ADMM}
\begin{algorithmic}
\State 
\begin{enumerate}[leftmargin=2.8em]
\item[\textbf{Step 0}]  Initialize $(\boldsymbol{z}^0,\mathbf{X}^0,\mathbf{Y}^0)\in\mathbb{R}^N\times\mathbb{H}^n_+\times\mathbb{H}^n_+$ and $\rho_0>0$. Give tolerance error $\varepsilon\geq 0$, $\alpha_r>0$, $\alpha_s>0$, $\beta>0$ and $\theta>1$. Set $k=0$.
\item[\textbf{Step 1}] Compute $\boldsymbol{z}^{k+1}:= \mathop{\arg}\mathop{\min}\limits_{\boldsymbol{z}}\Psi_{\rho_k}(\boldsymbol{z}, \mathbf{X}^k, \mathbf{Y}^k; \mathbf{U}^k ).$
\item[\textbf{Step 2}] Compute $\mathbf{X}^{k+1} := \mathop{\arg}\mathop{\min}\limits_{\mathbf{X}} \Psi_{\rho_k}(\boldsymbol{z}^{k+1}, \mathbf{X}, \mathbf{Y}^k; \mathbf{U}^k ).$
\item[\textbf{Step 3}] Compute $\mathbf{Y}^{k+1} := \mathop{\arg}\mathop{\min}\limits_{\mathbf{Y}} \Psi_{\rho_k}(\boldsymbol{z}^{k+1}, \mathbf{X}^{k+1}, \mathbf{Y}; \mathbf{U}^k ).$
\item[\textbf{Step 4}] Compute $\mathbf{U}^{k+1} := \mathbf{U}^{k+1}-\beta\rho_k(\mathcal{A}^*(\boldsymbol{z}^{k+1})+\alpha_r\mathbf{I}-\mathbf{X}^{k+1}+\mathbf{Y}^{k+1}).$ 
\item[\textbf{Step 5}] If $\max(\textbf{Pinf}_{k+1},\textbf{Dinf}_{k+1}) \leq \varepsilon$, stop. Else, set $k:=k+1$ and go to $\textbf{Step 1}$.
\end{enumerate}
\end{algorithmic}
\end{algorithm}
 



As presented in Algorithm \ref{algo_ADMM}, infeasibility of \eqref{eq.91} is used as the termination criterion and the basis for updating $\rho_k$. The primal infeasibility can be computed by $\textbf{Pinf}_{k+1} = \max(\eta_{\boldsymbol{z}}^{k+1}, \eta_{\mathbf{X}}^{k+1}, \eta_{\mathbf{Y}}^{k+1})$, where $(\eta_{\boldsymbol{z}}^{k+1}, \eta_{\mathbf{X}}^{k+1}, \eta_{\mathbf{Y}}^{k+1})$ are the residual of the first three KKT conditions:
\begin{equation}\nonumber
  \begin{aligned}
   &\eta_{\boldsymbol{z}}^{k+1} = \|\boldsymbol{z}^{k+1}+\boldsymbol{b}-\mathcal{A}(\mathbf{U}^{k+1})\|,\\
   & \eta_{\mathbf{X}}^{k+1} = \beta\rho_k\|\mathcal{A}^*(\boldsymbol{z}^{k+1})+\alpha_r\mathbf{I}-\mathbf{X}^{k+1}+\mathbf{Y}^{k+1}\|_F,\\
   &\eta_{\mathbf{Y}}^{k+1} = \|\mathbf{Y}^{k+1}-\mathbf{Y}^{k}+\beta\rho_k(\mathcal{A}^*(\boldsymbol{z}^{k+1})+\alpha_r\mathbf{I}-\mathbf{X}^{k+1}+\mathbf{Y}^{k+1})\|_F.
  \end{aligned}
\end{equation}
The dual infeasibility can be computed by $\textbf{Dinf}_{k+1} = \|\mathcal{A}^*(\boldsymbol{z}^{k+1})+\alpha_r\mathbf{I}-\mathbf{X}^{k+1}+\mathbf{Y}^{k+1}\|_F$. We update $\rho_{k+1}$  based on $(\textbf{Pinf}_{k+1},\textbf{Dinf}_{k+1})$: if $\textbf{Pinf}_{k+1}/\textbf{Dinf}_{k+1}\in (0.2,5)$, set $\rho_{k+1}:= \rho_{k}$.  Else, set $\rho_{k+1}:= \theta\rho_k$ if $\textbf{Pinf}_{k+1}/\textbf{Dinf}_{k+1}<0.2$, and set $\rho_{k+1}:= \rho_k/\theta$ if $\textbf{Pinf}_{k+1}/\textbf{Dinf}_{k+1}>5$. \par
\noindent\textbf{Data generation:}
We generate sparse phase retrieval instances as follows. First, draw $\bar{\boldsymbol{x}}\in\mathbb{C}^n$ from the multivariate complex standard normal distribution and set ninety percent entries in $\bar{\boldsymbol{x}}$ to zero randomly. The signal $\bar{\boldsymbol{x}}$ is normalized by $\bar{\boldsymbol{x}}= \frac{\bar{\boldsymbol{x}}}{\|\bar{\boldsymbol{x}}\|_{\infty}}$.  Then acquire $N = 10n$ phaseless measurements by $\boldsymbol{b}_i = \vert\boldsymbol{a}_i^{*}
\bar{\boldsymbol{x}}\vert^{2}, i = 1,\cdots, N$ with $\boldsymbol{a}_i\in \mathbb{C}^n$ drawing from the multivariate complex standard normal distribution. As a result, we can obtain $7$ different data sets with $n=400, 600, 800, 1000, 1200, 1400, 1600$.  \par
\noindent\textbf{Parameter setup:} In both SDCAM$_r$ and ADC-siDCA, we set $\mu=\mu_0/5^{t}$ with $\mu_0=50$, $\tau = 10^{5}$ and $\varepsilon_t=\varepsilon_0/1.5^t$ with $\varepsilon_0 = 10^{-4}$. The other parameters of ADC-siDCA and SDCAM$_r$ are set the same as in Subsection \ref{sec:5.2}. For PPALM method, we set the parameters as: $\rho_0 = 0.1$, $\sigma =2$ and $\varepsilon_t=10^{-4}/1.5^t$. For ADMM method, we set the parameters as follows: $\alpha_r= 1.5\times10^{-2}$, $\alpha_s = 2.0\times 10^{-4}$, $\rho_0 = 0.1$, $\beta =\frac{1+\sqrt{5}}{2}$ and $\theta=1.2$. We define the recovery error of $\widehat{\mathbf{U}}$ as $\text{RE} = \|\mathcal{A}(\widehat{\mathbf{U}})-\boldsymbol{b}\|/\max(1,\|\boldsymbol{b}\|)$.
In addition, we define the relatively phase recovery error of the vector $\hat{\boldsymbol{x}}$ as $\text{RPRE} = \|\hat{\boldsymbol{x}}-\bar{\boldsymbol{x}}\|/{\max(1,\|\bar{\boldsymbol{x}}\|)}$
where $\hat{\boldsymbol{x}} = e^{2\pi i\hat{\boldsymbol{z}}}\boldsymbol{x}'$, $\hat{\boldsymbol{z}}=\mathop{\arg}\mathop{\min}_{z\in[0,1]} \|e^{2\pi iz}\boldsymbol{x}'-\bar{\boldsymbol{x}}\|$, $\boldsymbol{x}'$ is obtained by performing spectral decomposition on $\widehat{\mathbf{U}}$. \par
\noindent\textbf{Termination criteria:} We set the termination conditions for  SDCAM$_r$ and ADC-siDCA as $\max(\textbf{Vio}_r^t,\textbf{Vio}_s^t)\leq 10^{-9}$ or $\mu_t\leq 10^{-10}$. For PPALM method, we terminate it when $\max(\textbf{Vio}_r^t,\textbf{Vio}_s^t)\leq 10^{-9}$ or $\rho_t>10^{10}$ holds. For ADMM method, we set the termination criteria as $\max(\textbf{Pinf}_{k+1},\textbf{Dinf}_{k+1})\leq 10^{-4}$.\par
\begin{figure}[H]
\centering 
\subfigure[RE]{
\includegraphics[width=0.65\linewidth]{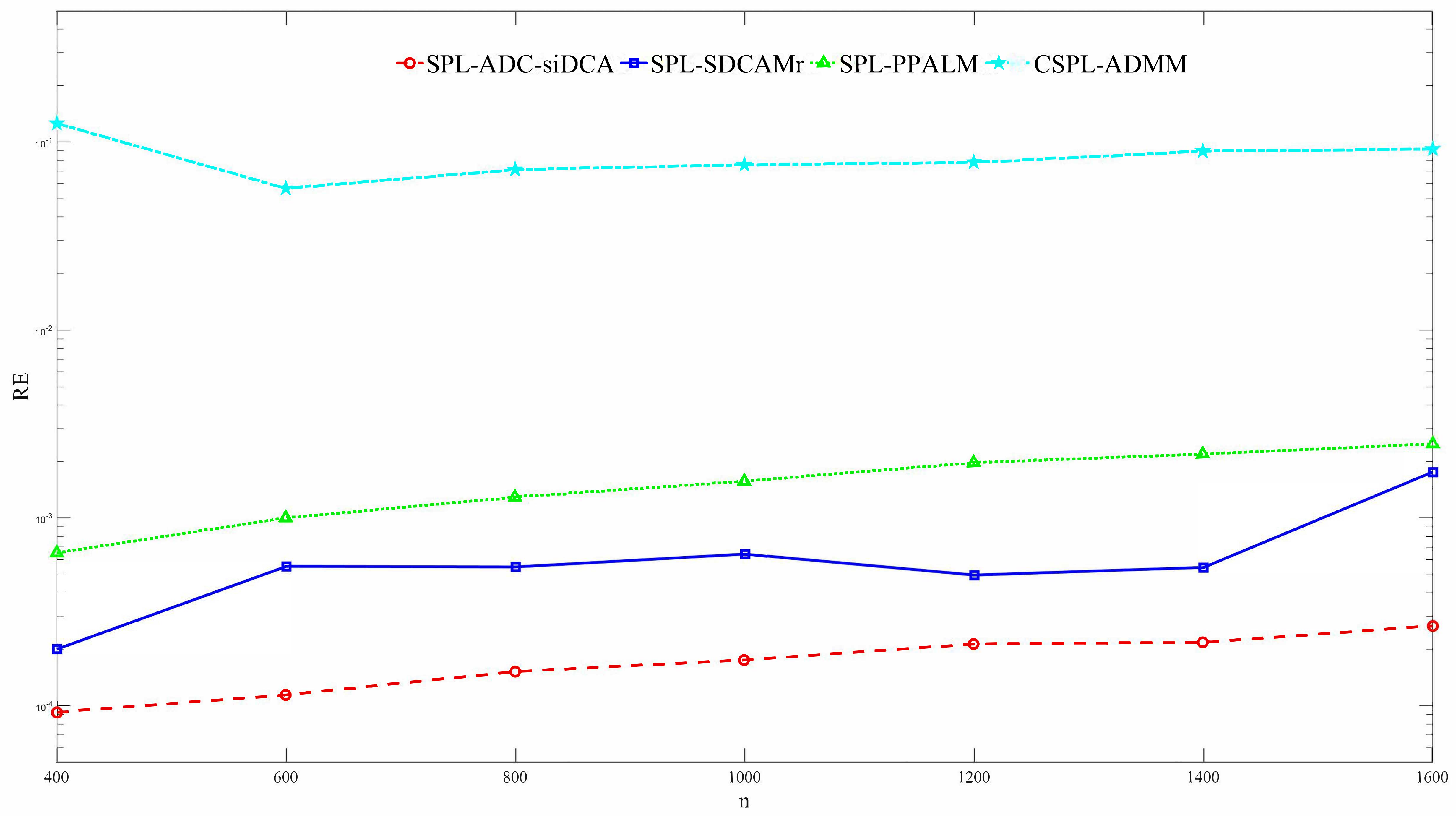}
}

\subfigure[RPRE]{
\includegraphics[width=0.65\linewidth]{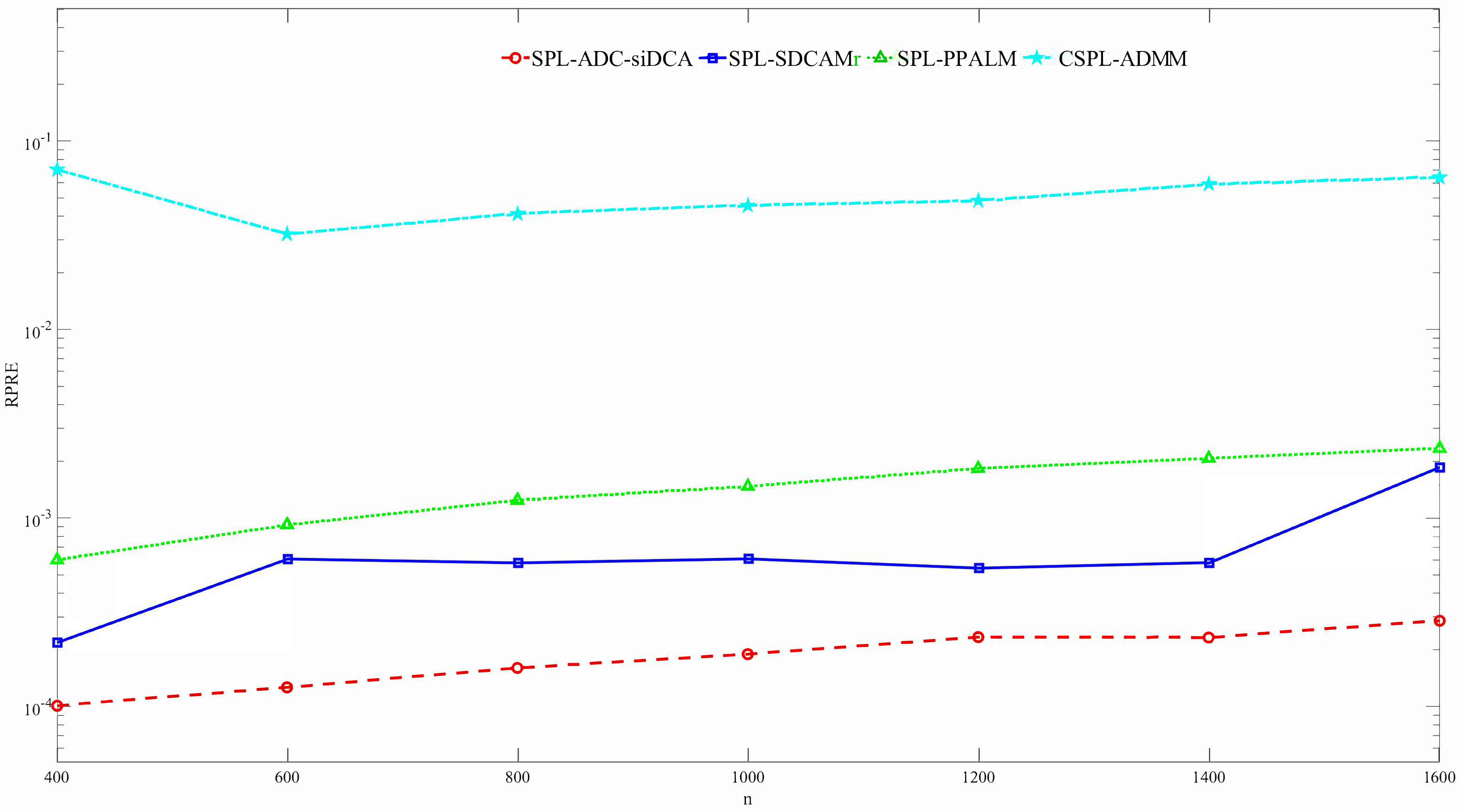}
}
 
\subfigure[Time/s]{
\includegraphics[width=0.65\linewidth]{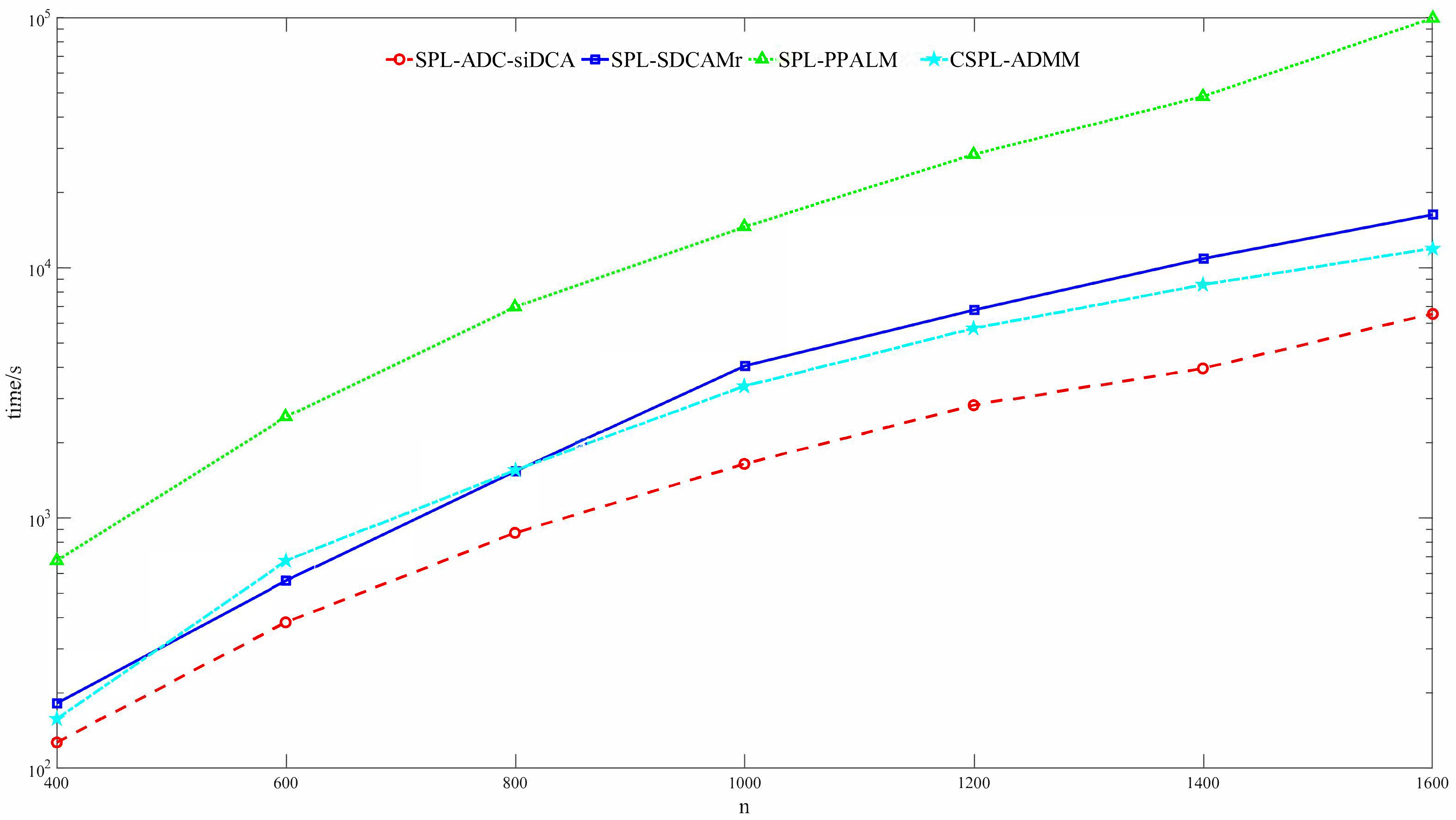}
}
\caption{ Comparison of the performance for sparse phase retrieval}\label{fig:7}
\end{figure}

In order to intuitively illustrate the numerical performance of various methods for solving sparse phase recovery, we present the change curves of RE, RPRE and CPU time (Time/s) of ADC-siDCA, SDCAM$_r$, PPALM method and ADMM method in Figure \ref{fig:7}. As shown in Figure \ref{fig:7}, with the increase of $n$, RE, RPRE and CPU time of these algorithms increase gradually. In addition, ADC-siDCA outperforms other algorithms in all aspects.\par
As a comparison, numerical results obtained by ADC-siDCA, SDCAM$_r$, PPALM method and ADMM method are presented in Table \ref{tab:3}. As a result from Table \ref{tab:3}, one can see that our ADC-siDCA outperforms the SDCAM$_r$ and PPALM method for solving SPL model. The recovery performance of SPL model is better than that of CSPL model in term of RE and RPRE. The RE and RPRE obtained by ADC-iDCA are smaller than those of other methods. Moreover, it takes much less solving time for ADC-siDCA than other methods.\par

\begin{sidewaystable}\setlength\tabcolsep{2pt}
\caption{ The performance of SPL-siDCA, SPL-SDCAM$_r$, SPL-PPALM and CSPL-ADMM for sparse phase retrieval}\centering
\label{tab:3}    
\begin{center}
\begin{minipage}{\textheight}

\begin{tabular}{cccccccccccccccccccccc}
\hline\noalign{\smallskip}
\multirow{2}{*}{n}& \multicolumn{4}{c}{SPL-ADC-siDCA}&\multicolumn{4}{c}{SPL-SDCAM$_r$}&\multicolumn{4}{c}{SPL-PPALM}&\multicolumn{4}{c}{CSPL-ADMM}\\
\cmidrule(lr){2-5}\cmidrule(lr){6-9}\cmidrule(lr){10-13}\cmidrule(lr){14-17}
&RE&RPRE &Spa&time/s&RE&RPRE &Spa&time/s&RE&RPRE &Spa&time/s&RE&RPRE &Spa&time/s \\[5pt]
\hline\noalign{\smallskip}
400 &\textbf{9.23e-5}&\textbf{1.01e-4}&20& \textbf{126.95}& 2.01e-4&2.19e-4& 20&181.91&6.52e-4&5.99e-4&20&676.05&1.26e-1&7.06e-2&20&157.63\\[15 pt]
600 &\textbf{1.14e-4}&\textbf{1.26e-4}&30&\textbf{381.91}& 5.53e-4&6.06e-4& 30&560.86&1.00e-3&9.17e-4&30&2533.50&5.67e-2&3.20e-2&30&673.46\\[15 pt]
800 &\textbf{1.52e-4}&\textbf{1.60e-4}&40& \textbf{870.34}& 5.49e-4&5.78e-4& 40&1538.97 &1.29e-3&1.24e-3&40&6977.66&7.14e-2&4.12e-2&40&1554.49\\[15 pt]
1000 &\textbf{1.75e-4}&\textbf{1.89e-4}&50& \textbf{1642.00}& 6.07e-4&6.43e-4& 50&4057.10&1.57e-3&1.47e-3&50&14569.75&7.58e-2&4.56e-2&50&3352.27\\[15 pt]
1200 &\textbf{2.13e-4}&\textbf{2.33e-4}&60& \textbf{2818.84}& 4.97e-4&5.42e-4& 60&6783.68 &1.97e-3&1.83e-3&60&28315.20&7.78e-2&4.81e-2&60&5716.29\\[15 pt]
1400 &\textbf{2.17e-4}&\textbf{2.32e-4}&70& \textbf{3951.30}& 5.46e-4&5.80e-4& 70&10878.55&2.19e-3&2.07e-3&70&48312.65&8.95e-2&5.89e-2&70&8543.38\\[15pt]
1600 &\textbf{2.66e-4}&\textbf{2.85e-4}&80& \textbf{6539.44}& 1.76e-3&1.86e-3& 80&16287.34&2.48e-3&2.34e-3&80&99268.64&9.21e-2&6.42e-2&80&11910.11\\[5pt]
\hline\noalign{\smallskip}
\end{tabular}
\end{minipage}
\footnotetext[*]{The marked result is the one with the best numerical performance in the same experimental group (excluding the sparsity of the matrix).}
\end{center}
\end{sidewaystable}

\section{Conclusion}\label{sec:6}
In this paper, the optimization problem of a sparse and low-rank matrix recovery model was considered, which entails dealing with a least squares problem along with both rank and cardinality constraints. To tackle the complexities arising from these constraints,  a Moreau smoothing method and an exact penalty method were adopted to transform the original problem into a difference-of-convex (DC) programming. Furthermore,  an asymptotic strategy was employed for updating the smoothing and penalty parameters, leading to the development of an asymptotic DC approximation (ADC) method. To solve this DC programming, we propose an efficient inexact DC algorithm with a sieving strategy (siDCA). A dual-based semismooth Newton method was used to solve the  subproblem of siDCA. We also prove the global convergence of the solution sequence generated by siDCA. To exhibit the efficacy of ADC-siDCA, we have compared it with two other methods, namely successive DC approximation minimization and penalty proximal alternating linearized minimization, by conducting matrix recovery experiments on both nonnegative matrices and positive semidefinite matrices. Our numerical results indicate that ADC-siDCA outperforms these methods in terms of recovery error and efficiency. Additionally, we conducted sparse phase retrieval experiments to illustrate the effectiveness of ADC-siDCA in Hermitian matrix recovery.

\appendix

\section{Computing Newton direction}\label{sec:3.3}
It should be noted that the main computation of SSN method is in solving the linear equations \eqref{eq.57}, then we need to given an efficient method to solve \eqref{eq.57}. We shall distinguish the following two cases to solve \eqref{eq.57}: $\mathbb{U}_+= \mathbb{S}_+^n$ and $\mathbb{U}_+= \mathbb{R}_+^{m\times n}$. Notice that when $\mathbb{U}_+= \mathbb{H}_+^n$, the Newton direction can be computed in same way as $\mathbb{U}_+= \mathbb{S}_+^n$.  \par
(1) When $\mathbb{U}_+= \mathbb{S}_+^n$, it holds that $\operatorname{Prox}_{\mu \delta_{\mathbb{S}_+^n}}(\mu\boldsymbol{t}^k(\boldsymbol{z}^j )) = \Pi_{\mathbb{S}_+^n}(\mu\boldsymbol{t}^k(\boldsymbol{z}^j ))$. Let $\lambda_1\geq\lambda_2\geq\cdots\geq\lambda_n$ be the eigenvalues of $\mu\boldsymbol{t}^k(\boldsymbol{z}^j)$ being arranged in nonincreasing order. Denote $\boldsymbol{\alpha}:=\left\{i\vert\lambda_i>0,i = 1,\cdots,n\right\}$ and $\boldsymbol{\beta}:=\left\{i\vert\lambda_i<0,i = 1,\cdots,n\right\}$. Then $\mu\boldsymbol{t}^k(\boldsymbol{z}^j)$ has the following spectral decomposition
\begin{equation}\nonumber
\mu\boldsymbol{t}^k(\boldsymbol{z}^j) = \mathbf{Q}\Lambda\mathbf{Q}^{\top},\quad \Lambda 
= \left[\begin{array}{ccc}
\Lambda_{\boldsymbol{\alpha}} & \mathbf{0} & \mathbf{0}\\
\mathbf{0} & \mathbf{0} & \mathbf{0}\\
\mathbf{0} & \mathbf{0} & \Lambda_{\boldsymbol{\beta}}
\end{array}\right],
\end{equation}
where $\Lambda$ is the diagonal matrix whose $i$-th diagonal entry is $\lambda_i$, the $i$-th column of $\mathbf{Q}$ is the eigenvector of $\mu\boldsymbol{t}^k(\boldsymbol{z}^j )$ corresponding to the eigenvalue $\lambda_i$.
Define an operator $\mathbb{W}_{\boldsymbol{z}}:\mathbb{S}^n\rightarrow \mathbb{S}^n$ as 
\[\mathbb{W}_{\boldsymbol{z}}(\mathbf{H}):=\mathbf{Q}\left(\mathbf{T} \circ\left(\mathbf{Q}^{\top} \mathbf{H} \mathbf{Q}\right)\right) \mathbf{Q}^{\top}, \; \mathbf{H} \in \mathbb{S}^n,\]
where ``$\circ$'' denotes the Hadamard product of two matrices, the matrix $\mathbf{T}$ is defined as 
\[\mathbf{T}:=\left[\begin{array}{cc}
\mathbf{E}_{\boldsymbol{\alpha} \boldsymbol{\alpha}} & \nu_{\boldsymbol{\alpha} \bar{\boldsymbol{\alpha}}} \\
\nu_{\boldsymbol{\alpha} \bar{\boldsymbol{\alpha}}}^{\top} & \mathbf{0}
\end{array}\right], \; \nu_{i j}:=\frac{\lambda_{i}}{\lambda_{i}-\lambda_{j}},\; i \in \boldsymbol{\alpha},\; j \in \bar{\boldsymbol{\alpha}}.\]
By Pang et al. \cite[Lemma 11]{Ref_pang2003semismooth}, it holds that
\[\mathbb{W}_{\boldsymbol{z}}\in\partial\Pi_{\mathbb{S}_+^n}(\mu\boldsymbol{t}^k(\boldsymbol{z})).\]
Consequently, \[\mathbb{I}+\mu\mathcal{A}\mathbb{W}_{\boldsymbol{z}}\mathcal{A}^*\in \widehat{\partial}^2\Theta_{\mu,c}^k(\boldsymbol{z}),\]
where $\mathbb{I}$ denotes identity operator.
From the definition of $\mathbb{W}_{\boldsymbol{z}}$, it follows that the matrix $\mathbf{S}^j$ defined in \eqref{eq.57} can be computed by
\[\mathbf{S}^j=\mu \mathbf{A} \mathbf{Q} \otimes \mathbf{Q} \operatorname{diag}(\mathbf{T}) \mathbf{Q}^{\top} \otimes \mathbf{Q}^{\top} \mathbf{A}^{\top}.\]
Since we solve \eqref{eq.57} by using a preconditioned conjugate gradients (PCG) method when $\mathbb{U}_+= \mathbb{S}_+^n$, then it is crucial to efficiently compute the matrix-vector product $\mathbf{S}^j\boldsymbol{d}_z^j$. Notice that $\mathbf{S}^j\boldsymbol{d}_z^j$ can be computed by $\mathbf{S}^j\boldsymbol{d}_z^j = \mathcal{A}(\mathbf{R})$, where $\mathbf{R}$ is defined as
\[\mathbf{R}:=\mathbf{Q}\left(\mathbf{T} \circ\left(\mathbf{Q}^{\top} \mathbf{D} \mathbf{Q}\right)\right) \mathbf{Q}^{\top}\]
with $\mathbf{D}=\mathcal{A}^{*}(\boldsymbol{d}_z^j)$. Thus efficient computation of $\mathbf{S}^j\boldsymbol{d}_z^j$
relies on the ability to efficiently compute the matrix $\mathbf{R}$. By noting that 
\begin{equation}\label{eq67}
\begin{aligned}
 \mathbf{R}=\left[\mathbf{Q}_{\boldsymbol{\alpha}}\; \mathbf{Q}_{\bar{\boldsymbol{\alpha}}}\right]\left[\begin{array}{cc}
\mathbf{Q}_{\boldsymbol{\alpha}}^{\top} \mathbf{D} \mathbf{Q}_{\boldsymbol{\alpha}} & \boldsymbol{\nu}_{\boldsymbol{\alpha} \bar{\boldsymbol{\alpha}}} \circ\left(\mathbf{Q}_{\boldsymbol{\alpha}}^{\top} \mathbf{D} \mathbf{Q}_{\bar{\boldsymbol{\alpha}}}\right) \\
\boldsymbol{\nu}_{\boldsymbol{\alpha} \bar{\boldsymbol{\alpha}}}^{\top} \circ\left(\mathbf{Q}_{\bar{\boldsymbol{\alpha}}}^{\top} \mathbf{D} \mathbf{Q}_{\boldsymbol{\alpha}}\right) & \mathbf{0}
\end{array}\right]\left[\begin{array}{c}
\mathbf{Q}_{\boldsymbol{\alpha}}^{\top} \\
\mathbf{Q}_{\bar{\boldsymbol{\alpha}}}^{\top}
\end{array}\right]=\mathbf{H}+\mathbf{H}^{\top}
,
\end{aligned}
\end{equation}
where $\mathbf{H}=\mathbf{Q}_{\boldsymbol{\alpha}}\left[\frac{1}{2}\left(\mathbf{U} \mathbf{Q}_{\boldsymbol{\alpha}}\right) \mathbf{Q}_{\boldsymbol{\alpha}}^{\top}+\left(\nu_{\boldsymbol{\alpha} \bar{\boldsymbol{\alpha}}} \circ\left(\mathbf{U} \mathbf{Q}_{\bar{\boldsymbol{\alpha}}}\right)\right) \mathbf{Q}_{\bar{\boldsymbol{\alpha}}}^{\top}\right]$ with $\mathbf{U} = \mathbf{Q}_{\boldsymbol{\alpha}}^{\top}\mathbf{D}$, it is easy to see
that $\mathbf{R}$ can be computed in at most $8\vert\boldsymbol{\alpha}\vert n^2$ floating point
operations (flops). On the other hand, it holds that
\begin{equation}\label{eq68}
  \mathbf{R}=\mathbf{D} - \mathbf{Q}\left((\mathbf{E}-\mathbf{T}) \circ\left(\mathbf{Q}^{\top} \mathbf{D} \mathbf{Q}\right)\right) \mathbf{Q}^{\top}
\end{equation}
and $\mathbf{R}$ can be computed in at most $8\vert\bar{\boldsymbol{\alpha}}\vert n^2$ flops. Consequently, by using the expressions \eqref{eq67} and \eqref{eq68}, we can compute $\mathbf{R}$ efficiently whenever $\vert\bar{\boldsymbol{\alpha}}\vert$ or $\vert\boldsymbol{\alpha}\vert$ is small.
\par
(2) When $\mathbb{U}_+ = \mathbb{R}_+^{m\times n}$,  $\Pi_{\mathbb{R}_+^{m\times n}}(\mu\boldsymbol{t}^k(\boldsymbol{z}^j ))$ can be computed by 
\[\Pi_{\mathbb{R}_+^{m\times n}}(\mu\boldsymbol{t}^k(\boldsymbol{z}^j )) = \max(\mu\boldsymbol{t}^k(\boldsymbol{z}^j) ,0).\]
Define an operator $\mathbb{V}_{\boldsymbol{z}}:\mathbb{R}^{m\times n}\rightarrow \mathbb{R}^{m\times n}$ as 
\[\mathbb{V}_{\boldsymbol{z}}(\mathbf{H}):=\mathbf{T} \circ \mathbf{H}, \quad \mathbf{H} \in \mathbb{R}^{m\times n},\]
where $\mathbf{T}$ is defined as 
\begin{equation}\nonumber
  \mathbf{T}_{i,j} = \left\{\begin{array}{ll}
    1, & \text{if}\;(\mu\boldsymbol{t}^k(\boldsymbol{z}))_{i,j}> 0,\\
    0, & \text{otherwise},
  \end{array}\right. 1\leq i\leq m , 1\leq j\leq n. 
\end{equation}
It is easy to see that $\mathbb{V}_{\boldsymbol{z}}\in\partial\Pi_{\mathbb{R}_+^{m\times n}}(\mu\boldsymbol{t}^k(\boldsymbol{z}^j ))$. On one hand, when PCG method is used to solve \eqref{eq.57}, the matrix-vector product $\mathbf{S}^j\boldsymbol{d}_z^j$ can be computed as follows:
\[\mathbf{S}^j\boldsymbol{d}_z^j =\mu\mathcal{A}\mathbb{V}_{\boldsymbol{z}}\mathcal{A}^*(\boldsymbol{d}_z^j)= \mu\mathcal{A}(\mathbf{T}^j\circ\mathcal{A}^*(\boldsymbol{d}_z^j)).\]
On the another hand, we can solve \eqref{eq.57} by using a direct method. Notice that the matrix $\mathbf{S}^j$ can be computed by 
\[\mathbf{S}^j = \mu\mathbf{A}\Sigma^j\mathbf{A}^{\top},\]
where $\Sigma^j:=\operatorname{diag}(\mathbf{T}^j)$. Evidently, $\Sigma^j$ is a diagonal matrix with element $1$ or $0$. Let $\mathbb{I}^j$ be the column index set of nonzero entries of $\Sigma^j$. Then it holds that
\[\mathbf{S}^j = \mu\mathbf{A}\Sigma^j\mathbf{A}^{\top}= \mu\mathbf{A}\Sigma^j\Sigma^j\mathbf{A}^{\top} = \mu\mathbf{A}\Sigma^j(\mathbf{A}\Sigma^j)^{\top} = \mu\mathbf{A}_{\mathbb{I}^j}\mathbf{A}_{\mathbb{I}^j}^{\top}.\]
Thus, we can compute the inverse of $\mathbf{I}_{N}+\mathbf{S}^j$ by the SMW formulation as follows:
\[(\mathbf{I}_{N}+\mathbf{S}^j)^{-1} = (\mathbf{I}_{N}+\mu\mathbf{A}_{\mathbb{I}^j}\mathbf{A}_{\mathbb{I}^j}^{\top})^{-1} = \mathbf{I}_{N} - \mathbf{A}_{\mathbb{I}^j}(\mu^{-1}\mathbf{I}_{\vert\mathbb{I}^j\vert}+\mathbf{A}_{\mathbb{I}^j}^{\top}\mathbf{A}_{\mathbb{I}^j})^{-1}\mathbf{A}_{\mathbb{I}^j}^{\top}.\]
As a result, \eqref{eq.57} can be solved by efficient matrix-vector multiplication as follows:
\[\boldsymbol{d}_z^j=-(\mathbf{I}_{N} - \mathbf{A}_{\mathbb{I}^j}(\mu^{-1}\mathbf{I}_{\vert\mathbb{I}^j\vert}+\mathbf{A}_{\mathbb{I}^j}^{\top}\mathbf{A}_{\mathbb{I}^j})^{-1}\mathbf{A}_{\mathbb{I}^j}^{\top})\nabla\Theta_{\mu,c}^k(\boldsymbol{z}^j).\]
It should be noted that when $\vert\mathbb{I}^j\vert\ll N$, the direct method discussed above can be more efficient than the PCG method for solving \eqref{eq.57}.

 \section{An efficient strategy for estimating inexact term} \label{appendix_A}
As shown in Algorithm \ref{alg_siDCA}, we need to check the inexact condition $\|\Delta^{k+1}\|_F\leq \epsilon_{k+1}$ and the sieving condition in \eqref{eq.48} in each iteration, then we need to compute or estimate $\Delta^{k+1}$. 
Let
\begin{equation}\nonumber
F_{\mu,c}^k(\mathbf{U}) := \frac{1}{2}\|\mathcal{A}(\mathbf{U})-\boldsymbol{b}\|^2+\frac{1}{2\mu}\|\mathbf{U}\|_F^2-\langle \Phi_c^{k},\mathbf{U}\rangle
\end{equation}
As a result, it holds that $G^k_{\mu,c}(\mathbf{U}) = F_{\mu,c}^k(\mathbf{U})+\delta_{\mathbb{U}_+}(\mathbf{U})$ and $F_{\mu,c}^k(\mathbf{U})$ is the smooth part with gradient 
\[\nabla F_{\mu,c}^k(\mathbf{U}) = \mathcal{A}^*(\mathcal{A}(\mathbf{U})-\boldsymbol{b})-\Phi_c^k+\frac{1}{\mu}\mathbf{U}.\]
If we obtain an approximate solution $\mathbf{U}_f$ by solving $\min_{\mathbf{U}}G^k_{\mu,c}(\mathbf{U})$, then there exists a $\Delta_f$ such that \begin{equation}\label{eq.61}
  \Delta_f\in\nabla F_{\mu,c}^k(\mathbf{U}_f) +\partial \delta_{\mathbb{U}_+}(\mathbf{U}_f).
\end{equation}
 According to the Second Prox Theorem \cite[Theorem 6.39]{Ref_beck2017first}, it follows that \eqref{eq.61} is equivalent to 
 \begin{equation}\label{eq.62}
  \mathbf{U}_f = \Pi_{\mathbb{U}_+}(\mathbf{U}_f-\mu\nabla F_{\mu,c}^k(\mathbf{U}_f)+\mu\Delta_f).
\end{equation}
Since $\delta_{\mathbb{U}_+}(\mathbf{U})$ is nonsmooth, it is impossible to directly obtain $\Delta_f$ from \eqref{eq.61} and \eqref{eq.62}. To address this issue, we introduce an auxiliary variable $\widetilde{\mathbf{U}}_f$ defined as \begin{equation}\label{eq.63}
\widetilde{\mathbf{U}}_f :=  \Pi_{\mathbb{U}_+}(\mathbf{U}_f-\mu\nabla F_{\mu,c}^k(\mathbf{U}_f)).
\end{equation}
 Then \eqref{eq.63} can be rewritten as
\begin{equation}\label{eq.64}
\widetilde{\mathbf{U}}_f = \Pi_{\mathbb{U}_+}( \widetilde{\mathbf{U}}_f-\mu\nabla F_{\mu,c}^k(\widetilde{\mathbf{U}}_f)+ \mu\widetilde{\Delta}_f),
\end{equation}
where 
\begin{equation}\label{eq.65}
  \widetilde{\Delta}_f = \tfrac{1}{\mu}(\mathbf{U}_f-\widetilde{\mathbf{U}}_f)-\nabla F_{\mu,c}^k(\mathbf{U}_f)+\nabla F_{\mu,c}^k(\widetilde{\mathbf{U}}_f).
\end{equation}
It is noted that \eqref{eq.64} is equivalent to 
\begin{equation}\nonumber
  \widetilde{\Delta}_f\in\nabla F_{\mu,c}^k(\widetilde{\mathbf{U}}_f)+\partial \delta_{\mathbb{U}_+}(\widetilde{\mathbf{U}}_f),
\end{equation}
which implies that $\widetilde{\mathbf{U}}_f$ is an approximate solution of \eqref{eq.45} with $\widetilde{\Delta}_f$. Hence, for a feasible solution $\mathbf{U}_f$, the approximate solution $\widetilde{\mathbf{U}}_f$ and  $\widetilde{\Delta}_f$ can be obtained by using the expressions in \eqref{eq.63} and \eqref{eq.65}, respectively. By using this technique,
we give the following two strategies to compute or estimate $\Delta^{k+1}$.\par
\text{(1)} At the $j$-th iteration of Algorithm \ref{algo_SSN}, we obtain a feasible solution $\mathbf{U}^{(j)}$ as defined in \eqref{eq.56}, shown as $\mathbf{U}^j = \mu(\Phi_c^{k}-\mathcal{A}^*(\boldsymbol{z}^j)-\mathbf{Y}^j)$. By using the expressions in \eqref{eq.63} and \eqref{eq.65}, we compute the approximate solution $\widetilde{\mathbf{U}}^{(j)}$ and  $\widetilde{\Delta}^{(j)}\in\partial G^k_{\mu,c}(\widetilde{\mathbf{U}}^{(j)})$.
If $\|\widetilde{\Delta}^{(j)}\|\leq \epsilon_{k+1}$,
 we obtain the approximate solution $\mathbf{V}^{k+1} := \widetilde{\mathbf{U}}^{(j)}$ and  $\Delta^{k+1} := \widetilde{\Delta}^{(j)}$. Clearly, $\Delta^{k+1}\in \partial G^k_{\mu,c}(\mathbf{V}^{k+1})$ and $\|\Delta^{k+1}\|< \epsilon_{k+1}$ hold.  However, one need to compute $\widetilde{\mathbf{U}}^{(j)}$ and $\widetilde{\Delta}^{(j)}$ in each iteration of Algorithm \ref{algo_SSN}, which would add additional computation cost, especially the computation of $\operatorname{Prox}_{\mu q}$ defined in \eqref{eq.63}.\par
 
 \text{(2)} We use the KKT residual $\boldsymbol{\gamma}^j$ of the dual problem to terminate Algorithm \ref{algo_SSN} in practical numerical experiments. If an upper bound estimate of $\|\widetilde{\Delta}^{(j)}\|_F$ from $\|\boldsymbol{\gamma}^j\|$ according to the relationship between $\widetilde{\Delta}^{(j)}$ and $\boldsymbol{\gamma}^j$ is obtained, we only need to check whether it is smaller than the given inexact bound $\epsilon_{k+1}$, instead of computing $\|\widetilde{\Delta}^{(j)}\|_F$ at each step of Algorithm \ref{algo_SSN}. This
will further improve the efficiency of the whole algorithm framework.\par
 From $\mathbf{U}^{(j)} = \Pi_{\mathbb{U}_+}(\mu\boldsymbol{t}^k(\boldsymbol{z}^j))$ and the definition of $\boldsymbol{\gamma}^j$ in \eqref{eq.60}, it holds that
\begin{equation}\label{eq.66}
\mathcal{A}(\mathbf{U}^{(j)}) -\boldsymbol{b}=\boldsymbol{z}^j-\boldsymbol{\gamma}^j.
\end{equation}
Since  $\widetilde{\mathbf{U}}^{(j)}$ is computed by using the expression in \eqref{eq.63}, it holds that
\begin{equation}\nonumber
\begin{aligned}
 \widetilde{\mathbf{U}}^{(j)} 
 &= \Pi_{\mathbb{U}_+}(\mathbf{U}^{(j)}-\mu\nabla F_{\mu,c}^k(\mathbf{U}^{(j)}))\\
 & =\Pi_{\mathbb{U}_+}(\mathbf{U}^{(j)}-\mu\mathcal{A}^*(\mathcal{A}(\mathbf{U}^{(j)}) -\boldsymbol{b}) +\mu\Phi_c^{k}-\mathbf{U}^{(j)}))\\
&= \Pi_{\mathbb{U}_+}(\mathbf{U}^{(j)}+\mu\mathcal{A}^*(\boldsymbol{\gamma}^j-\boldsymbol{z}^j) +\mu\Phi_c^{k}-\mu(\Phi_c^{k}-\mathcal{A}^*(\boldsymbol{z}^j)-\mathbf{Y}^j))\\
&= \Pi_{\mathbb{U}_+}(\mathbf{U}^{(j)}+\mu\mathcal{A}^*(\boldsymbol{\gamma}^j) +\mu\mathbf{Y}^j),
\end{aligned}
\end{equation}
 where the third equality follows from \eqref{eq.66}. Then a $\widetilde{\Delta}^{(j)}$ can be computed by
\begin{equation}\nonumber
\begin{aligned}
 \widetilde{\Delta}^{(j)} = \tfrac{1}{\mu}( \mathbf{U}^{(j)}-\widetilde{\mathbf{U}}^{(j)})-\nabla F_{\mu,c}^k(\mathbf{U}^{(j)})+\nabla F_{\mu,c}^k(\widetilde{\mathbf{U}}^{(j)}) 
 =\mathcal{A}^*\mathcal{A}(\widetilde{\mathbf{U}}^{(j)}-\mathbf{U}^{(j)}).
\end{aligned}
\end{equation}
From the non-expansiveness of the proximal operator, it follows that
\begin{equation}\nonumber
\begin{aligned}
&\|\widetilde{\mathbf{U}}^{(j)}-\mathbf{U}^{(j)}\|_F\\
 =& \|\Pi_{\mathbb{U}_+}(\mathbf{U}^{(j)}+\mu\mathcal{A}^*(\boldsymbol{\gamma}^j) +\mu\mathbf{Y}^j)-\mathbf{U}^{(j)}\|_F\\
=&\|\Pi_{\mathbb{U}_+}(\mathbf{U}^{(j)}+\mu\mathcal{A}^*(\boldsymbol{\gamma}^j) +\mu\mathbf{Y}^j)-\Pi_{\mathbb{U}_+}(\mathbf{U}^{(j)}+\mu\mathbf{Y}^j)\|_F\\
\leq& \mu\|\mathcal{A}^*(\boldsymbol{\gamma}^j)\|_F,
\end{aligned}
\end{equation}
where the third equality holds from \eqref{eq.59}. Consequently,
\begin{equation}\nonumber
\begin{aligned}
 \|\widetilde{\Delta}^{(j)}\|_F & =\| \mathcal{A}^*\mathcal{A}(\widetilde{\mathbf{U}}^{(j)}-\mathbf{U}^{(j)})\|_F\leq\mu\| \mathbf{A}^{\top}\mathbf{A}\|_F\|\mathcal{A}^*(\boldsymbol{\gamma}^j)\|_F\leq \mu\|\mathbf{A}\|_F\|\mathbf{A}^{\top}\mathbf{A}\|_F\|\boldsymbol{\gamma}^j\| = \frac{\|\boldsymbol{\gamma}^j\|\epsilon_{k+1}}{\zeta_{k+1}},
\end{aligned}
\end{equation}
where $\zeta_{k+1}$ is defined as 
\begin{equation}\nonumber
\zeta_k := \mu^{-1}\|\mathbf{A}\|_F^{-1}\|\mathbf{A}^{\top}\mathbf{A}\|_F^{-1}\epsilon_{k+1} .
\end{equation} 
As a result, we only need to check the condition $\|\boldsymbol{\gamma}^j\|< \zeta_{k+1}$, instead of $\|\widetilde{\Delta}^{(j)}\|_F<\epsilon_{k+1}$ in each step of Algorithm \ref{algo_SSN}. When $\|\boldsymbol{\gamma}^j\|< \zeta_{k+1}$ holds, that means that $\|\widetilde{\Delta}^{(j)}\|_F<\epsilon_{k+1}$, then we set $\mathbf{V}^{k+1} := \widetilde{\mathbf{U}}^{(j)}$ and $\Delta^{k+1} := \widetilde{\Delta}^{(j)}$. In this strategy, the calculations in \eqref{eq.63} and \eqref{eq.64} need to be performed only once, which saves a lot of computation. As a result, it is more efficient than the first strategy in practical numerical experiments, although it may overestimate the inexactness of $\widetilde{\mathbf{U}}^{(j)}$.\par

\section{Proof of convergence of siDCA}\label{appendix_b}
 \subsection*{Proof of Proposition \ref{prop_4.1}}
\begin{proof}
For statement (1), if $\mathbf{U}^{k+1}\in\mathbb{U}_+$ is generated in null step, i.e., $\mathbf{U}^{k+1} = \mathbf{U}^{k}$, it holds immediately that $J_{\mu,c}(\mathbf{U}^{k+1})\leq J_{\mu,c}(\mathbf{U}^{k})$.
If $\mathbf{U}^{k+1}\in\mathbb{U}_+$ is generated in serious step, i.e., $\mathbf{U}^{k+1} = \mathbf{V}^{k+1}$, 
it follows that
\[\mathbf{U}^{k+1} = \mathop{\arg}\mathop{\min}\limits_{\mathbf{U}} G_{\mu,c}^k(\mathbf{U})-\langle \Delta^{k+1},\mathbf{U}\rangle.\]
Since $G_{\mu,c}^k(\mathbf{U})$ is strongly convex, then the following inequality holds:
\begin{equation}\label{eq.67}
\begin{aligned}
 &f_{\mu,c}(\mathbf{U}^{k+1})-\langle \mathbf{U}^{k+1},\mathbf{W}^{k}\rangle-\langle\mathbf{U}^{k+1},\Delta^{k+1}\rangle\\
\leq &f_{\mu,c}(\mathbf{U}^{k})-\langle\mathbf{U}^{k},\mathbf{W}^{k}\rangle -\langle\mathbf{U}^{k},\Delta^{k+1}\rangle-\frac{1}{2\mu}\|\mathbf{U}^{k+1}-\mathbf{U}^{k}\|_F^2.
\end{aligned}
\end{equation}
Thanks to the convexity of $h_{\mu,c}(\mathbf{U})$, we have 
\[h_{\mu,c}(\mathbf{U}^{k+1})\geq h_{\mu,c}(\mathbf{U}^{k}) + \langle\mathbf{U}^{k+1}-\mathbf{U}^{k},\mathbf{W}^{k}\rangle.\]
This, together with \eqref{eq.67}, yields that
\begin{equation}\label{eq.68}
\begin{aligned}
 &\frac{1}{2\mu}\|\mathbf{U}^{k+1}-\mathbf{U}^{k}\|_F^2-\langle\mathbf{U}^{k+1}-\mathbf{U}^{k},\Delta^{k+1}\rangle\\
\leq& \left[f_{\mu,c}(\mathbf{U}^{k})-h_{\mu,c}(\mathbf{U}^{k})\right]-\left[f_{\mu,c}(\mathbf{U}^{k+1})-h_{\mu,c}(\mathbf{U}^{k+1})\right]\\
=&J_{\mu,c}(\mathbf{U}^{k}) -J_{\mu,c}(\mathbf{U}^{k+1}).
\end{aligned}
\end{equation}
Since $\mathbf{U}^{k+1}= \mathbf{V}^{k+1}$ is generated in serious step of siDCA, then it holds that sieving condition \eqref{eq.48} holds:
\begin{equation}\nonumber 
\|\Delta^{k+1}\|_F\leq (1-\kappa)\frac{1}{2\mu}\|\mathbf{U}^{k+1}-\mathbf{U}^{k}\|_F.
\end{equation}
Consequently, 
\begin{equation}\nonumber
\begin{aligned}
 &\frac{1}{2\mu}\|\mathbf{U}^{k+1}-\mathbf{U}^{k}\|_F^2-\langle\mathbf{U}^{k+1}-\mathbf{U}^{k},\Delta^{k+1}\rangle\\
 \geq& \frac{1}{2\mu}\|\mathbf{U}^{k+1}-\mathbf{U}^{k}\|_F^2-\|\mathbf{U}^{k+1}-\mathbf{U}^{k}\|_F\|\Delta^{k+1}\|_F\\
 \geq &\frac{\kappa}{2\mu}\|\mathbf{U}^{k+1}-\mathbf{U}^{k}\|_F^2.
\end{aligned}
\end{equation}
By applying this to \eqref{eq.68}, it follows that
\begin{equation}\label{eq.69}
\frac{\kappa}{2\mu}\|\mathbf{U}^{k+1}-\mathbf{U}^{k}\|_F^2\leq J_{\mu,c}(\mathbf{U}^{k}) -J_{\mu,c}(\mathbf{U}^{k+1})
\end{equation}
and the sequence $\{J_{\mu,c}(\mathbf{U}^{k})\}$ is non-increasing. \par
 For statement (2), according to statement (1), we have $J_{\mu,c}(\mathbf{U}^k)\leq J_{\mu,c}(\mathbf{U}^0)$, i.e., $\{J_{\mu,c}(\mathbf{U}^k)\}$ is bounded. This, together with the level-boundedness of $J_{\mu,c}$ from Theorem \ref{theorem_7}, yields that $\{\mathbf{U}^{k}\}$ is bounded. This completes the proof.
\end{proof}

 \subsection*{Proof of Theorem \ref{theo_4.1}}
\begin{proof}
For statement (1), 
since $\Delta^{\bar{k}+1}=\mathbf{0}$, then we obtain that $\mathbf{V}^{\bar{k}+1}$ is an optimal solution of \eqref{eq.45}. As a result, it follows from $\mathbf{V}^{\bar{k}+1}= \mathbf{U}^{\bar{k}}$ that $\mathbf{0} \in \partial G_{\mu,c}^{\bar{k}}(\mathbf{U}^{\bar{k}})$, i.e.,
\[\mathbf{0} \in \nabla f_{\mu,c}(\mathbf{U}^{\bar{k}}) +\partial \delta_{\mathbb{U}_+}(\mathbf{U}^{\bar{k}})-\mathbf{W}^{\bar{k}}.\] 
This, together with $\mathbf{W}^{\bar{k}}\in\partial h_{\mu,c}(\mathbf{U}^{\bar{k}})$, implies that $\mathbf{U}^{\bar{k}}$ is a stationary point of \eqref{eq.44}.

For statement (2), we first prove $\lim_{k\rightarrow \infty}\mathbf{V}^{k+1} = \mathbf{U}^{\hat{k}+1}$. Since the sieving condition in \eqref{eq.48} does not hold for any $k>\hat{k}$, i.e., 
\begin{equation}\label{eq72}
(1-\kappa)\frac{1}{2\mu}\|\mathbf{V}^{k+1}-\mathbf{U}^{\hat{k}+1}\|_F< \|\Delta^{k+1}\|_F\leq \epsilon_{k+1}
\end{equation}
holds for any $k>\hat{k}$.
 Thanks to the monotonic descent property of $\left\{\epsilon_{k}\right\}$ and $\lim_{k\rightarrow\infty}\epsilon_{k} = 0$, by taking limit on both sides of inequality \eqref{eq72}, we have
 \begin{equation}\nonumber
\lim_{k\rightarrow \infty}(1-\kappa)\frac{1}{2\mu}\|\mathbf{V}^{k+1}-\mathbf{U}^{\hat{k}+1}\|_F = 0
\end{equation}
and $\lim_{k\rightarrow \infty}\mathbf{V}^{k+1} = \mathbf{U}^{\hat{k}+1}$. Next, we will prove that $\mathbf{U}^{\hat{k}+1}$ is a stationary point of \eqref{eq.44}. Since $\mathbf{V}^{k+1}$ is an approximate solution of \eqref{eq.45} with $\Delta^{k+1}$, it holds that for any $k > \hat{k}$, 
\[\Delta^{k+1} \in \nabla f_{\mu,c}(\mathbf{V}^{k+1}) +\partial \delta_{\mathbb{U}_+}(\mathbf{V}^{k+1})- \mathbf{W}^{\hat{k}+1}.\]
Then there exists a $\zeta^{k+1}\in\partial \delta_{\mathbb{U}_+}(\mathbf{V}^{k+1})$ such that
\begin{equation}\nonumber
 \nabla f_{\mu,c}(\mathbf{V}^{k+1})+\zeta^{k+1} - \mathbf{W}^{\hat{k}+1}-\Delta^{k+1} = 0.
\end{equation}
Evidently,
 \begin{equation}\label{eq73} 
\begin{aligned}
 \|\nabla f_{\mu,c}(\mathbf{V}^{k+1})+\zeta^{k+1} - \mathbf{W}^{\hat{k}+1}\|_F&\leq \|\Delta^{k+1}\|_F\leq\epsilon_{k+1}.
\end{aligned}
\end{equation}
From $\lim_{k\rightarrow\infty}\mathbf{V}^{k+1} = \mathbf{U}^{\hat{k}+1}$, it holds that $\left\{\mathbf{V}^{k}\right\}$ is bounded. From \eqref{eq73}, we can obtain the boundedness of $\left\{\zeta^{k}\right\}$. As a consequence of \cite[Proposition 4.1.1]{Ref_hiriart1993Convex}, there exists a subset $\mathbb{K}^{\prime}\subset\mathbb{K} = \left\{0,1,2,\cdot\cdot\cdot\right\}$ such that $\lim_{k\in\mathbb{K}^{\prime}}\zeta^{k+1} = \widehat{\zeta}\in \partial \delta_{\mathbb{U}_+}(\mathbf{U}^{\hat{k}+1}) $. Hence, it follows from $\lim_{k\in\mathbb{K}^{\prime}}\epsilon_k = 0$ that
 \begin{equation}\nonumber 
\lim_{k\in\mathbb{K}^{\prime}}\|\nabla f_{\mu,c}(\mathbf{V}^{k+1}) - \mathbf{W}^{\hat{k}+1}+\zeta^{k+1}\|_F=\|\nabla f_{\mu,c}(\mathbf{U}^{\hat{k}+1}) - \mathbf{W}^{\hat{k}+1}+\widehat{\zeta}\|_F=0,
\end{equation}
which implies that
\begin{equation} \nonumber
\mathbf{0} \in \nabla f_{\mu,c}(\mathbf{U}^{\hat{k}+1}) +\partial \delta_{\mathbb{U}_+}(\mathbf{U}^{\hat{k}+1})-\partial h_{\mu,c}(\mathbf{U}^{\hat{k}+1}) 
\end{equation}
and $\mathbf{U}^{\hat{k}+1}$ is a stationary point of \eqref{eq.44}. This completes the proof.
\end{proof}

 \subsection*{Proof of Theorem \ref{thm:1}}\par
\begin{proof}
For statement (1), since $\mathbf{U}^{k_{l}}$ and $\mathbf{U}^{k_{l+1}}$ are the stability centers generated in two adjacent serious steps,  by replacing the $\mathbf{U}^{k}$ and $\mathbf{U}^{k+1}$ in \eqref{eq.69} with $\mathbf{U}^{k_{l}}$ and $\mathbf{U}^{k_{l+1}}$, respectively, we obtain
\begin{equation}\label{eq.74}
\frac{\kappa}{2\mu}\|\mathbf{U}^{k_{l+1}}-\mathbf{U}^{k_{l}}\|_F^2\leq J_{\mu,c}(\mathbf{U}^{k_{l}}) -J_{\mu,c}(\mathbf{U}^{k_{l+1}}).
\end{equation}
From statement $(1)$ of Proposition \ref{prop_4.1}, it holds that $J_{\mu,c}(\mathbf{U})$ is lower bounded and the sequence $\left\{J_{\mu,c}(\mathbf{U}^{k_{l}})\right\}$ is nonincreasing and lower bounded. As a result, by summing both sides of \eqref{eq.74} from $l = 0$ to $\infty$, we obtain
\begin{equation}\nonumber
\begin{aligned}
 \sum_{l= 0}^{\infty}\frac{\kappa}{2\mu}\|\mathbf{U}^{k_{l+1}}-\mathbf{U}^{k_l}\|_F^2\leq J_{\mu,c}(\mathbf{U}^{0}) -\liminf_{l\rightarrow \infty} J_{\mu,c}(\mathbf{U}^{k_{l+1}})<\infty.
\end{aligned}
\end{equation}
Then 
$\lim_{l\rightarrow \infty}\|\mathbf{U}^{k_{l+1}}-\mathbf{U}^{k_l}\|_F = 0$ holds.\par
For statement $(2)$, since $\mathbf{U}^{k_{l+1}}$ is an approximate solution of \eqref{eq.45} with $\Delta^{k_{l+1}}$, then we have $\Delta^{k_{l+1}} \in \partial G_{\mu,c}^{k_l} (\mathbf{U}^{k_{l+1}})$, i.e.,
\begin{equation} \nonumber
\Delta^{k_{l+1}} \in \nabla f_{\mu,c}(\mathbf{U}^{k_{l+1}}) +\partial \delta_{\mathbb{U}_+}(\mathbf{U}^{k_{l+1}})- \mathbf{W}^{k_l}.
\end{equation}
Hence, there exists a $\zeta^{k_{l+1}}\in\partial \delta_{\mathbb{U}_+}(\mathbf{U}^{k_{l+1}})$ such that 
\[\nabla f_{\mu,c}(\mathbf{U}^{k_{l+1}}) +\zeta^{k_{l+1}}- \mathbf{W}^{k_l}-\Delta^{k_{l+1}}=0.\]
Due to $\|\Delta^{k_{l+1}}\|_F\leq \epsilon_{k_{l+1}}$, it holds that 
\begin{equation}\label{eq.75}
\|\nabla f_{\mu,c}(\mathbf{U}^{k_{l+1}}) - \mathbf{W}^{k_l}+\zeta^{k_{l+1}}\|_F\leq \epsilon_{k_{l+1}}.
\end{equation}
From the boundedness of $\left\{\mathbf{U}^{k}\right\}$, it follows that $\left\{\mathbf{U}^{k_l}\right\}$ is also bounded and there exists a subset $\mathbb{L}^{\prime} \subset \mathbb{L}=\{0,1,\ldots\}$ such that $\left\{\mathbf{U}^{k_l}\right\}_{\mathbb{L}^{\prime}}$ converges to an accumulation point $\overline{\mathbf{U}}\in\left\{\mathbf{U}^{k_l}\right\}_{\mathbb{L}}$. 
This, together with the fact that $h_{\mu,c}$ is a finite-valued convex function, implies that the subsequence $\left\{\mathbf{W}^{k_l}\right\}_{\mathbb{L}^{\prime}}$ is also bounded. By using this and \eqref{eq.75}, one can obtain that $\left\{\zeta^{k_{l}}\right\}_{\mathbb{L}^{\prime}}$ is bounded. Let $\overline{\mathbf{W}}$ and $\overline{\zeta}$ be an accumulation point of $\left\{\mathbf{W}^{k_l}\right\}_{\mathbb{L}^{\prime}}$ and $\left\{\zeta^{k_{l}}\right\}_{\mathbb{L}^{\prime}}$, respectively. As a consequence of \cite[Proposition 4.1.1]{Ref_hiriart1993Convex}, we may assume without loss of generality that there exists a subset $\mathbb{L}^{\prime \prime} \subset \mathbb{L}^{\prime}$ such that $\lim_{l\in \mathbb{L}^{\prime \prime}}\mathbf{W}^{k_l} = \overline{\mathbf{W}}\in \partial h_{\mu,c}(\overline{\mathbf{U}})$, $\lim_{l\in \mathbb{L}^{\prime \prime}}\zeta^{k_{l}} = \overline{\zeta}\in \partial \delta_{\mathbb{U}_+}(\overline{\mathbf{U}})$. 
Taking limit on the two sides of inequality in \eqref{eq.75} with $l\in \mathbb{L}^{\prime \prime}$, it follows from $\lim_{k\rightarrow\infty}\epsilon_{k} = 0$ that
 \begin{equation}\nonumber
\begin{aligned}
 \|\nabla f_{\mu,c}(\overline{\mathbf{U}}) +\overline{\zeta}- \overline{\mathbf{W}}\|_F=\lim_{l\in\mathbb{L}^{\prime\prime}}\|\nabla f_{\mu,c}(\mathbf{U}^{k_{l+1}}) - \mathbf{W}^{k_l}+\zeta^{k_{l+1}}\|_F \leq \lim_{l\in \mathbb{L}^{\prime \prime}}\epsilon_{k_{l+1}}= 0,
\end{aligned}
\end{equation}
which implies that 
\begin{equation} \nonumber
\mathbf{0} \in \nabla f_{\mu,c}(\overline{\mathbf{U}}) +\partial \delta_{\mathbb{U}_+}(\overline{\mathbf{U}})-\partial h_{\mu,c}(\overline{\mathbf{U}}) 
\end{equation}
and that any accumulation point of $\left\{\mathbf{U}^{k_l}\right\}$ is a stationary point of \eqref{eq.44}.
This completes the proof.
\end{proof}

\subsection*{Proof of Proposition \ref{prop_4.2}}\par
\begin{proof}
For statement (1), since $\mathbf{U}^{k_{l+1}} = \mathbf{V}^{k+1}\in\mathbb{U}_+$ is the stability center generated in serious step, then condition \eqref{eq.48} holds, shown as
\begin{equation}\label{eq___C12}
\|\Delta^{k_{l+1}}\|_F\leq(1-\kappa)\frac{1}{2\mu}\|\mathbf{U}^{k_{l+1}}-\mathbf{U}^{k_l}\|_F.
\end{equation}
Then we have
\begin{equation}\label{eq___C13}
\frac{\kappa}{2\mu}\|\mathbf{U}^{k_{l+1}}-\mathbf{U}^{k_l}\|_F^2\leq \frac{1}{2\mu}\|\mathbf{U}^{k_{l+1}}-\mathbf{U}^{k_l}\|_F^2-\langle\Delta^{k_{l+1}},\mathbf{U}^{k_{l+1}}-\mathbf{U}^{k_l}\rangle.
\end{equation}
 Consequently,
\begin{equation}\nonumber 
\begin{aligned}
&\Psi(\mathbf{U}^{k_{l+1}},\mathbf{W}^{k_l},\mathbf{U}^{k_l},\Delta^{k_{l+1}})\\
=& f_{\mu,c}(\mathbf{U}^{k_{l+1}})-\langle \mathbf{U}^{k_{l+1}},\mathbf{W}^{k_l}\rangle + h^*_{\mu,c}(\mathbf{W}^{k_l})\\
&+\frac{1}{2\mu}\|\mathbf{U}^{k_{l+1}} -\mathbf{U}^{k_l}\|^2_F\\
=& f_{\mu,c}(\mathbf{U}^{k_{l+1}})-\langle \mathbf{U}^{k_{l+1}}-\mathbf{U}^{k_l},\mathbf{W}^{k_l}\rangle\\
& - h_{\mu,c}(\mathbf{U}^{k_l}) +\frac{1}{2\mu}\|\mathbf{U}^{k_{l+1}} -\mathbf{U}^{k_l}\|^2_F\\
\geq& f_{\mu,c}(\mathbf{U}^{k_{l+1}})- h_{\mu,c}(\mathbf{U}^{k_{l+1}})+\frac{1}{2\mu}\|\mathbf{U}^{k_{l+1}} -\mathbf{U}^{k_l}\|^2_F\\
\geq& f_{\mu,c}(\mathbf{U}^{k_{l+1}})- h_{\mu,c}(\mathbf{U}^{k_{l+1}}) =J_{\mu,c}(\mathbf{U}^{k_{l+1}}),
\end{aligned}
\end{equation}
where the second equality follows from the convexity of $h_{\mu,c}$ and the fact that $\mathbf{W}^{k_l}\in\partial h_{\mu,c}(\mathbf{U}^{k_l})$, the first inequality follows from the convexity of $h_{\mu,c}$.\par
For statement (2), since $\mathbf{U}^{k_{l+1}}\in \mathbb{U}_+$ is an approximate solution of \eqref{eq.44} with $\Delta^{k_{l+1}}$, we have
\[\mathbf{U}^{k_{l+1}} = \mathop{\arg}\mathop{\min}\limits_{\mathbf{U}} G_{\mu,c}^{k_l}(\mathbf{U})-\langle \Delta^{k_{l+1}},\mathbf{U}\rangle.\]
Since $G_{\mu,c}^{k_l}(\mathbf{U})$ is strongly convex with parameter $\frac{1}{\mu}$, then the following inequality holds: 
\begin{equation}\label{eq___C14}
\begin{aligned}
 G_{\mu,c}^{k_l}(\mathbf{U}^{k_{l+1}})-\langle\Delta^{k_{l+1}},\mathbf{U}^{k_{l+1}}\rangle
\leq G_{\mu,c}^{k_l}(\mathbf{U}^{k_l}) -\langle\Delta^{k_{l+1}},\mathbf{U}^{k_l}\rangle-\frac{1}{2\mu}\|\mathbf{U}^{k_{l+1}}-\mathbf{U}^{k_l}\|_F^2.
\end{aligned}
\end{equation}
Evidently,
\begin{equation} \nonumber
\begin{aligned}
&\Psi(\mathbf{U}^{k_{l+1}},\mathbf{W}^{k_l},\mathbf{U}^{k_l},\Delta^{k_{l+1}})\\
=& f_{\mu,c}(\mathbf{U}^{k_{l+1}})+\delta_{\mathbb{U}_+}(\mathbf{U}^{k_{l+1}})-\langle \mathbf{U}^{k_{l+1}},\mathbf{W}^{k_l}\rangle + h^*_{\mu,c}(\mathbf{W}^{k_l})\\
&+\frac{1}{2\mu}\|\mathbf{U}^{k_{l+1}} -\mathbf{U}^{k_l}\|^2_F
-\langle\Delta^{k_{l+1}},\mathbf{U}^{k_{l+1}}-\mathbf{U}^{k_l}\rangle\\
 \leq& f_{\mu,c}(\mathbf{U}^{k_l})+\delta_{\mathbb{U}_+}(\mathbf{U}^{k_{l}})-\langle \mathbf{U}^{k_l},\mathbf{W}^{k_l}\rangle + h^*_{\mu,c}(\mathbf{W}^{k_l})\\
=& f_{\mu,c}(\mathbf{U}^{k_l})+\delta_{\mathbb{U}_+}(\mathbf{U}^{k_{l}})- h_{\mu,c}(\mathbf{U}^{k_l}),
\end{aligned}
\end{equation}
where the last equality follows from the convexity of $h_{\mu,c}$ and the fact that $\mathbf{W}^{k_l}\in\partial h_{\mu,c}(\mathbf{U}^{k_l})$. Similar to  \eqref{eq___C13}, the following inequality holds:
\begin{equation}\label{eq___C15}
\frac{\kappa}{2\mu}\|\mathbf{U}^{k_l}-\mathbf{U}^{k_{l-1}}\|_F^2\leq \frac{1}{2\mu}\|\mathbf{U}^{k_l}-\mathbf{U}^{k_{l-1}}\|_F^2-\langle\Delta^{k_{l}},\mathbf{U}^{k_l}-\mathbf{U}^{k_{l-1}}\rangle.
\end{equation}
Consequently, 
\begin{equation}\nonumber 
\begin{aligned}
&\Psi(\mathbf{U}^{k_{l+1}},\mathbf{W}^{k_l},\mathbf{U}^{k_l},\Delta^{k_{l+1}}) \\
 \leq& f_{\mu,c}(\mathbf{U}^{k_l}) +\delta_{\mathbb{U}_+}(\mathbf{U}^{k_{l}})-\langle \mathbf{U}^{k_l},\mathbf{W}^{k_{l-1}}\rangle + h^*_{\mu,c}(\mathbf{W}^{k_{l-1}})\\
 = &\Psi(\mathbf{U}^{k_l},\mathbf{W}^{k_{l-1}},\mathbf{U}^{k_{l-1}},\Delta^{k_{l}})- \frac{1}{2\mu}\|\mathbf{U}^{k_l} -\mathbf{U}^{k_{l-1}}\|^2_F+ \langle\Delta^{k_{l}},\mathbf{U}^{k_l}-\mathbf{U}^{k_{l-1}}\rangle\\
\leq& \Psi(\mathbf{U}^{k_l},\mathbf{W}^{k_{l-1}},\mathbf{U}^{k_{l-1}},\Delta^{k_{l}})- \frac{\kappa}{2\mu}\|\mathbf{U}^{k_l} -\mathbf{U}^{k_{l-1}}\|^2_F,
\end{aligned}
\end{equation}
where the second inequality follows from the convexity of $h_{\mu,c}$ and the Young’s inequality applied to $h_{\mu,c}$. The last inequality is due to \eqref{eq___C15}.\par
For statement $(3)$, it holds from Proposition \ref{prop_4.1} that $\{\mathbf{U}^{k_l}\}$ is bounded. The boundedness of $\{\mathbf{W}^{k_l}\}$ follows immediately from the fact that $h_{\mu,c}$ is finite-valued convex function and that $\mathbf{W}^{k_l}\in\partial h_{\mu,c}(\mathbf{U}^{k_l})$. The boundedness of $\{\Delta^{k_{l}} \}$ follows from the fact that $\lim_{l\rightarrow\infty}\epsilon_{k_{l+1}} = 0$. Then, the bounded sequence $\{(\mathbf{U}^{k_{l+1}}, \mathbf{W}^{k_l}, \mathbf{U}^{k_l},\Delta^{k_{l+1}})\}$ has nonempty accumulation point set $\Gamma$.\par
For statement (4), since $J_{\mu,c}$ is bounded below, from \eqref{eq.77} and \eqref{eq.78}, we have that $\left\{\Psi(\mathbf{U}^{k_{l+1}},\mathbf{W}^{k_l},\mathbf{U}^{k_l},\Delta^{k_{l+1}}) \right\}$ is nonincreasing and bounded below. Thus, the limit $\Upsilon = \lim_{l\rightarrow\infty}\Psi(\mathbf{U}^{k_{l+1}},\mathbf{W}^{k_l},\mathbf{U}^{k_l},\Delta^{k_{l+1}})$ exists. Next, we will prove that $\Psi \equiv \Upsilon$ on $\Gamma$. Take any $(\widehat{\mathbf{U}},\widehat{\mathbf{W}},\widehat{\mathbf{U}},\widehat{\Delta})\in \Gamma$. As a result, it follows that there exists a subset $\mathbb{L}^{\prime}\subset \mathbb{L}$ such that 
\[\lim_{l\in\mathbb{L}^{\prime}}(\mathbf{U}^{k_{l+1}},\mathbf{W}^{k_l},\mathbf{U}^{k_l},\Delta^{k_{l+1}}) =(\widehat{\mathbf{U}},\widehat{\mathbf{W}},\widehat{\mathbf{U}},\widehat{\Delta}).\]
From the compactness of $\mathbb{U}_+$, it holds that $\widehat{\mathbf{U}}\in\mathbb{U}_+$.
This, together with the optimality of $\mathbf{U}^{k_{l+1}}$ for solving $\min_{\mathbf{U}}G_{\mu,c}^{k_l}(\mathbf{U})$, yields that
\begin{equation}\label{eq___C16}
\begin{aligned}
  f_{\mu,c}(\mathbf{U}^{k_{l+1}})-\langle \mathbf{U}^{k_{l+1}},\mathbf{W}^{k_l}\rangle-\langle \mathbf{U}^{k_{l+1}},\Delta^{k_{l+1}}\rangle
\leq f_{\mu,c}(\widehat{\mathbf{U}})-\langle \widehat{\mathbf{U}},\mathbf{W}^{k_l}\rangle -\langle \widehat{\mathbf{U}},\Delta^{k_{l+1}}\rangle.
\end{aligned}
\end{equation}
By rearranging terms in the above inequality, we obtain 
\begin{equation}\label{eq___C17}
\begin{aligned}
 f_{\mu,c}(\mathbf{U}^{k_{l+1}})-\langle \mathbf{U}^{k_{l+1}}-\widehat{\mathbf{U}},\mathbf{W}^{k_l}+\Delta^{k_{l+1}}\rangle\leq f_{\mu,c}(\widehat{\mathbf{U}}).
\end{aligned}
\end{equation}
From the boundedness of $\{\mathbf{U}^{k_l}\}$, $\{\mathbf{W}^{k_l}\}$  and $\{\Delta^{k_{l}}\}$, we have 
\[\lim_{l\in\mathbb{L}^{\prime}}\langle \mathbf{U}^{k_{l+1}}-\widehat{\mathbf{U}},\mathbf{W}^{k_l}\rangle = 0,\quad \lim_{l\in\mathbb{L}^{\prime}}\langle \mathbf{U}^{k_{l+1}}-\widehat{\mathbf{U}},\Delta^{k_{l+1}}\rangle = 0.\]
Then, we have
\begin{equation}\nonumber
\begin{aligned}
 \Upsilon& = \lim_{l\in\mathbb{L}^{\prime}}\Psi(\mathbf{U}^{k_{l+1}},\mathbf{W}^{k_l},\mathbf{U}^{k_l},\Delta^{k_{l+1}})\\
&= \lim_{l\in\mathbb{L}^{\prime}}f_{\mu,c}(\mathbf{U}^{k_{l+1}})+\delta_{\mathbb{U}_+}(\mathbf{U}^{k_{l+1}})-\langle \mathbf{U}^{k_{l+1}},\mathbf{W}^{k_l}\rangle + h^*_{\mu,c}(\mathbf{W}^{k_l})\\
&\quad+\frac{1}{2\mu}\|\mathbf{U}^{k_{l+1}} -\mathbf{U}^{k_l}\|^2_F
-\langle\Delta^{k_{l+1}},\mathbf{U}^{k_{l+1}}-\mathbf{U}^{k_l}\rangle\\
&= \lim_{l\in\mathbb{L}^{\prime}} f_{\mu,c}(\mathbf{U}^{k_{l+1}})-\langle \mathbf{U}^{k_{l+1}}-\widehat{\mathbf{U}},\mathbf{W}^{k_l}+\Delta^{k_{l+1}}\rangle-\langle \mathbf{U}^{k_{l+1}},\mathbf{W}^{k_l}\rangle\\
&+\frac{1}{2\mu}\|\mathbf{U}^{k_{l+1}}-\mathbf{U}^{k_l}\|_F^2+ h^*_{\mu,c}(\mathbf{W}^{k_l})-\langle\Delta^{k_{l+1}},\mathbf{U}^{k_{l+1}}-\mathbf{U}^{k_l}\rangle\\
&\leq \lim\sup_{l\in\mathbb{L}^{\prime}} f_{\mu,c}(\widehat{\mathbf{U}})+\frac{1}{2\mu}\|\mathbf{U}^{k_{l+1}}-\mathbf{U}^{k_l}\|_F^2-\langle \mathbf{U}^{k_{l+1}},\mathbf{W}^{k_l}\rangle \\
&\quad+ h^*_{\mu,c}(\mathbf{W}^{k_l})-\langle\Delta^{k_{l+1}},\mathbf{U}^{k_{l+1}}-\mathbf{U}^{k_l}\rangle\\
&= \lim\sup_{l\in\mathbb{L}^{\prime}} f_{\mu,c}(\widehat{\mathbf{U}})+\frac{1}{2\mu}\|\mathbf{U}^{k_{l+1}}-\mathbf{U}^{k_l}\|_F^2- h_{\mu,c}(\mathbf{U}^{k_l}) \\
&\quad-\langle \mathbf{U}^{k_{l+1}}-\mathbf{U}^{k_l},\mathbf{W}^{k_l}\rangle
-\langle\Delta^{k_{l+1}},\mathbf{U}^{k_{l+1}}-\mathbf{U}^{k_l}\rangle\\
&= f_{\mu,c}(\widehat{\mathbf{U}})-h_{\mu,c}(\widehat{\mathbf{U}}) =J_{\mu,c}(\widehat{\mathbf{U}}) \leq \Psi(\widehat{\mathbf{U}},\widehat{\mathbf{W}},\widehat{\mathbf{U}}, \widehat{\Delta}),
\end{aligned}
\end{equation}
where the fourth equality follows from the convexity of $h_{\mu,c}$ and $\mathbf{W}^{k_l}\in\partial h_{\mu,c}(\mathbf{U}^{k_l})$, the last inequality holds from \eqref{eq.77} with $l$ trending to infinity and the lower semicontinuity of $\Psi$. Since $\Psi$ is lower semicontinuous, it holds that
 \[\Psi(\widehat{\mathbf{U}},\widehat{\mathbf{W}},\widehat{\mathbf{U}},\widehat{\Delta})= \lim\inf_{l\in\mathbb{L}^{\prime}}\Psi(\mathbf{U}^{k_{l+1}},\mathbf{W}^{k_l},\mathbf{U}^{k_l},\Delta^{k_{l+1}})=\Upsilon\]
and $\Psi \equiv \Upsilon$ on $\Gamma$.\par
 For statement (5), since the subdifferential of the function $\Psi$ at the point $(\mathbf{U}^{k_{l+1}},\mathbf{W}^{k_l},\mathbf{U}^{k_l},\Delta^{k_{l+1}})$ is 
 \begin{equation}\nonumber
\begin{aligned}
 &\partial \Psi(\mathbf{U}^{k_{l+1}},\mathbf{W}^{k_l},\mathbf{U}^{k_l},\Delta^{k_{l+1}})\\
=&\left[\begin{array}{c}
\nabla f_{\mu,c}(\mathbf{U}^{k_{l+1}})+\partial\delta_{\mathbb{U}_+}(\mathbf{U}^{k_{l+1}})- \mathbf{W}^{k_l}+\frac{1}{\mu}(\mathbf{U}^{k_{l+1}}-\mathbf{U}^{k_l})-\Delta^{k_{l+1}}\\
-\mathbf{U}^{k_{l+1}}+\partial H_{c}^{*}(\mathbf{W}^{k_l}) \\
-\frac{1}{\mu}(\mathbf{U}^{k_{l+1}}-\mathbf{U}^{k_l})+\Delta^{k_{l+1}}\\
\mathbf{U}^{k_l}-\mathbf{U}^{k_{l+1}}
\end{array}\right].
\end{aligned}
\end{equation}
Since $\mathbf{U}^{k_{l+1}}$ is the optimal solution of \eqref{eq.45}, we have
\[\mathbf{0}\in\nabla f_{\mu,c}(\mathbf{U}^{k_{l+1}})+\partial\delta_{\mathbb{U}_+}(\mathbf{U}^{k_{l+1}})- \mathbf{W}^{k_l} -\Delta^{k_{l+1}}.\] 
Since $\mathbf{U}^{k_l}\in\partial H_{c}^{*}(\mathbf{W}^{k_l})$, then
 \begin{equation}\label{eq___C18}
\left[\begin{array}{c}
\frac{1}{\mu}(\mathbf{U}^{k_{l+1}}-\mathbf{U}^{k_l})\\
\mathbf{U}^{k_l}-\mathbf{U}^{k_{l+1}} \\
-\frac{1}{\mu}(\mathbf{U}^{k_{l+1}}-\mathbf{U}^{k_l})+\Delta^{k_{l+1}}\\
\mathbf{U}^{k_l}-\mathbf{U}^{k_{l+1}}
\end{array}\right]\in\partial \Psi(\mathbf{U}^{k_{l+1}},\mathbf{W}^{k_l},\mathbf{U}^{k_l},\Delta^{k_{l+1}})
\end{equation}
Since $\mathbf{U}^{k_{l+1}}$ is the stability center of serious step, then it holds that \[\|\Delta^{k_{l+1}}\|_F\leq (1-\kappa)\frac{1}{\mu}\|\mathbf{U}^{k_{l+1}}-\mathbf{U}^{k_l}\|_F.\] Thus there exists a constant $\rho$ such that the following inequality holds:
\begin{equation}\label{eq___C19}
\operatorname{dist}(\mathbf{0},\partial \Psi(\mathbf{U}^{k_{l+1}},\mathbf{W}^{k_l},\mathbf{U}^{k_l},\Delta^{k_{l+1}}))\leq \rho\|\mathbf{U}^{k_{l+1}}-\mathbf{U}^{k_l}\|_F
\end{equation}
This completes the proof.
\end{proof}

\subsection*{Proof of Theorem \ref{thm:6}}\par
For simplicity of notation, we set $\Psi^{k_{l}} = \Psi(\mathbf{U}^{k_{l}},\mathbf{W}^{k_{l-1}},\mathbf{U}^{k_{l-1}},\Delta^{k_{l}})$ for each $l>0$. 
\begin{proof}
From Proposition \ref{prop_4.2}, we have that $\{\Psi^{k_l}\}$ is nonincreasing and its limit $\Upsilon$ exists. Then, we get $\Psi^{k_l}\geq\Upsilon$ for any $l>0$.
Next, we will prove that for any $l>0$, $\Psi^{k_l} > \Upsilon.$
To this end, we suppose that $\exists L>0$ such that $\Psi^{k_L} = \Upsilon$, then $\Psi^{k_l} = \Upsilon$ holds for all $l>L$. From \eqref{eq.78}, we have $\mathbf{U}^{k_l} =\mathbf{U}^{k_L}$ for each $l \geq L$. This implies that only finite serious steps are performed in siDCA, which is contrary to the assumption.\par
Since $\Psi$ satisfies the KŁ property at each point in the compact set $\Gamma\subset\operatorname{dom} \partial \Psi$ and $\Psi\equiv\Upsilon$ on $\Gamma$, thus it satisfies the uniform KŁ property \cite{Ref_bolte2014proximal}. Then there exist $\epsilon>0$ and a continuous concave function $\varphi:[0,a)\rightarrow \mathbb{R}_+$ being continuously differentiable and monotonically increasing on $(0,a)$ and satisfying $\varphi(0) =0$ with $a > 0$ such that for any $ (\mathbf{U},\mathbf{W},\mathbf{V},\mathbf{Z})\in \Theta$,
\begin{equation}\nonumber
\varphi^{\prime}(\Psi(\mathbf{U},\mathbf{W},\mathbf{V}, \mathbf{Z})-\Upsilon) \cdot \operatorname{dist}(\mathbf{0}, \partial \Psi(\mathbf{U},\mathbf{W},\mathbf{V},\mathbf{Z})) \geq 1,
\end{equation} where
\begin{equation}\nonumber
\begin{aligned}
\Theta=&\left\{(\mathbf{U},\mathbf{W},\mathbf{V},\mathbf{Z}): \operatorname{dist}((\mathbf{U},\mathbf{W},\mathbf{V},\mathbf{Z}), \Gamma)<\epsilon\right\} 
\cap\left\{(\mathbf{U},\mathbf{W},\mathbf{V},\mathbf{Z}): \Upsilon<\Psi(\mathbf{U},\mathbf{W},\mathbf{V},\mathbf{Z})<\Upsilon+a\right\}.
\end{aligned}
\end{equation}
Since $\Gamma$ is the set of accumulation points of $\{(\mathbf{U}^{k_{l+1}}, \mathbf{W}^{k_l}, \mathbf{U}^{k_l},\Delta^{k_{l+1}})\}$, we have
\begin{equation}\nonumber
\lim_{l\rightarrow\infty}\operatorname{dist} ((\mathbf{U}^{k_{l+1}}, \mathbf{W}^{k_l}, \mathbf{U}^{k_l},\Delta^{k_{l+1}}),\Gamma) = 0.
\end{equation}
Thus, there exists a constant $\bar{L}>0$ such that for any $l>\bar{L}-2$,
\[\operatorname{dist} ((\mathbf{U}^{k_{l+1}}, \mathbf{W}^{k_l}, \mathbf{U}^{k_l},\Delta^{k_{l+1}}),\Gamma)<\epsilon.\]
 From Proposition \ref{prop_4.2}, we have that the sequence $\{\Psi^{k_l} \}$ converges to $\Upsilon$, then exists a constant $\bar{\bar{L}}>0$ such that for any $l>\bar{\bar{L}}-2$,
 \[\Upsilon<\Psi^{k_{l+1}}<\Upsilon+a.\]
 Let $\tilde{L} = \max\left\{\bar{L},\bar{\bar{L}}\right\}$, then $\forall\,l>\tilde{L}$, we have $(\mathbf{U}^{k_{l-1}}, \mathbf{W}^{k_{l-2}}, \mathbf{U}^{k_{l-2}},\Delta^{k_{l-1}} )\in\Theta$ and 
\begin{equation}\label{eq___57}
\varphi^{\prime}(\Psi^{k_{l-1}}-\Upsilon) \cdot \operatorname{dist}(\mathbf{0}, \partial \Psi^{k_{l-1}} )\geq 1.
\end{equation}
 By using the concavity of $\varphi$, it holds that for each $l>\tilde{L}$, 
\begin{equation}\label{eq___58}
\begin{aligned}
&\left[\varphi(\Psi^{k_{l-1}}-\Upsilon) -\varphi(\Psi^{k_{l+1}}-\Upsilon)\right] \cdot \operatorname{dist}(\mathbf{0}, \partial \Psi^{k_{l-1}}) \\
\geq&\varphi^{\prime}(\Psi^{k_{l-1}}-\Upsilon)\cdot \operatorname{dist}(\mathbf{0}, \partial \Psi^{k_{l-1}})\cdot (\Psi^{k_{l-1}}-\Psi^{k_{l+1}}) \\
\geq& \Psi^{k_{l-1}}-\Psi^{k_{l+1}},
\end{aligned}
\end{equation}
where the last inequality follows from \eqref{eq___57}.
Let $\pi^{k_{l}} = \varphi(\Psi^{k_{l}}-\Upsilon)$ for each $l>0$. Since $\varphi$ is monotone increasing on $(0,a)$ and $\{\Psi^{k_l}\}$ is nonincreasing, then we get that $\{\pi^{k_{l}}\}$ is nonincreasing. Combining the results in \eqref{eq.78}, \eqref{eq.79} and \eqref{eq___58}, we have 
\begin{equation}\label{eq___59}
\begin{aligned}
\|\mathbf{U}^{k_l} -\mathbf{U}^{k_{l-1}}\|^2_F+\|\mathbf{U}^{k_{l-1}} -\mathbf{U}^{k_{l-2}}\|^2_F
\leq \frac{\rho}{\kappa\alpha}(\pi^{k_{l-1}} -\pi^{k_{l+1}})\|\mathbf{U}^{k_{l-1}} -\mathbf{U}^{k_{l-2}}\|_F.
\end{aligned}
\end{equation}
According to the above inequality, we have \[\frac{\rho}{\kappa\alpha}(\pi^{k_{l-1}} -\pi^{k_{l+1}})-\|\mathbf{U}^{k_{l-1}} -\mathbf{U}^{k_{l-2}}\|_F\geq 0.\]
Then by applying the arithmetic mean–geometric mean inequality, we obtain
\begin{equation}\nonumber
\begin{aligned}
\|\mathbf{U}^{k_l} -\mathbf{U}^{k_{l-1}}\|_F
\leq& \sqrt{\frac{\rho}{2\kappa\alpha}(\pi^{k_{l-1}} -\pi^{k_{l+1}})-\frac{1}{2}\|\mathbf{U}^{k_{l-1}} -\mathbf{U}^{k_{l-2}}\|_F}\cdot\sqrt{2\|\mathbf{U}^{k_{l-1}} -\mathbf{U}^{k_{l-2}}\|_F}\\
\leq& \frac{\rho}{4\kappa\alpha}(\pi^{k_{l-1}} -\pi^{k_{l+1}})-\frac{1}{4}\|\mathbf{U}^{k_{l-1}} -\mathbf{U}^{k_{l-2}}\|_F+ \|\mathbf{U}^{k_{l-1}} -\mathbf{U}^{k_{l-2}}\|_F.
\end{aligned}
\end{equation}
Consequently, we have
\begin{equation}\label{eq___60}
\begin{aligned}
\frac{1}{4}\|\mathbf{U}^{k_l} -\mathbf{U}^{k_{l-1}}\|_F
\leq \frac{\rho}{4\kappa\alpha}(\pi^{k_{l-1}} -\pi^{k_{l+1}})+ \frac{3}{4}(\|\mathbf{U}^{k_{l-1}} -\mathbf{U}^{k_{l-2}}\|_F-\|\mathbf{U}^{k_l} -\mathbf{U}^{k_{l-1}}\|_F).
\end{aligned}
\end{equation}
Summing both sides of \eqref{eq___60} from $l=\tilde{L}$ to $\infty$, we have
\begin{equation}\nonumber
\begin{aligned}
\frac{1}{4}\sum_{l=\tilde{L}}^{\infty}\|\mathbf{U}^{k_l} -\mathbf{U}^{k_{l-1}}\|_F\leq& \frac{\rho}{4\kappa\alpha}(\pi^{k_{\tilde{L}-1}} +\pi^{k_{\tilde{L}}})- \lim_{l\rightarrow\infty}\frac{\rho}{4\kappa\alpha}(\pi^{k_{l}} +\pi^{k_{l+1}})\\
&+ \frac{3}{4}(\|\mathbf{U}^{k_{\tilde{L}-1}} -\mathbf{U}^{k_{\tilde{L}-2}}\|_F-\lim_{l\rightarrow\infty}\|\mathbf{U}^{k_l} -\mathbf{U}^{k_{l-1}}\|_F).
\end{aligned}
\end{equation}
From $\lim_{l\rightarrow\infty}\frac{\rho}{4\kappa\alpha}(\pi^{k_{l}} +\pi^{k_{l+1}})=0$ and $\lim_{l\rightarrow\infty}\|\mathbf{U}^{k_l} -\mathbf{U}^{k_{l-1}}\|_F=0$, we obtain
\begin{equation}\nonumber
\begin{aligned}
\frac{1}{4}\sum_{l=\tilde{L}}^{\infty}\|\mathbf{U}^{k_l} -\mathbf{U}^{k_{l-1}}\|_F
\leq \frac{\rho}{4\kappa\alpha}(\pi^{k_{\tilde{L}-1}} +\pi^{k_{\tilde{L}}})+ \frac{3}{4}\|\mathbf{U}^{k_{\tilde{L}-1}} -\mathbf{U}^{k_{\tilde{L}-2}}\|_F<\infty.
\end{aligned}
\end{equation}
Thus the subsequence $\{\mathbf{U}^{k_l}\}$ is convergent and $\sum_{l=0}^{\infty}\|\mathbf{U}^{k_{l+1}}-\mathbf{U}^{k_l}\|_F<\infty$. From Theorem \ref{thm:1}, we have that the sequence $\{\mathbf{U}^{k_l}\}$ generated by siDCA converges to a stationary point of \eqref{eq.44}. This completes the proof.
\end{proof}
\section{Some figures for matrix recovery}
\begin{figure}[H]
\centering 
  \includegraphics[width=0.45\linewidth]{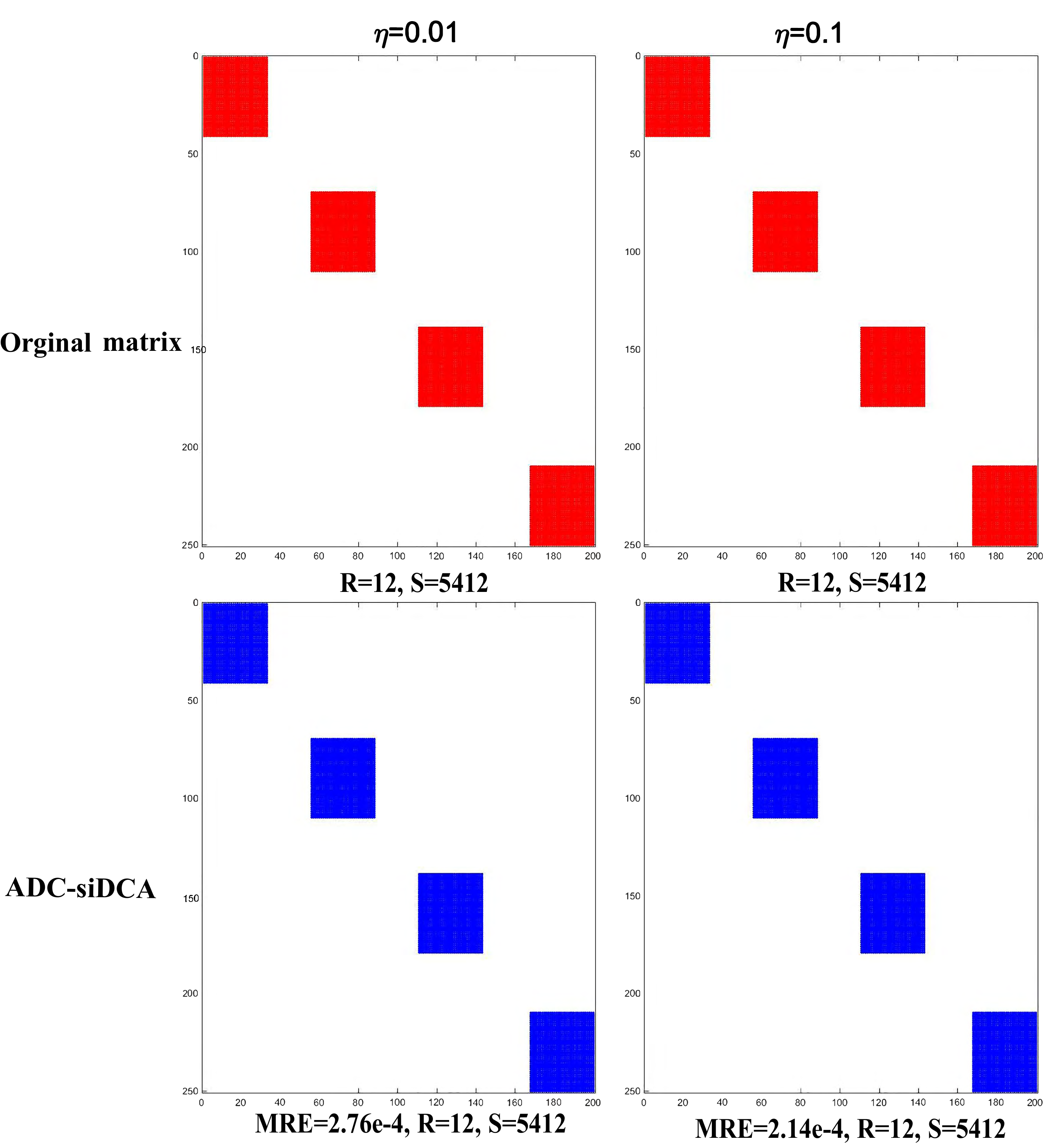}\\
  \includegraphics[width=0.45\linewidth]{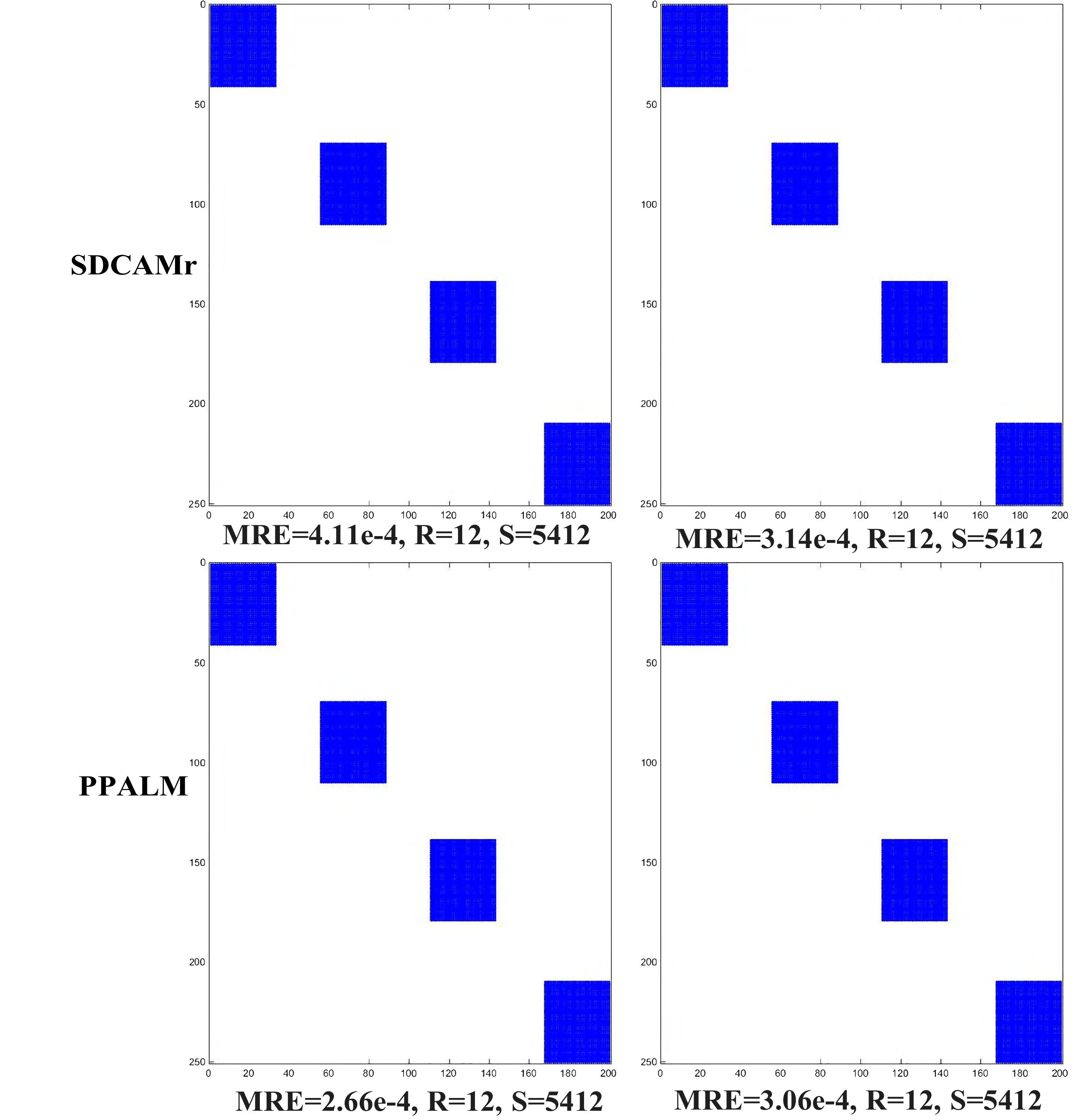}
  \caption{Comparison of sparse and low-rank matrix recovery performance of nonnegative Cliq model}\label{fig:1}
\end{figure}

\begin{figure}[H]
\centering  
  \includegraphics[width=0.45\linewidth]{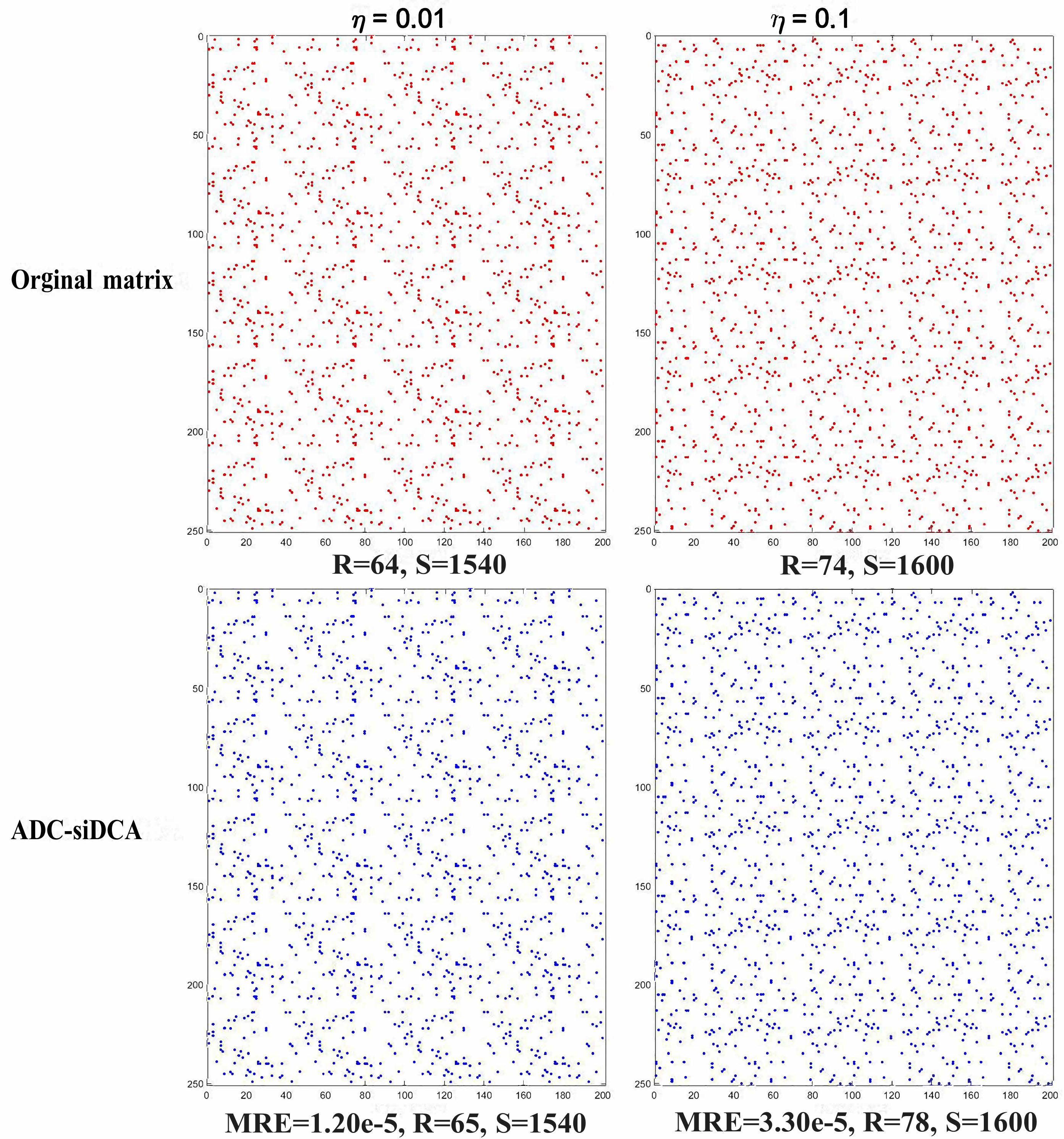}\\
  \includegraphics[width=0.45\linewidth]{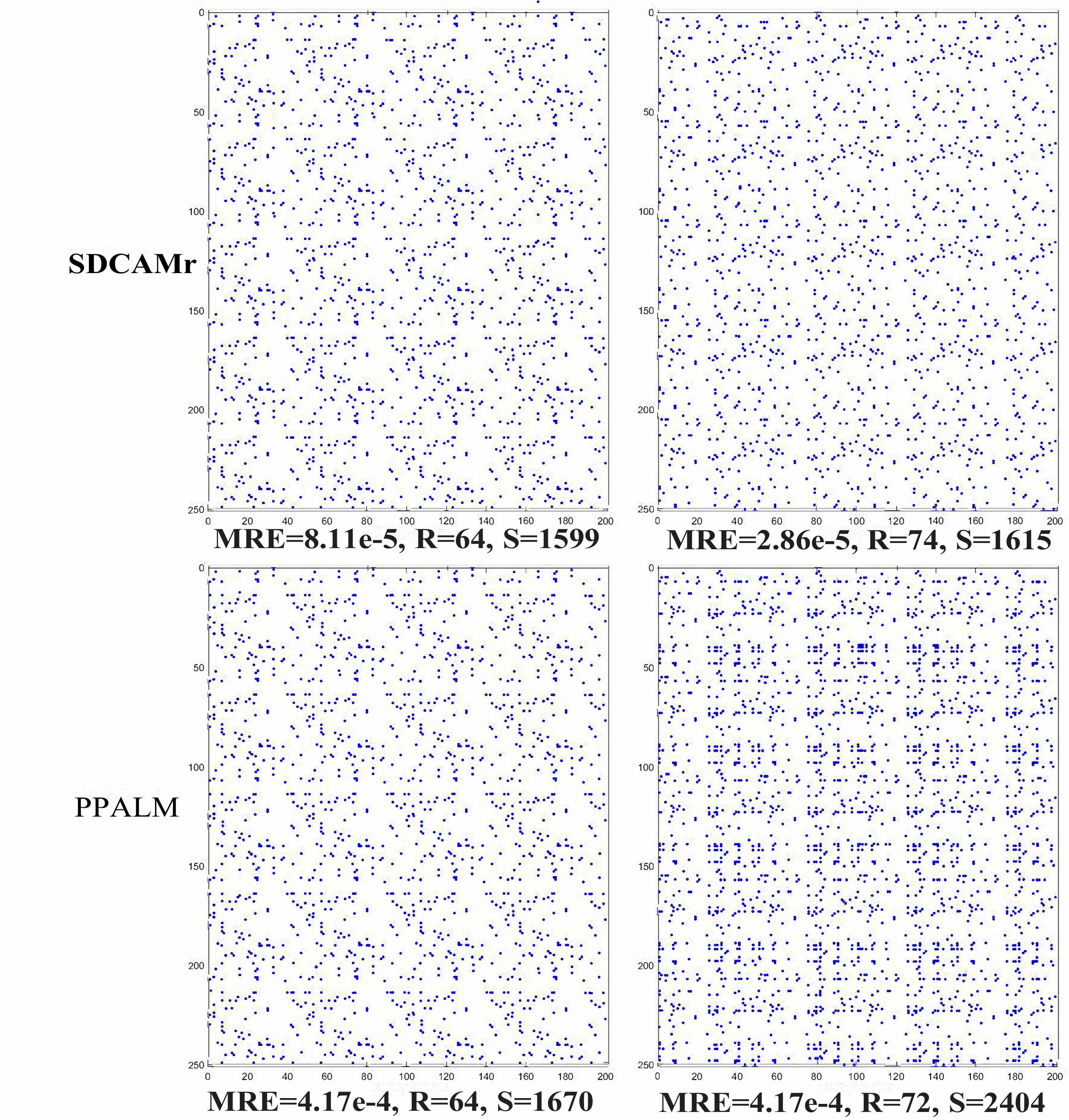}
  \caption{Comparison of sparse and low-rank matrix recovery performance of Rand1 model}\label{fig:2}
\end{figure}

\begin{figure}[H]
\centering  
  \includegraphics[width=0.45\linewidth]{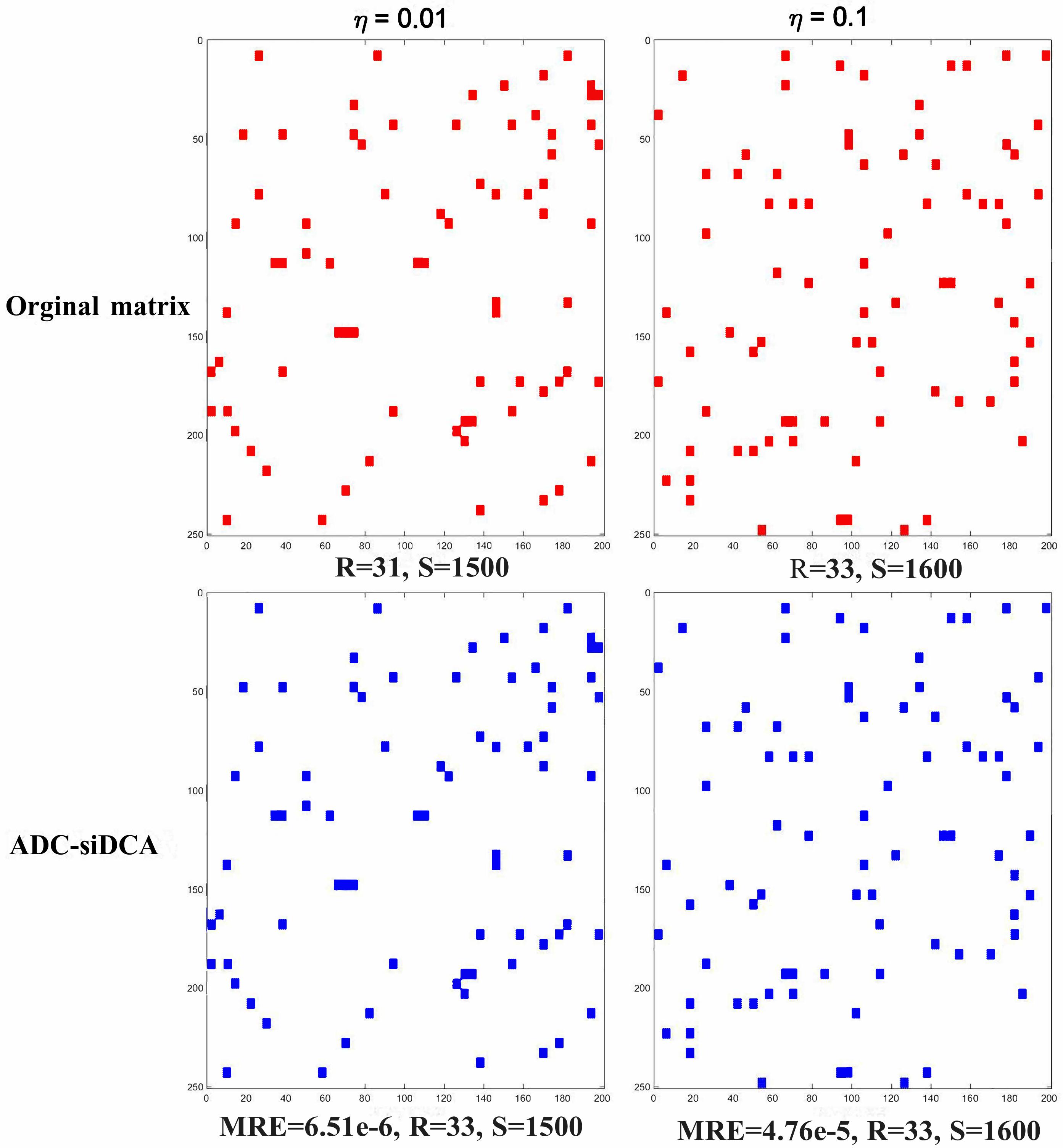}\\
  \includegraphics[width=0.45\linewidth]{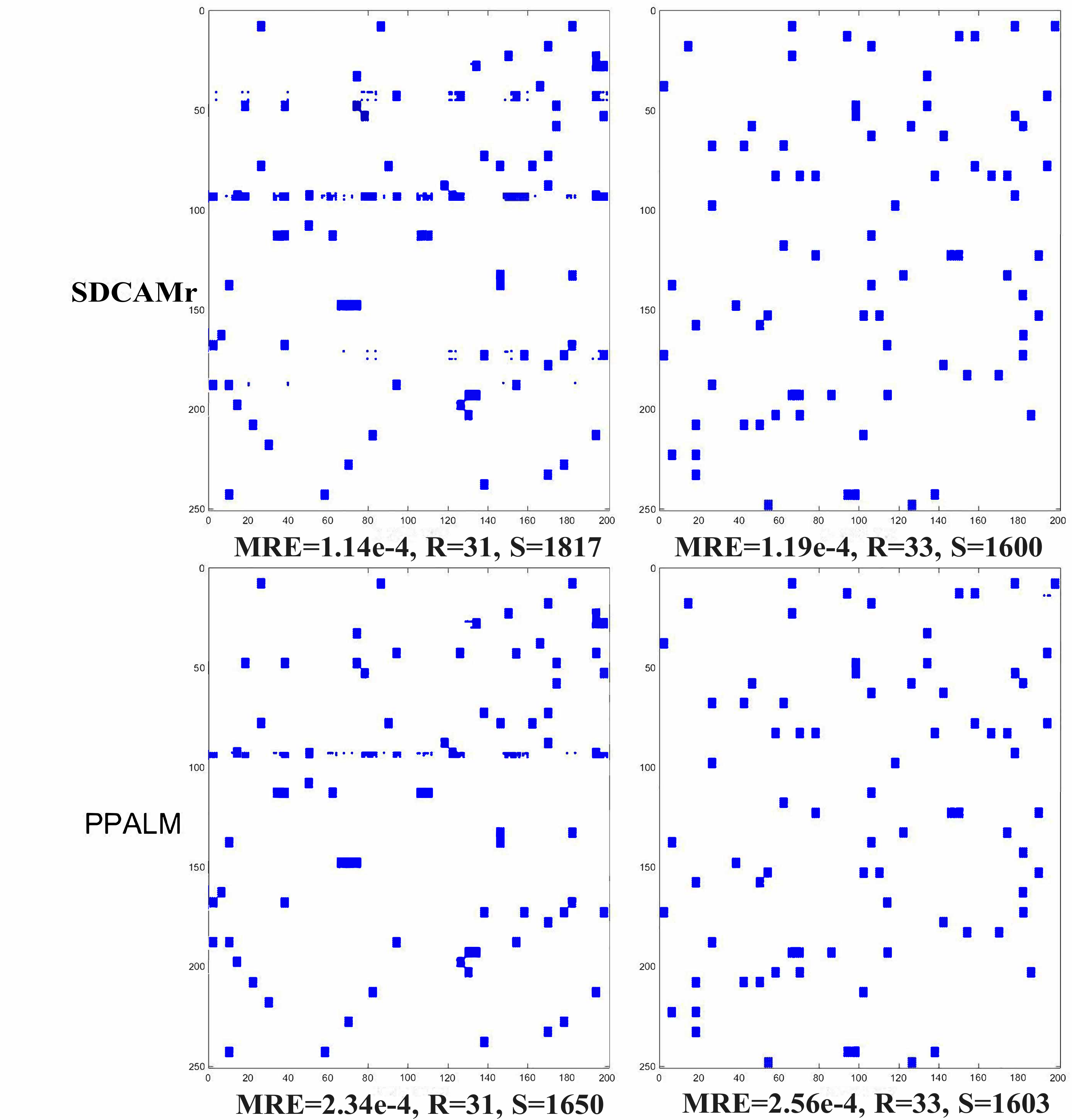}
  \caption{Comparison of sparse and low-rank matrix recovery performance of Rand2 model}\label{fig:3}
\end{figure}

In order to intuitively demonstrate the effectiveness of ADC-siDCA for solving sparse and low-rank matrix recovery problem when $\mathbb{U}_+=\mathbb{R}_+^{m\times n}$, we present the recovery results of these three methods for Cliq model, Rand1 model and Rand2 model  in Figure \ref{fig:1}, Figure \ref{fig:2} and Figure \ref{fig:3}, respectively. As shown in Figure \ref{fig:1}- Figure \ref{fig:3}, we can see that the distribution, rank and sparsity of the matrices obtained by these three methods are very close to the original matrices, which indicates that these three methods are effective for solving these three models. For complex models, such as Rand1 model and Rand2 model, compared with the other methods, the matrices recovered by ADC-siDCA satisfy rank constraint and cardinality constraint better.\par

\begin{figure}[H]
\centering 
  \includegraphics[width=0.65\linewidth]{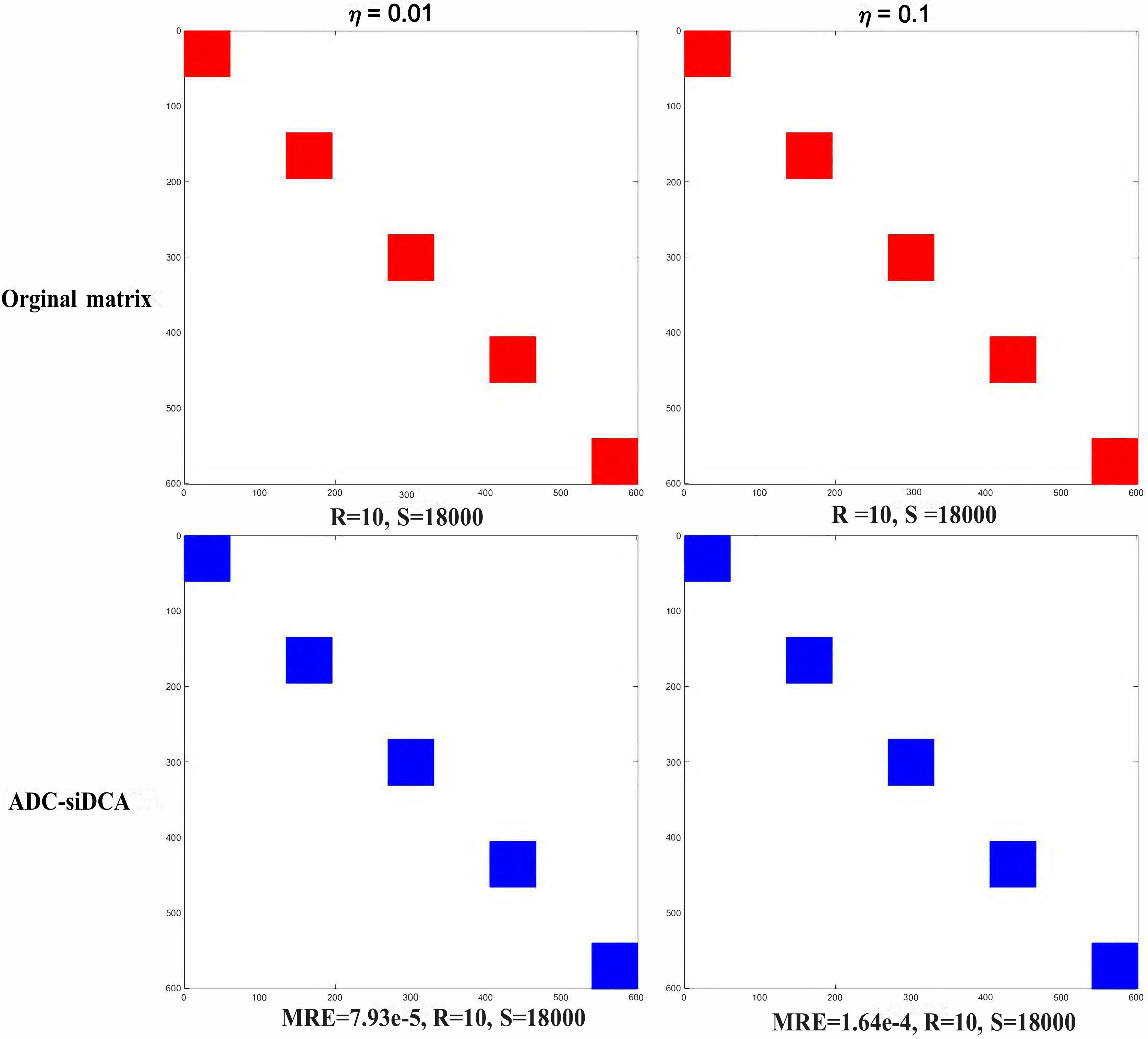}\\
  \includegraphics[width=0.65\linewidth]{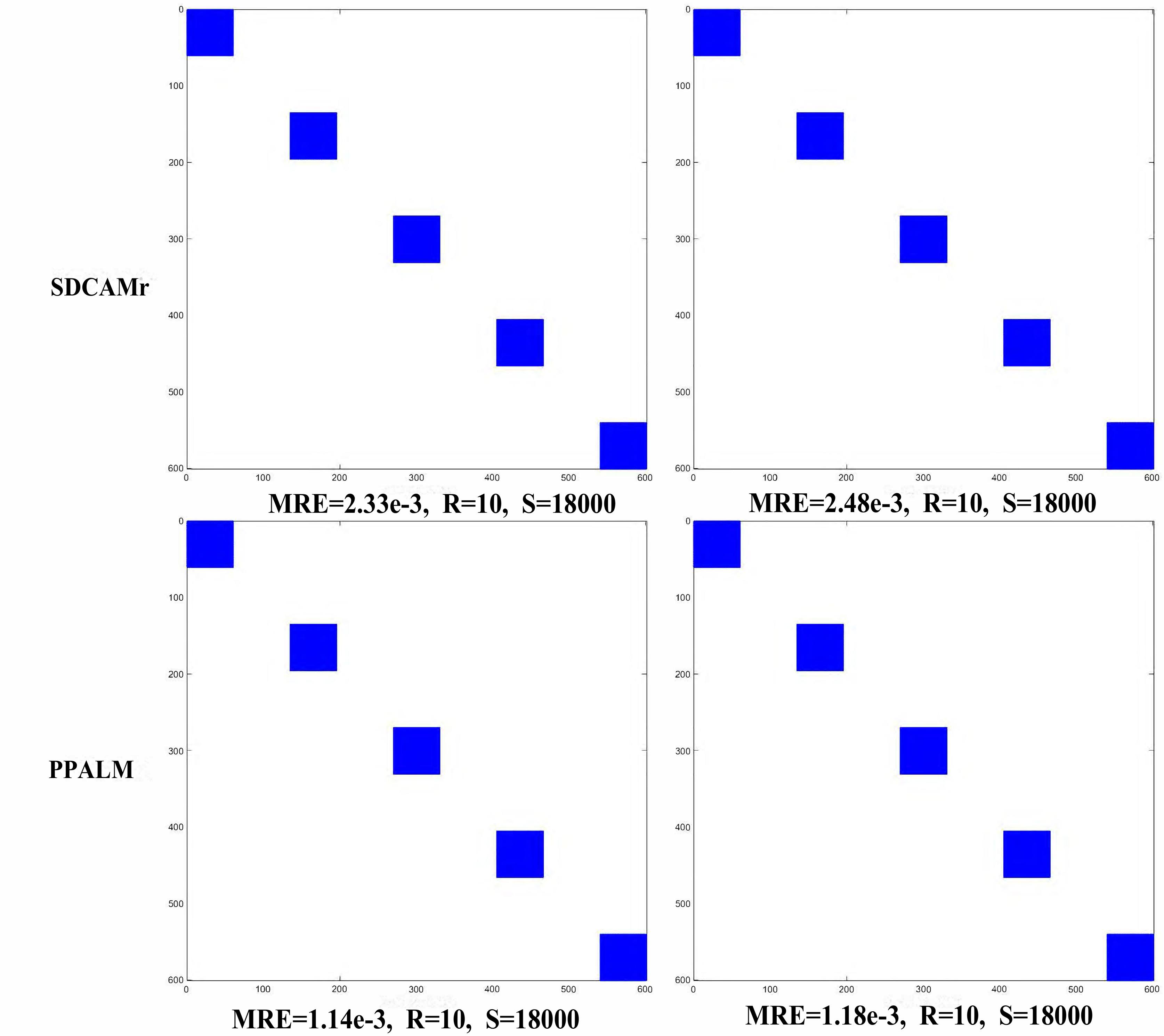}
  \caption{Comparison of sparse and low-rank recovery performance of positive semidefinite clique model}\label{fig:4}
\end{figure}

\begin{figure}[H]
\centering 
  \includegraphics[width=0.65\linewidth]{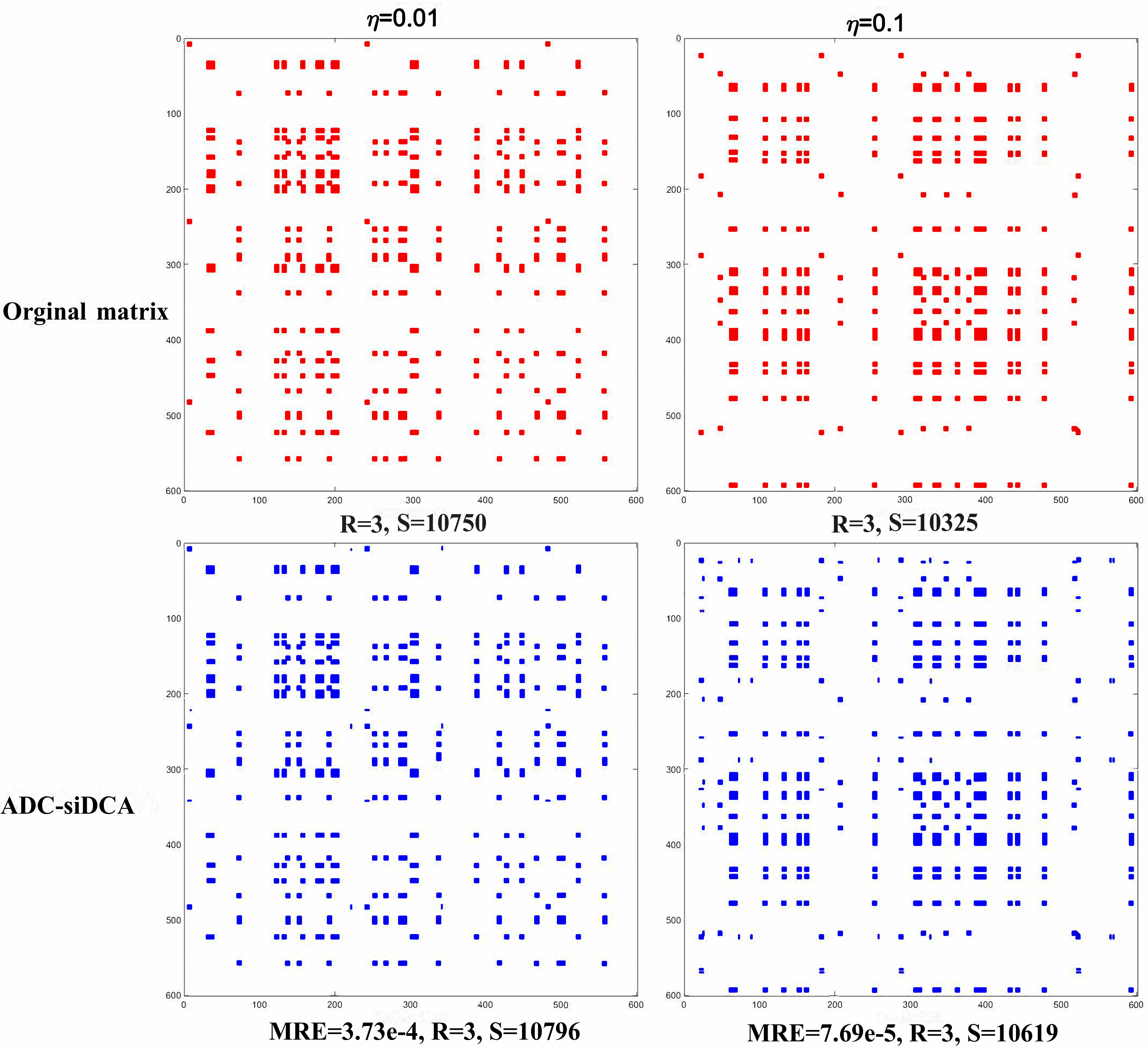}\\
  \includegraphics[width=0.65\linewidth]{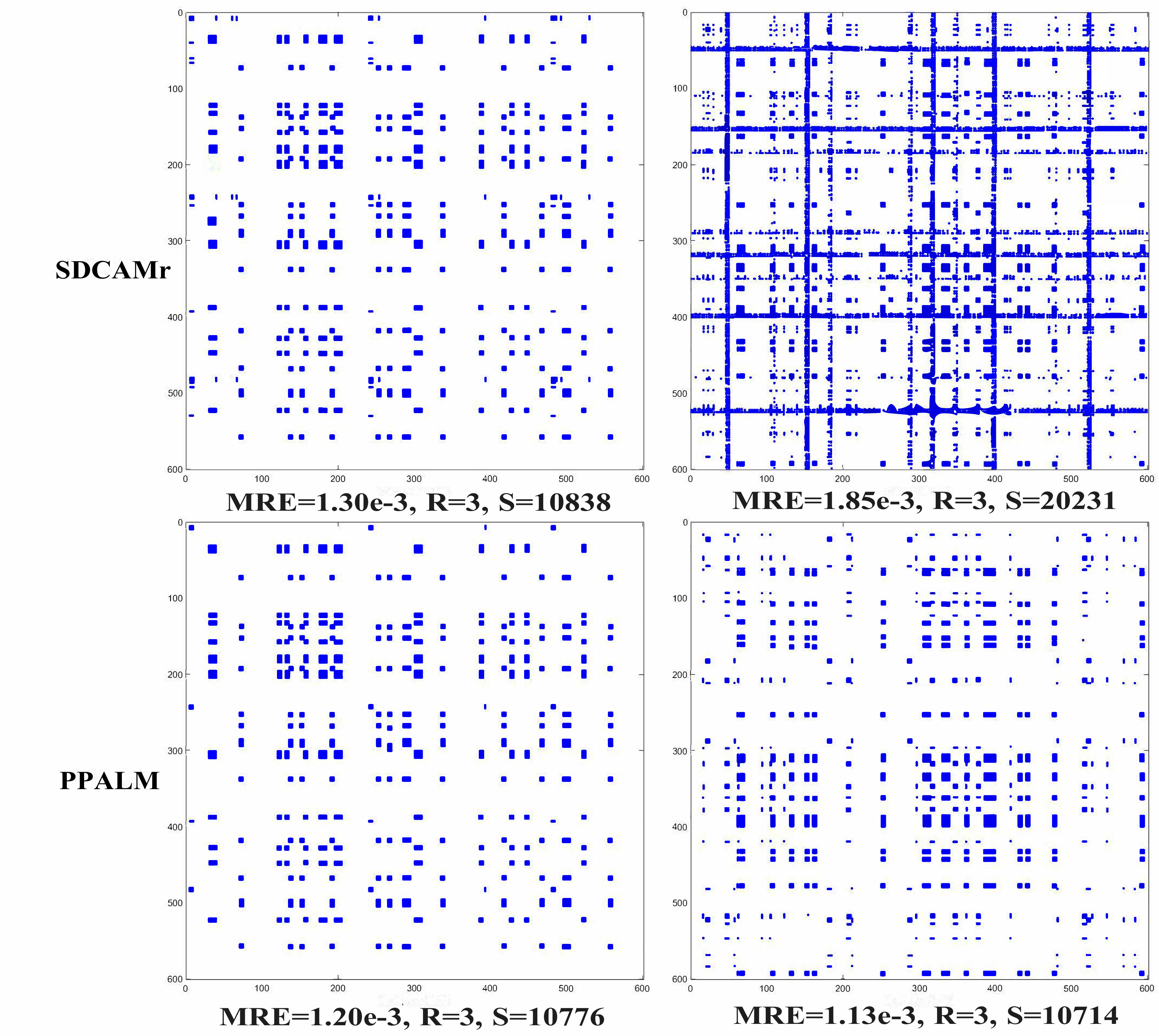}
  \caption{Comparison of sparse and low-rank recovery performance of positive semidefinite random model}\label{fig:5}
\end{figure}

\begin{figure}[H]
\centering 
  \includegraphics[width=0.65\linewidth]{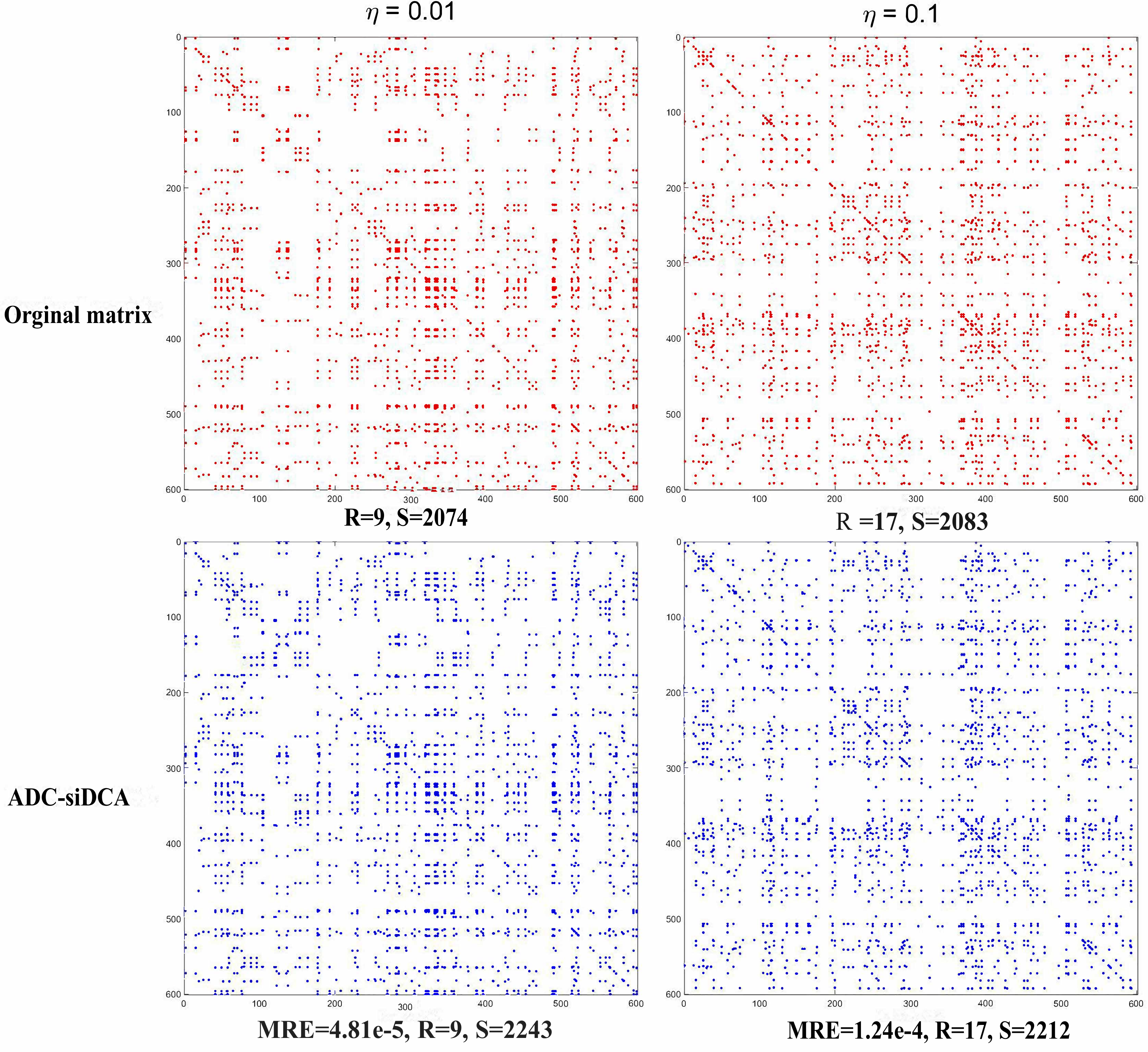}\\
  \includegraphics[width=0.65\linewidth]{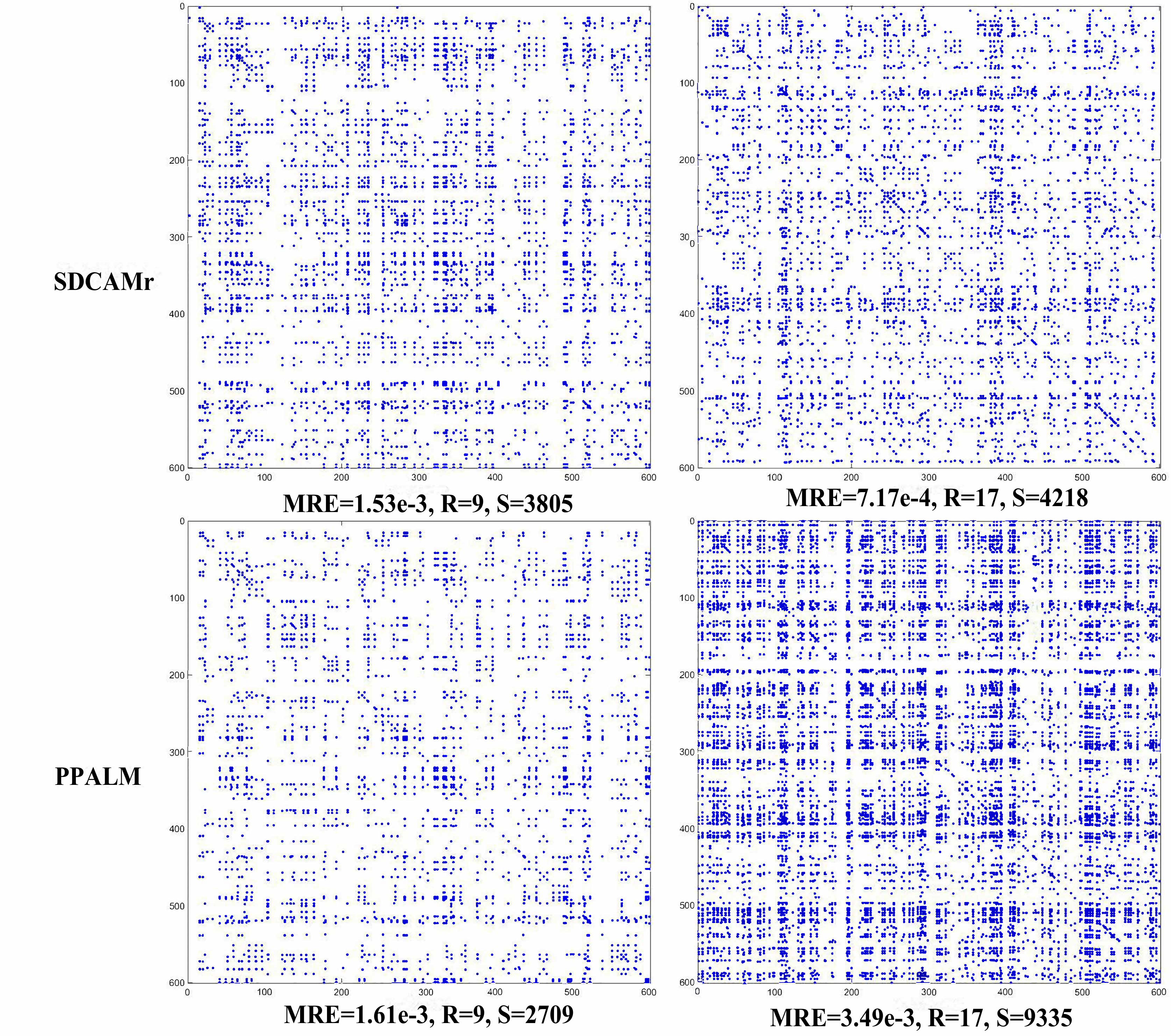}
  \caption{Comparison of sparse and low-rank recovery performance of sparsity random positive semidefinite model}\label{fig:6}
\end{figure}

In order to more intuitively demonstrate the effectiveness of ADC-siDCA for solving sparse and low-rank matrix recovery when $\mathbb{U}_+ = \mathbb{S}_+^n$, we present the recovery results of Cliq model, Rand model and Spr model of positive semidefinite matrix by the three methods in Figure \ref{fig:4}, Figure \ref{fig:5} and Figure \ref{fig:6}, respectively. As shown in Figure \ref{fig:4} - Figure \ref{fig:6}, we can see that the matrix recovered by ADC-siDCA is highly consistent with the distribution, rank and sparsity of the original matrix. In addition, for Cliq model, the matrices obtained by these three methods is highly consistent with the original matrices in terms of distribution, rank and sparsity. However, for more complex models, such as Rand model and Spr model, the recovery performance of SDCAM$_r$ and PPALM method is not satisfactory.\par



\section*{Funding}
The National Natural Science Foundation of China (Grand No. 11971092); The Fundamental Research Funds for the Central Universities (Grand No. DUT20RC(3)079)

\bibliographystyle{abbrvnat}
\bibliography{ref} 

\end{document}